\documentclass[A4paper, 12pt]{article}
\usepackage[utf8]{inputenc}
\usepackage[OT2, T1]{fontenc}
\usepackage{cite}
\usepackage{amsmath}
\usepackage{amsfonts}
\usepackage{amssymb}
\usepackage{amsxtra}
\usepackage{mathrsfs}
\usepackage{tensor}
\usepackage[thmmarks, amsmath, thref, amsthm]{ntheorem}
\usepackage[nottoc]{tocbibind}
\usepackage[pagebackref=true]{hyperref}

\usepackage{graphicx}
\usepackage[all]{xy}
\usepackage[portrait, top=2cm, bottom=2cm, left=2cm, right=2cm]{geometry}
\usepackage{enumerate}

\numberwithin{equation}{section}

\theoremstyle{plain}
\newtheorem{thm}{Theorem}[section]
\newtheorem{lemm}[thm]{Lemma}
\newtheorem{prop}[thm]{Proposition}
\newtheorem{coro}[thm]{Corollary}

\newtheorem{customthm}{Theorem}

\newtheorem{customconj}{Conjecture}

\theoremstyle{definition}
\newtheorem{defn}[thm]{Definition}
\newtheorem{remk}[thm]{Remarks}

\newcommand{\ol}[1]{\overline{#1}}
\newcommand{\wh}[1]{\widehat{#1}}
\newcommand{\pair}[1]{\left<#1\right>}

\newcommand{\ab}{\operatorname{ab}}
\newcommand{\Aut}{\operatorname{Aut}}
\newcommand{\BM}{\operatorname{BM}}
\newcommand{\CH}{\operatorname{CH}}
\newcommand{\cores}{\operatorname{cor}}
\newcommand{\cts}{\operatorname{cts}}
\newcommand{\BMHP}{\mathbf{(BMHP)}}
\newcommand{\BMWA}{\mathbf{(BMWA)}}
\newcommand{\Br}{\operatorname{Br}}
\newcommand{\Coker}{\operatorname{Coker}}
\newcommand{\dtype}{\mathbf{dtype}}
\newcommand{\DType}{\operatorname{DType}}
\newcommand{\et}{\textnormal{\'et}}
\newcommand{\Ext}{\operatorname{Ext}}
\newcommand{\Frac}{\operatorname{Frac}}
\newcommand{\Gal}{\operatorname{Gal}}
\newcommand{\Hom}{\operatorname{Hom}}
\newcommand{\id}{\operatorname{id}}
\newcommand{\Img}{\operatorname{Im}}
\newcommand{\iExt}{\mathscr{E}\kern -.5pt xt}
\newcommand{\iHom}{\mathscr{H}\kern -.5pt om}
\newcommand{\Ker}{\operatorname{Ker}}
\newcommand{\lien}{\mathbf{lien}}
\newcommand{\loc}{\operatorname{loc}}
\newcommand{\nr}{\operatorname{nr}}
\newcommand{\Pic}{\operatorname{Pic}}
\newcommand{\PT}{\operatorname{PT}}
\newcommand{\SAut}{\operatorname{SAut}}
\newcommand{\sconnect}{\operatorname{sc}}
\newcommand{\Spec}{\operatorname{Spec}}
\newcommand{\ssimple}{\operatorname{ss}}
\newcommand{\rec}{\operatorname{rec}}
\newcommand{\red}{\operatorname{red}}
\newcommand{\Res}{\operatorname{Res}}
\newcommand{\res}{\operatorname{res}}
\newcommand{\tor}{\operatorname{tor}}
\newcommand{\Tot}{\operatorname{Tot}}
\newcommand{\uni}{\operatorname{u}}
\newcommand{\Unit}{\operatorname{U}}
\newcommand{\UPic}{\operatorname{UPic}}

\newcommand{\Abb}{\mathbb{A}}
\newcommand{\Cbb}{\mathbb{C}}
\newcommand{\Gbb}{\mathbb{G}}
\newcommand{\Hbb}{\mathbb{H}}
\newcommand{\Lbb}{\mathbb{L}}
\newcommand{\Pbb}{\mathbb{P}}
\newcommand{\Qbb}{\mathbb{Q}}
\newcommand{\Rbb}{\mathbb{R}}
\newcommand{\Zbb}{\mathbb{Z}}

\newcommand{\Ccal}{\mathcal{C}}
\newcommand{\Dcal}{\mathcal{D}}
\newcommand{\Ical}{\mathcal{I}}
\newcommand{\Ocal}{\mathcal{O}}

\newcommand{\f}{\mathrm{f}}
\renewcommand{\H}{\mathrm{H}}
\newcommand{\N}{\mathrm{N}}
\newcommand{\Z}{\mathrm{Z}}

\newcommand{\Hscr}{\mathscr{H}}
\newcommand{\Tscr}{\mathscr{T}}
\newcommand{\Uscr}{\mathscr{U}}

\DeclareSymbolFont{cyrletters}{OT2}{wncyr}{m}{n}
\DeclareMathSymbol{\Be}{\mathalpha}{cyrletters}{"42}
\DeclareMathSymbol{\Sha}{\mathalpha}{cyrletters}{"58}

\newcommand{\del}{\partial}

\title{On the descent conjecture for rational points and zero-cycles}
\author{Nguy$\tilde{\hat{\text{e}}}$n M$\d{\text{a}}$nh Linh}
\date{\today}

\begin{document}
\maketitle

\begin{abstract}
    The descent method is one of the approaches to study the Brauer--Manin obstruction to the local--global principle and to weak approximation on varieties over number fields, by reducing the problem to ``descent varieties''. In recent lecture notes by Wittenberg, he formulated a ``descent conjecture'' for torsors under linear algebraic groups. The present article gives a proof of this conjecture in the case of connected groups, generalizing the toric case from the previous work of Harpaz--Wittenberg. As an application, we deduce directly from Sansuc's work the theorem of Borovoi for homogeneous spaces of connected linear algebraic groups with connected stabilizers. We are also able to reduce the general case to the case of finite (\'etale) torsors. When the set of rational points is replaced by the Chow group of zero-cycles, an analogue of the above conjecture for arbitrary linear algebraic groups is proved.
\end{abstract}


{\em 2020 Mathematics subject classifications (MSC).} 11R34, 14G05, 14G12.

\tableofcontents

\section{Introduction} \label{sec:Intro}

\subsection{Brauer--Manin obstruction and the descent conjecture} \label{subsec:Points}

Let $X$ be a smooth variety defined over a number field $k$. We say that $X$ satisfies the {\em Hasse principle} (or the {\em local--global principle}) if either $X(k_v) = \varnothing$ for some completion $k_v$ of $k$ or $X(k) \neq \varnothing$. It was Manin who introduced a general method to study the failure of the Hasse principle on $X$, as follows. Assuming $X$ geometrically integral, we associate to $X$ the cohomological Brauer group $\Br(X):= \H^2_{\et}(X,\Gbb_m)$ (which generalizes the Brauer group $\Br(k) = \H^2(k,\bar{k}^\times)$ of the field $k$), and consider the unramified subgroup $\Br_{\nr}(X) \subseteq \Br(X)$ (see Subsection \ref{subsec:Notation} below). One then defines the {\em Brauer--Manin pairing} \cite{Manin1971Brauer}
    \begin{equation} \label{eq:BrauerManinPairingNr} \tag{BM${}_{\nr}$}
        \pair{-,-}_{\BM}: X(k_\Omega) \times \Br_{\nr}(X) \to \Qbb/\Zbb, \quad ((x_v)_v,\alpha) \mapsto \sum_v \alpha(x_v),
    \end{equation}
where $X(k_\Omega)$ is the product of all the $X(k_v)$'s, and where $\Br(k_v) \hookrightarrow \Qbb/\Zbb$ {\em via} the invariant map from local class field theory. By the global reciprocity law (Albert--Brauer--Hasse--Noether), the pairing \eqref{eq:BrauerManinPairingNr} vanishes on the image of $\Br(k)$ in $\Br_{\nr}(X)$ and on the diagonal image of $X(k)$ in $X(k_\Omega)$. If we denote by $X(k_\Omega)^{\Br_{\nr}(X)}$ the subset of families of local points of $X$ orthogonal to $\Br_{\nr}(X)$, then its vacuity is an obstruction to the existence of $k$-points on $X$, the so-called {\em Manin obstruction}. If either $X(k_\Omega)^{\Br_{\nr}(X)} = \varnothing$ or $X(k) \neq \varnothing$ (that is, the Manin obstruction to the Hasse principle for $X$ is the only one), we shall say that $X$ {\em satisfies the property $\BMHP$}.

We are also interested in the question of weak approximation, namely, the density of $X(k)$ in $X(k_\Omega)$ (equipped with the product of $v$-adic topologies). By continuity of \eqref{eq:BrauerManinPairingNr}, the closure $\ol{X(k)}$ is contained in the (closed) set $X(k_\Omega)^{\Br_{\nr}(X)}$. Thus, a family of local points on $X$ not orthogonal to $\Br_{\nr}(X)$ cannot be approximated by $k$-rational points. This is known as the {\em Brauer--Manin obstruction} to weak approximation. If $\ol{X(k)} = X(k_\Omega)^{\Br_{\nr}(X)}$ (that is, the Brauer--Manin obstruction to weak approximation for $X$ is the only one), then we say that $X$ {\it satisfies the property $\BMWA$}. Note that we included here the case where both sides are empty (hence, the property $\BMWA$ implies the property $\BMHP$). The properties $\BMHP$ and $\BMWA$ are stable birational invariants, thanks to the stable birational invariance of the unramified Brauer group, the theorem of Nishimura--Lang \cite{Nishimura1955Rational}, and the Artin--Whaples approximation lemma.

\begin{defn} [Campana--Koll\'ar--Miyaoka--Mori] \label{defn:RationallyConnected}
    A smooth variety $X$ over a field $k$ is said to be {\em rationally connected} if, for every algebraically closed field $K$ containing $k$, two general points $x_0$, $x_1 \in X(K)$ can be connected by a $K$-rational map $\Abb^1_K \dashrightarrow X \times_k K$.
\end{defn}

Rational connectedness is a stable birational invariant. Examples of rationally connected varieties include (geometrically) unirational varieties. For example, in characteristic $0$, this is the case for homogeneous spaces of connected linear algebraic groups (indeed, the celebrated theorem of Chevalley asserts that every such groups are geometrically rational, and even unirational over their field of definition \cite{Chevalley1954Group}). An important conjecture concerning arithmetics of rationally connected varieties over number fields is the following.

\begin{customconj} [Colliot-Th\'el\`ene] \label{conj:RationallyConnected}
    Every smooth rationally connected variety over a number field satisfies the properties $\BMHP$ and $\BMWA$.
\end{customconj}

One of the reasons why Conjecture \ref{conj:RationallyConnected} is so strong is that its truth would imply that any finite group is a Galois group over $k$. This is the inverse Galois problem over number fields, long known to be open. The reader is invited to take a look at the recent notes \cite{Wittenberg2024ParkCity} by Wittenberg. One of the techniques to attack Conjecture \ref{conj:RationallyConnected} is the {\em descent theory}. Namely, let $f: Y \to X$ be a {\em torsor} under a linear algebraic group $G$ defined over a number field $k$, that is, $f$ is an fppf morphism between smooth $k$-varieties, where $Y$ equipped with a right action $Y \times_k G \to Y$ of $G$ such that $f$ is $G$-equivariant and $G(\bar{k})$ acts simply transitively on the fibre $f^{-1}(x)$ over each geometric point $x \in X(\bar{k})$. For each Galois cocycle $\sigma: \Gal(\bar{k}/k) \to G(\bar{k})$, we have a {\em twisted torsor} $f^{\sigma}: Y^\sigma \to X$ under an inner $k$-form $G^{\sigma}$ of $G$. Then, we have the following partition of $X(k)$ (see {\em e.g.} \cite[p. 22]{Skorobogatov2001Torsors}):
    \begin{equation*}
        X(k) = \bigsqcup_{[\sigma] \in \H^1(k,G)} f^{\sigma}(Y^{\sigma}(k)).
    \end{equation*}
The idea of the descent theory is to establish an equality of the form
    \begin{equation*}
        X(k_\Omega)^{\Br_{\nr}(X)} = \ol{\bigcup_{[\sigma] \in \H^1(k,G)} f^{\sigma}(Y^{\sigma}(k_\Omega)^{{\Br_{\nr} (Y^{\sigma})}})}
    \end{equation*}
in the topological space $X(k_\Omega)$. In particular, if the ``descent varieties'' $Y^{\sigma}$ satisfy $\BMHP$ (resp. $\BMWA$), then it is also the case for $X$. In the very same notes, Wittenberg formulated the following conjecture, which we shall call the ``descent conjecture'' \cite[Conjecture 3.7.4]{Wittenberg2024ParkCity}.

\begin{customconj} \label{conj:Descent} 
    With the above notations and assumptions, suppose in addition that $Y$ (hence also $X$) is rationally connected. If the varieties $Y^{\sigma}$ satisfy the property $\BMHP$ (resp. $\BMWA$), then $X$ itself satisfies the property $\BMHP$ (resp. $\BMWA$).
\end{customconj}

Descent theory was originally developed in the foundational work \cite{CTS1987Descente} of Colliot-Th\'el\`ene and Sansuc for torsors under a torus $T$, over a variety $X$ with $\bar{k}[X]^\times = \bar{k}^\times$. From its birth, this approach has given many nice results on the Brauer--Manin obstruction for (geometrically) rational varieties, notably for Ch\^atelet surfaces as in the works \cite{CTSSD1987I, CTSSD1987II} by Colliot-Th\'el\`ene, Sansuc, and Swinnerton-Dyer. As noted by Skorbogatov (see \cite[\S{6.1}]{Skorobogatov2001Torsors}), the torus $T$ can be replaced by a group of multiplicative type $S$. The torsors under consideration ({\em i.e.}, the twists of a given torsor) are the torsors of the same ``type''. The notion of extended type, introduced later by Harari and Skorobogatov \cite{HS2013Descent}, enabled descent theory (under a group of multiplicative type) for varieties with nonconstant invertible functions. Other generalizations in this direction are found in the various works of Stoll \cite{Stoll2006Descent}, Demarche \cite{Demarche2009Etale}, Skorobogatov \cite{Skorobogatov2009Descent}, Poonen \cite{Poonen2010Points}, Wei \cite{Wei2016Open}, and Cao--Demarche--Xu \cite{CDX2019Compare}..., where there is notably the notion of {\em \'etale Brauer--Manin obstruction} for torsors under a finite (\'etale) group scheme.

Conjecture \ref{conj:Descent} was established by Harpaz and Wittenberg for torsors under a torus \cite{HW2020Galois} or a finite supersolvable group (scheme) \cite{HW2024Supersolvable}. In the present article, we prove this conjecture for torsors under an arbitrary connected linear algebraic group. More precisely, we shall prove

\begin{customthm} [Theorem \ref{thm:DescentConnected}] \label{customthm:DescentConnected}
    Let $k$ be a number field, $X$ a smooth rationally connected $k$-variety, and $G$ a connected linear algebraic $k$-group. Let $Y \to X$ be a torsor under $G$, then
 	\begin{equation*}
 		X(k_\Omega)^{\Br_{\nr}(X)} = \ol{\bigcup_{[\sigma] \in \H^1(k,G)} f^{\sigma}(Y^{\sigma}(k_\Omega)^{\Br_{\nr}(Y^{\sigma})})}
 	\end{equation*}
    in the topological space $X(k_\Omega)$. In particular, if the varieties $Y^{\sigma}$ satisfy $\BMHP$ (resp. $\BMWA$) for all $[\sigma] \in \H^1(k,G)$, then it is also the case for $X$.
\end{customthm}

The proof of Theorem \ref{customthm:DescentConnected} relies on the technique of Harpaz--Wittberg in \cite{HW2020Galois} and Borovoi's general machinery of {\em abelianization of nonabelian cohomology} \cite{Borovoi1998Reductive}. This latter allows one to pass from the nonabelian Galois cohomology set $\H^1(k,G)$ to a hypercohomology group of a $2$-term complex of tori. We shall discuss this in Section \ref{sec:Preliminary} before proving Theorem \ref{customthm:DescentConnected} in Section \ref{sec:DescentConnected}. We also prove the following result, which reduces Conjecture \ref{conj:Descent} to the case of finite (\'etale) torsors.

\begin{customthm} [Theorem \ref{thm:DescentFinite}] \label{customthm:DescentFinite}
    Let $k$ be a number field, $X$ a smooth geometrically integral $k$-variety, and $G$ a linear algebraic $k$-group. Suppose that $f: Y \to X$ is a torsor under $G$, with $Y$ rationally connected (the existence of such a torsor requires in particular that $X$ is rationally connected). Let $G^{\circ}$ denote the identity component of $G$ and $Z:=Y/G^{\circ}$, which gives a torsor $g: Z \to X$ under the finite $k$-group $G/G^{\circ}$. Then
 	\begin{equation*}
            \ol{\bigcup_{[\sigma] \in \H^1(k,G)} f^{\sigma} (Y^{\sigma}(k_\Omega)^{\Br_{\nr}(Y^{\sigma})})} = \ol{\bigcup_{[\tau] \in \H^1(k,G/G^{\circ})} g^{\tau}(Z^{\tau}(k_\Omega)^{\Br_{\nr}(Z^{\tau})})}
	\end{equation*}
    in $X(k_\Omega)$. In particular, Conjecture \ref{conj:Descent} holds for the torsor $Y \to X$ if it holds for $Z \to X$.
\end{customthm}

What we shall show is in fact slightly stronger. Theorems \ref{customthm:DescentConnected} and \ref{customthm:DescentFinite} hold true even when we do not have a torsor {\em a priori}, but only a {\em descent type} (see Theorems \ref{thm:DescentConnectedType} and \ref{thm:DescentFiniteType}). This is a definition introduced in the recent work \cite{Linh2025Type}, which unifies the various notions of type by Colliot-Th\'el\`ene--Sansuc, Harari--Skorobogatov, and Harpaz--Wittenberg. This new notion enables descent theory in the case where $X(k)$ is not necessarily nonempty (let us say, when $X$ is a homogeneous space of a linear algebraic group), which is useful for the study of the property $\BMHP$. 

As an application, we give a new proof for Borovoi's theorem, which states that the property $\BMHP$ and $\BMWA$ hold for homogeneous spaces of connected linear algebraic groups with connected geometric stabilizers ({\em cf.} Corollary \ref{coro:BorovoiThm}), using the theorem of Sansuc for {\em principal} homogeneous spaces.

\subsection{Analogue for zero-cycles} \label{subsec:Cycles}

In Section \ref{sec:DescentCycle}, we consider the Brauer--Manin obstruction to the Hasse principle and to weak approximation for {\em zero-cycles}. Let $X$ be a smooth {\em proper} geometrically integral variety over a number field $k$. Let $\Z_0(X)$ be the free abelian groups of zero-cycles on $X$, whose elements are (finite) formal sum $x = \sum_{P} n_P [P]$, where $P$ runs through the closed points of $X$, and $n_P \in \Zbb$. The {\em degree} of such a cocycle is $\deg_k(x) := \sum_P n_P [k(P):k]$. For any $\alpha \in \Br(X)$ and place $v$ of $k$, we may evaluate $\alpha$ at the closed points $P_v \in X_{k_v}$ to obtain $\alpha(P_v) \in \Br(k_v(P_v)) \hookrightarrow \Qbb/\Zbb$. By linearity, we obtain a local pairing
    \begin{equation*}
        \pair{-,-}_v: \Z_0(X_{k_v}) \times \Br(X) \to \Qbb/\Zbb.
    \end{equation*}
Let $\CH_0(X)$ denote the quotient of $\Z_0(X)$ by rational equivalence (see {\em e.g.} \cite[Chapter 1]{Fulton1984Intersection}). Since $X$ is proper, the above pairing factors through a pairing $\CH_0(X_{k_v}) \times \Br(X) \to \Qbb/\Zbb$. If $v$ is an Archimedean place, this latter also factors through $\CH_0'(X_{k_v}):=\Coker(\N_{\bar{k}_v/k_v})$, where $\N_{\bar{k}_v/k_v}: \CH_0(X_{\bar{k}_v}) \to \CH_0(X_{k_v})$ denotes the norm (or corestriction) map. By reassembling these local pairings, we obtain a Brauer--Manin type pairing
    \begin{equation*}
		\pair{-,-}_{\BM}: \CH_{0,\Omega}(X) \times \Br(X) \to \Qbb/\Zbb,
    \end{equation*}
where $\CH_{0,\Omega}(X):=\prod_{v \nmid \infty} \CH_0(X_v) \times \prod_{v | \infty} \CH_0'(X_v)$. The pairing $\pair{-,-}_{\BM}$ vanishes on the image of the localization (or restriction) map $\CH_0(X) \to \CH_{0,\Omega}(X)$. Hence, we obtain a complex
    \begin{equation*}
        \CH_0(X) \to \CH_{0,\Omega}(X) \to \Hom(\Br(X),\Qbb/\Zbb).
    \end{equation*}
Since the group $\Br(X)$ is torsion, we have a complex
    \begin{equation} \label{eq:ConjectureE} \tag{E}
        \varprojlim_n (\CH_0(X)/n) \to \varprojlim_n (\CH_{0,\Omega}(X)/n) \to \Hom(\Br(X),\Qbb/\Zbb).
    \end{equation}
    
\begin{customconj} [Conjecture (E)] \label{conj:E}
    The complex \eqref{eq:ConjectureE} is exact for any smooth proper geometrically integral variety $X$ over a number field $k$.
\end{customconj}

The above conjecture is due to Colliot-Th\'el\`ene, Sansuc, Kato, and Saito \cite{CTS1981Chow, KS1986Arithmetic, CT1995Chow}, and has been reformulated into this form by van Hamel \cite{vanHamel2003Cycle} and Wittenberg \cite{Wittenberg2012Fibration}. It entails two previous conjectures, namely Conjecture (E${}_0$) and Conjecture (E${}_1$). 

\begin{customconj} [Conjecture (E${}_1$)] \label{conj:E1}
    Let $X$ be a smooth proper geometrically integral variety over a number field $k$. If there exists a family of local zero-cycles of degree $1$ on $X$ that is orthogonal to $\Br(X)$ relative to the pairing  $\pair{-,-}_{\BM}$, then there exists a zero-cycle of degree $1$ on $X$.
\end{customconj}

The property in the statement of Conjecture \ref{conj:E1} can be seen as an analogue of $\BMHP$.

\begin{customconj} [Conjecture (E${}_0$)] \label{conj:E0}
    Let $X$ be a smooth proper geometrically integral variety over a number field $k$. Let $A_0(X) \subseteq \CH_0(X)$ be the subgroup of rationally equivalence classes of zero-cycles of degree $0$. Then we have an exact sequence
        \begin{equation} \label{eq:ConjectureE0} \tag{E${}_0$}
            \varprojlim_n (A_0(X)/n) \to \prod_v\left (\varprojlim_n (A_0(X_{k_v})/n) \right) \to \Hom(\Br(X),\Qbb/\Zbb).
	\end{equation}
\end{customconj}
As was pointed out by Wittenberg in \cite[Remarques 1.1 (ii)]{Wittenberg2012Fibration}, we have
	\begin{equation*}
		\varprojlim_n (A_0(X_{k_v})/n) = \Coker(\N_{\bar{k}_v/k_v}: A_0(X_{k_v}) \to A_0(X_{\bar{k}_v}))
	\end{equation*}
for Archimedean $v$. Thus, exactness of \eqref{eq:ConjectureE} does indeed imply that of \eqref{eq:ConjectureE0}. Liang has shown that the exactness of  \eqref{eq:ConjectureE0} and the existence of a global zero-cycle of degree $1$ imply some sort Brauer--Manin obstruction to weak approximation for zero-cycles \cite[Proposition 2.2.1]{Liang2013Arithmetic} (as soon as a suitable topology on the groups $\Z_0(X_{k_v})$ is defined, see Subsection \ref{subsec:DescentCycleBM} below). Thus, the property in Conjecture \ref{conj:E0} is an analogue of $\BMWA$ (assuming the existence of a $k$-rational point).

In the recent work \cite{BB2024Descent}, Balestrieri and Berg introduced a formalism for the descent theory for zero-cycles. This allows us to formulate an analogue of Conjecture \ref{conj:Descent} in this context. It turns out that, under the rationally connected assumption, this version of the descent conjecture is true for torsors under an arbitrary linear algebraic group.

\begin{customthm} [Theorem \ref{thm:DescentCycle}] \label{customthm:DescentCycle}
    Let $k$ be a number field, $X$ a smooth geometrically integral $k$-variety, and $G$ a linear algebraic $k$-group. Suppose that $Y \to X$ is a torsor under $G$, with $Y$ rationally connected (in particular, $X$ is rationally connected). If Conjecture \ref{conj:E} holds for the smooth compactifications of the twisted torsors $Y^{\sigma}_K \to X_K$, for all finite extensions $K/k$ and all $[\sigma] \in \H^1(K,G)$, then it also holds for the smooth compactifications of $X$.
\end{customthm}

The proof of Theorem \ref{customthm:DescentCycle} combines Liang's method \cite{Liang2013Arithmetic}, the fibration theorems of Harpaz and Wittenberg \cite{HW2016Fibration}, and Theorem \ref{customthm:DescentConnected} on ``connected descent''.

\subsection{Notations and conventions} \label{subsec:Notation} 

The following notations shall be deployed throughout the article.

{\bf Abelian groups.} When $A$ is a topological abelian group which is either discrete torsion, profinite, or of finite exponent, we denote $A^D:=\Hom_{\cts}(A,\Qbb/\Zbb)$ its Pontrjagin dual.

{\bf Fields.} When $k$ is a field, we use $\bar{k}$ to denote a fixed {\em separable} closure of $k$ and $\Gamma_k:=\Gal(\bar{k}/k)$ to denote its absolute Galois group. By convention, a $k$-{\em variety} is a separated $k$-scheme of finite type. We usually denote $\bar{X}:=X \times_k \bar{k}$, but $\bar{X}$ does not always mean the base change to $\bar{k}$ of some $k$-variety. A {\em smooth compactification} of $X$ is a smooth proper variety containing $X$ as a dense open subset. In characteristic $0$, such a compactification (assuming $X$ smooth) always exists by resolution of singularities \cite{Hironaka1964Resolution}.

A {\em local field} is, by definition, the field $\Rbb$, $\Cbb$, or a field $k$ which is complete with respect to a discrete valuation, with finite residue field. If $k$ is a number field, $\Omega_k$ (resp. $\Omega_k^{\f}$, resp. $\Omega_k^\infty$) denotes the set of places (resp. finite places, resp. Archimedean places) of $k$. For $v \in \Omega_k$, $k_v$ denotes the $v$-adic completion of $k$, and $\Ocal_v$ denotes the ring of integers of $k_v$ if $v \in \Omega_k^{\f}$. If $X$ is a $k$-variety, let $X_v:=X \times_k k_v$, and write $\loc_v: \H^i(X,-) \to \H^i(X_v,-)$ for the localization maps in cohomology. Equip the set $X(k_v)$ with the natural topology induced by the topology on $k_v$ (see {\em e.g.} \cite[Proposition 3.1]{Conrad2012Adelic}), and the set $X(k_\Omega):=\prod_{v \in \Omega_k} X(k_v)$ with the product topology. On the subset $X(\Abb_k) \subseteq X(k_\Omega)$ of {\em adelic points} of $X$, we always consider the topology induced by the product topology\footnote{There is also the {\em adelic} topology, related to strong approximation, but this does not concern us.} on $X(k_\Omega)$.

{\bf Cohomology.} Unless state otherwise, all (hyper-)cohomology groups will be \'etale or Galois. If $X$ is a scheme, $\Dcal^+(X)$ denotes the derived category of bounded below complexes of sheaves of abelian groups on $X_{\et}$. If $k$ is a field, $\Dcal^+(k):=\Dcal^+(\Spec(k))$ is equivalent to the bounded below derived category of $\Gamma_k$-modules.

Let $k$ be a number field, and let $C$ be a bounded below complex of either $\Gamma_k$-modules which are finitely generated (as abelian groups) or commutative algebraic $k$-groups. Let $\Ccal$ be a complex of commutative group schemes over $\Spec(\Ocal_{k,\Sigma})$ extending $C$, where $\Sigma \subseteq \Omega_k$ is some finite set containing $\Omega_k^{\infty}$. For $i \ge 0$, the {\em restricted product} $\Pbb^i(k,C)$ is the subset of $\prod_{v \in \Omega_k} \H^i(k_v,C)$ consisting of families $(c_v)_{v \in \Omega_k}$ for which $c_v$ comes from $\H^i(\Ocal_v,\Ccal)$ for all but finitely many $v \notin \Sigma$. We define the {\em $i$th Tate--Shafarevich group} $\Sha^i(k,C)$ to be the kernel of the localization $\H^i(k,C) \to \Pbb^i(k,C)$.

{\bf Nonabelian cohomology.} Let $X$ be a scheme and let $G \to X$ be a {\em smooth} (not necessarily commutative) $X$-group scheme. A (right) {\em $X$-torsor under $G$} is, by definition, an fppf $X$-scheme $Y$ equipped with a right $X$-action $Y \times_X G \to Y$ of $G$, such that the morphism
    \begin{equation*}
		Y \times_X G \to Y \times_X Y, \quad (y,g) \mapsto (y, y \cdot g)
    \end{equation*}
is an isomorphism. When $G$ is affine, these torsors are classified by the nonabelian $\check{\text{C}}$ech cohomology (pointed) set $\H^1(X,G)$ (see \cite[Chapter III, \S 4]{Milne1980Etale} or \cite[p. 18]{Skorobogatov2001Torsors}). The class in $\H^1(X,G)$ of such an $X$-torsor $Y$ under $G$ is denoted by $[Y]$.

{\bf Twisting.} Suppose that $G$ is a smooth algebraic group over a perfect field $k$. If $\sigma \in \Z^1(k,G)$ is a Galois cocycle (that is, a continuous map $\sigma: \Gamma_k \to G(\bar{k})$ satisfying $\sigma_{st} = \sigma_s \tensor[^s]{\sigma}{_t}$ for all $s,t \in \Gamma_k$), we may {\em twist} $G$ by it to obtain an inner $k$-form $G^{\sigma}$, which can be defined by Galois descent as follows. As a $\bar{k}$-group, $\bar{G}^{\sigma}:=\bar{G}$, but the twisted Galois action on $G^{\sigma}(\bar{k}) = G(\bar{k})$ is given by
    \begin{equation*}
        \Gamma_k \times G^{\sigma}(\bar{k}) \to G^{\sigma}(\bar{k}), \quad (s,g) \mapsto \sigma_s \tensor[^s]{g}{} \sigma_s^{-1}.
    \end{equation*}
We may also twist the $k$-{\em variety} $G$ to obtain a $k$-form $G_{\sigma}$, whose Galois action is given by
    \begin{equation*}
        \Gamma_k \times G_{\sigma}(\bar{k}) \to G_{\sigma}(\bar{k}), \quad (s,g) \mapsto \tensor[^s]{g}{}  \sigma_s^{-1}.
    \end{equation*}
Then, one has $G_\sigma(k) \neq \varnothing$ if and only if $[\sigma] = 1$, the distinguished element of $\H^1(k,G)$. Equipped with the obvious translation actions (which are {\em a priori} defined over $\bar{k}$), the variety $G_{\sigma}$ is a {\em left} torsor under $G$ and a right torsor under $G^{\sigma}$. When $X$ is a $k$-variety and $f: Y \to X$ is a torsor under $G$ (more precisely, under $G_X:=G \times_k X$), we define the twisted torsor $f^{\sigma}: Y^{\sigma} \to X$ to be the contracted product $Y \times_k^G G_{\sigma}$, {\em i.e.}, the quotient of $Y \times_k G_{\sigma}$ by the left action $g \cdot (y,h) := (y \cdot g^{-1}, gh)$ of $G$. This is a right torsor under $G^{\sigma}$. Thus, we obtain a twisting map
    \begin{equation*}
        \theta_\sigma: \H^1(X,G) \to \H^1(X,G^{\sigma}), \quad [Y] \mapsto [Y^\sigma].
    \end{equation*}
For $X = \Spec(k)$, the map $\theta_\sigma$ takes $[\sigma]$ to $1$ and takes $1$ to $[G_{\sigma}]$. In general, we may ``evaluate'' $Y$ at any point $x \in X(k)$ to obtain a class $[Y](x):=[Y_x]=[f^{-1}(x)] \in \H^1(k,G)$. Then, we have $x \in f^{\sigma}(Y^{\sigma}(k))$ if and only if $[Y^{\sigma}](x) = 1$ in $\H^1(k,G^{\sigma})$, or $[Y](x) = [\sigma]$ in $\H^1(k,G)$. If $\tau \in \Z^1(k,G)$ is another cocycle, then $G^\tau = (G^{\sigma})^{\theta_\sigma(\tau)}$, $Y^{\tau} = (Y^{\sigma})^{\theta_\sigma(\tau)}$, and $f^{\tau} = (f^{\sigma})^{\theta_\sigma(\tau)}$. If $G$ is commutative, then $G^{\sigma} = G$, $[G_\sigma] = -[\sigma]$ in $\H^1(k,G)$, and $[Y^\sigma] = [Y] - p^\ast[\sigma]$ in $\H^1(X,G)$, where $p: X \to \Spec(k)$ is the structure morphism. 

Finally, if $k$ is number field, we define the pointed set $\Pbb^1(k,G) \subseteq \prod_{v \in \Omega_k} \H^1(k_v,G)$ as a restricted product in a similar way as in the commutative case. We also put
    \begin{equation*}
        \Sha(G) := \Ker(\H^1(k,G) \to \Pbb^1(k,G)).
    \end{equation*}
This means, for any $k$-torsor $P$ under $G$, one has $[P] \in \Sha(G)$ if and only if $P(\Abb_k) \neq \varnothing$.

{\bf Cartier duality.} Let $k$ be a field. If $G$ is an algebraic $k$-group, we denote by $\wh{G} := \iHom_k(G,\Gbb_m)$ its Galois module of geometric characters, that is, $\wh{G} = \Hom_{\bar{k}}(\bar{G},\Gbb_{m,\bar{k}})$ as a (finitely generated) abelian group, equipped with the action of $\Gamma_k$ defined by $(\tensor[^s]{\chi}{})(g) := \tensor[^s]{\chi}{}(\tensor[^{s^{-1}}]{g}{})$ for all $s \in \Gamma_k, \chi \in \wh{G}$, and $g \in \bar{G}$. A $k$-{\em group of multiplicative type} is a $k$-form of $\Gbb_m^n \times_k \prod_{i=1}^m \mu_{n_i}$, where $n \ge 0$ and $n_i \ge 1$. A $k$-{\em torus} is a connected $k$-group of multiplicative type ({\em i.e.}, a $k$-form of $\Gbb_m^n$). The formation $G \mapsto \wh{G}$ is an additive exact anti-equivalence of categories between $k$-groups of multiplicative type and finitely generated abelian groups equipped with a continuous action of $\Gamma_k$, called {\em Cartier duality}. This extends to a duality between groups of multiplicative type and twisted constant group schemes over an arbitrary base. We refer to \cite{Grothendieck1970Multiplicatif1,Grothendieck1970Multiplicatif2,Grothendieck1970Multiplicatif3} for a complete exposition.

A $k$-torus $T$ is said to be {\em quasi-trivial} if it is isomorphic to $\Res_{A/k}(\Gbb_{m,A})$ for some \'etale $k$-algebra $A$, where $\Res_{A/k}$ denotes the restriction  of scalars {\em \`a la Weil}. This is equivalent to saying that $\wh{T}$ is a $\Gamma_k$-permutation module, {\em i.e.}, it has a $\Zbb$-basis which is $\Gamma_k$-stable. In this case, one has $\H^1(k,T) = 0$ by Shapiro's lemma and Hilbert's Theorem 90.

{\bf Brauer group.} By convention, the Brauer group of a scheme $X$ is always the Brauer--Grothendieck group $\Br(X) := \H^2(X,\Gbb_m)$. For a smooth integral variety $X$ over a field $k$ of characteristic $0$, Grothendieck's purity theorem \cite[Theorem 3.7.7]{CTS2021Brauer} asserts that
    \begin{equation} \label{eq:GrothendieckPurity}
        \Br(X) = \bigcap_{\xi \in X^{(1)}} \Ker(\del_\xi: \Br (k(X)) \to \H^1(k(\xi),\Qbb/\Zbb)),
    \end{equation}
where $\del_\xi$ are the {\em residue} maps, and $X^{(1)}$ is the set of points of codimension $1$ of $X$. In particular, $\Br(X)$ is torsion. The {\em unramified} Brauer group $\Br_{\nr}(X)$ is the subgroup of $\Br(k(X))$ consisting of elements $\alpha$ whose residue $\del_\Ocal(\alpha)$ at every discrete valuation ring $\Ocal \supseteq k$ with $\Frac(\Ocal) = k(X)$ vanishes (in particular, $\Br_{\nr}(X) \subseteq \Br(X)$). If $X^c$ is a smooth compactification of $X$, then $\Br_{\nr}(X) = \Br(X^c)$. The group $\Br_{\nr}$ is a stable birational invariant of smooth integral varieties \cite[Proposition 6.2.7, Corollaries 6.2.10 and 6.2.11]{CTS2021Brauer}. 

When $k$ is an (Archimedean or nonarchimedean) local field, we identify $\Br(k)$ to a subgroup of $\Qbb/\Zbb$ using the local invariant map (this group is $\Qbb/\Zbb$ in the nonarchimedean case, $\{0\}$ if $k = \Cbb$, and $\Zbb/2$ if $k = \Rbb$). If $K/k$ is a finite extension, the local invariants of $K$ and $k$ are compatible with the corestriction map $\cores_{K/k}: \Br(K) \to \Br(k)$ (thus, the restriction $\res_{K/k}: \Br(k) \to \Br(K)$ corresponds to the multiplication by $[K:k]$ in $\Qbb/\Zbb$).

The {\em algebraic} Brauer group of $X$ is $\Br_1(X):=\Ker(\Br(X) \to \Br(\bar{X}))$ (recall that $\bar{X} = X \times_k \bar{k}$). If $k$ is a number field, we define $\Be(X)$ to be the subgroup of $\Br_1(X)$ consisting of ``everywhere locally constant'' elements, {\em i.e.}
    \begin{equation} \label{eq:FirstObstruction}
        \Be(X):=\Ker\left(\Br_1(X) \to \prod_{v \in \Omega_k} \frac{\Br_1(X_v)}{\Img(\Br(k_v) \to \Br(X_v))} \right).
    \end{equation}
A result of Harari \cite{Harari1994Fibration} assures that $\Be(X) \subseteq \Br_{\nr}(X)$ (for $X$ smooth geometrically integral). A variant of the product Brauer--Manin pairing \eqref{eq:BrauerManinPairingNr} is the {\em adelic Brauer--Manin pairing}
    \begin{equation} \label{eq:BrauerManinPairingAdelic} \tag{BM${}_{\Abb}$}
        \pair{-,-}_{\BM}: X(\Abb_k) \times \Br(X) \to \Qbb/\Zbb, \quad ((x_v)_{v \in \Omega_k},\alpha) \mapsto \sum_{v \in \Omega_k} \alpha(x_v),
    \end{equation}
For any subgroup $A \subseteq \Br(X)$, we denote by $X(\Abb_k)^{A}$ the subset of adelic points orthogonal to $A$ relative to the pairing \eqref{eq:BrauerManinPairingAdelic}. By the global reciprocity law, every adelic point is orthogonal to $\Img(\Br(k) \to \Br(X))$, and $X(k) \subseteq X(\Abb_k)^{\Br(X)}$. When $X$ is proper, one has $X(\Abb_k) = X(k_\Omega)$ by the valuative criterion for properness, $\Br_{\nr}(X) = \Br(X)$, and the pairings \eqref{eq:BrauerManinPairingNr} and \eqref{eq:BrauerManinPairingAdelic} coincide.

\section{Preliminary remarks} \label{sec:Preliminary} 

\subsection{Nonabelian descent types} \label{subsec:Type}

The notion of (nonabelian) descent types enables descent theory in the context where no torsors are given {\em a priori}. An example where the need of this object arises is as follows. If $G$ is a linear algebraic group over a field $k$, $H \subseteq G$ is a Zariski closed subgroup, and $X:=H \backslash G$, then the projection $G \to X$ is a torsor under $H$. But if $X$ is only a homogeneous space of $G$ (with $X(k) = \varnothing$, then we do not have a torsor $G \to X$. Worse, the stabilizer $\bar{H}$ of a point $x \in X(\bar{k})$ needs not be defined over $k$.

A complete theory of nonabelian descent types is available in  \cite{Linh2025Type}, which builds on the previous works of various authors on {\em kernels}\footnote{or {\em lien}, or {\em band}.} and {\em nonabelian Galois $2$-cohomology}. We refer to \cite[\S 1 and \S 2]{Borovoi1993H2}, \cite[\S 1]{FSS1998H2}, and \cite[\S 2.2]{DLA2019Reduction} for this subject. To summarize, a $k$-kernel (where $k$ is a perfect field) is a smooth algebraic $\bar{k}$-group $\bar{G}$ equipped with an {\em outer semilinear action} $\kappa$ of $\Gamma_k$ ({\em grosso modo}, it is a semilinear Galois action modulo conjugation) satisfying certain continuity and algebraicity conditions. If $L = (\bar{G},\kappa)$ is a $k$-kernel, the set $\H^2(k,L)$ of nonabelian Galois $2$-cohomology with coefficients in $L$ classifies the extensions of topological groups
    \begin{equation*}
        1 \to \bar{G}(\bar{k}) \to E \to \Gamma_k \to 1
    \end{equation*}
{\em bound by $L$}, that is, the outer action of $\Gamma_k$ on $\bar{G}(\bar{k})$ induced by this extension coincides with the one induced by $\kappa$. Here, we require that $\bar{G}(\bar{k})$ is discrete in $E$ and that the map $E \to \Gamma_k$ is open. If such an extension $E$ has a (continuous) section $\varsigma$, we say that its class $[E] \in \H^2(k,L)$ is {\em neutral}. The section $\varsigma$ defines a semilinear action on $\bar{G}$, hence also a $k$-form $G$ of $\bar{G}$ (by Galois descent). In this case, we call $G$ a $k$-form of $L$, and write $L = \lien(G)$. It might happen that the set $\H^2(k,L)$ does not have any neutral class; and a neutral class is not necessarily unique if it exists. If $\bar{G}$ is commutative, then the outer action $\kappa$ is in fact an action, hence $L$ has unique $k$-form $G$. In this case, the set $\H^2(k,G):=\H^2(k,\lien(G))$ is the usual Galois $2$-cohomology group with coefficients in $G$, and its unique neutral class is $0$.

Let $X$ be a nonempty $k$-variety. If $\pi: \bar{Y} \to \Spec(\bar{k})$ is a nonempty $\bar{k}$-variety equipped with a morphism $\bar{Y} \to \bar{X}:=X_{\bar{k}}$, we say that an automorphism $\alpha \in \Aut(\bar{Y}/X)$ is {\em semilinear} if there exists $s \in \Gamma_k$ such that $\pi \circ \alpha = \Spec(s^{-1}) \circ \pi$. Such an $s$ is necessarily unique. If $\SAut(\bar{Y}/X) \subseteq \Aut(\bar{Y}/X)$ denotes the subgroup of semilinear $X$-automorphisms of $\bar{Y}$, then we have an exact sequence
    \begin{equation} \label{eq:DescentTypeExactSequenceSAut}
        1 \to \Aut(\bar{Y}/\bar{X}) \to \SAut(\bar{Y}/X) \to \Gamma_k.
    \end{equation}
Let $L = (\bar{G},\kappa)$ be a $k$-kernel. A {\em descent type on $X$ bound by $L$} \cite[Definition 2.1]{Linh2025Type} is the equivalent class $\lambda = [\bar{Y},E]$ of a pair $(\bar{Y},E)$, where $\bar{Y} \to \bar{X}$ is a torsor under $\bar{G}$ and $E \subseteq \SAut(\bar{Y}/X)$ is a subgroup fitting in an exact sequence
    \begin{equation*}
	\xymatrix{
            1 \ar[r] & \bar{G}(\bar{k}) \ar[r] \ar@{_{(}->}[d] & E \ar[r] \ar@{_{(}->}[d] & \Gamma_k \ar[r] \ar@{=}[d] & 1 \\
            1 \ar[r] & \Aut(\bar{Y}/\bar{X}) \ar[r] & \SAut(\bar{Y}/X) \ar[r] & \Gamma_k
	}
    \end{equation*} 
of topological groups, where the top row is bound by $L$ and the bottom row is \eqref{eq:DescentTypeExactSequenceSAut}. Here, we impose that two such pairs $(\bar{Y},E)$ and $(\bar{Y}',E')$ are equivalent if there exists a $\bar{G}$-equivariant $\bar{X}$-isomorphism $\iota: \bar{Y} \to \bar{Y}'$ such that $\iota \circ \alpha \circ \iota^{-1} \in E'$ for all $\alpha \in E$. The set of descent types on $X$ bound by $L$ is denoted by $\DType(X,L)$.

Thus, a descent type $\lambda \in \DType(X,L)$ entails not only a torsor $\bar{Y} \to \bar{X}$ under $\bar{G}$, but also a class $\del(\lambda):=[E] \in \H^2(k,L)$. If $Y \to X$ is a torsor under a $k$-form $G$ of $L$, it defines a descent type $\dtype(Y) \in \DType(X,G):=\DType(X,\lien(G))$, whose attached $2$-cohomology class in $\H^2(k,G)$ is neutral. The converse also holds.

\begin{lemm} \label{lemm:DescentType}
    Let $X$ be a quasi-projective, reduced, and geometrically connected variety over a perfect field $k$. Let $L = (\bar{G},\kappa)$ be a $k$-kernel, and $\lambda = [\bar{Y},E] \in \DType(X,L)$.

    \begin{enumerate}
        \item \label{lemm:DescentType1} The class $[E] \in \H^2(k,L)$ is neutral if and only if there exists a torsor $Y \to X$ under a $k$-form $G$ of $L$, of descent type $\dtype(Y) = \lambda$.

        \item \label{lemm:DescentType2} If $X(k) \neq \varnothing$, then the class $[E]$ is neutral.

        \item \label{lemm:DescentType3} Let $Y \to X$ be a torsor under a $k$-form $G$ of $L$, of descent type $\dtype(Y) = \lambda$. If $Y' \to X$ is a torsor under a $k$-form $G'$ of $L$, then $\dtype(Y') = \lambda$ if and only if there exists a  cocycle $\sigma \in \Z^1(k,G)$ such that $G' = G^{\sigma}$ and the $X$-torsors $Y'$ and $Y^{\sigma}$ under $G^{\sigma}$ are isomorphic.
    \end{enumerate}
\end{lemm}
\begin{proof}
    See \cite[Proposition 2.4]{Linh2025Type}.
\end{proof}

The set $\DType(X,L)$ enjoys a functoriality property in $X$ (see \cite[Construction 2.5]{Linh2025Type}) as follows. Let $f: X' \to X$ be a morphism of nonempty $k$-varieties. If $\lambda = [\bar{Y},E] \in \DType(X,L)$, we may define the pullback type $f^\ast \lambda = [\bar{Y}',E'] \in \DType(X',L)$, where $\bar{Y}' := \bar{Y} \times_{\bar{X}} \bar{X}'$ and $\del(f^\ast \lambda) = [E'] = [E] = \del(\lambda)$ in $\H^2(k,L)$. Furthermore, for any overfield $K/k$, the $k$-kernel $L$ restricts to a $K$-kernel $L_K$, and there is a restriction map
    \begin{equation*}
        \res_{K/k}: \DType(X,L) \to \DType(X_K,L_K), \quad \lambda \mapsto \lambda_K,
    \end{equation*}
which is compatible with the restriction map $\H^2(k,L) \to \H^2(K,L_K)$ on $2$-cohomology.

There is also a partial functoriality with respect to $L$, which allows some sort d\'evissage argument on descent types. Let $L' = (\bar{G}',\kappa')$ be a second $k$-kernel, and $\varphi: L \to L'$ a surjective morphism of $k$-kernels (that is, a surjective morphism $\bar{G} \to \bar{G}'$ of algebraic $\bar{k}$-groups which is compatible with the outer Galois actions $\kappa$ and $\kappa'$). Let $\bar{H}:= \Ker(\bar{G} \xrightarrow{\varphi} \bar{G}')$.

\begin{lemm} \label{lemm:DescentTypePushforward}
    With the above notations, the following are true.
    \begin{enumerate}
        \item \label{lemm:DescentTypePushforward1} There is a canonical pushforward map $\varphi_\ast: \DType(X,L) \to \DType(X,L')$. If a descent type $\lambda \in \DType(X,L)$ has the corresponding torsor $\bar{Y} \to \bar{X}$ under $\bar{G}$, then the $\bar{X}$-torsor under $\bar{G}'$ associated with $\varphi_\ast \lambda \in \DType(X,L')$ is $\bar{Y}/\bar{H}$.

        \item \label{lemm:DescentTypePushforward2} Let $\lambda = [\bar{Y},E] \in \DType(X,L)$. Then, every torsor $Z \to X$ of descent type $\dtype(Z) = \varphi_\ast \lambda$ gives rise to a $k$-kernel $L_Z = (\bar{H},\kappa_Z)$ and a descent type $\lambda_Z = [\bar{Y},E_Z] \in \DType(Z,L_Z)$.

        \item \label{lemm:DescentTypePushforward3} Let $Z$, $G'$, $L_Z$, and $\lambda_Z$ be as in \ref{lemm:DescentTypePushforward2}. Then, for any torsor $Y \to Z$ under a $k$-form $H$ of $L_Z$, of descent type $\lambda_Z$, the composite $Y \to Z \to X$ is a torsor under a $k$-form $G$ of $L$, of descent type $\lambda$. In particular, one has $G' = G/H$ and $Z = Y/H$.
    \end{enumerate}
\end{lemm}
\begin{proof}
    See \cite[Construction 2.6, Lemma 2.7, Proposition 2.11]{Linh2025Type}.
\end{proof}

For abelian descent types, the set $\DType$ has a very nice description.

\begin{lemm} \label{lemm:DescentTypeAbelian}
    Let $G$ be a commutative algebraic $k$-group. Then, the set $\DType(X,G)$ is equipped with a structure of abelian group, functorial in $X$ and $G$ (even for nonsurjective morphisms $G \to G'$ of commutative algebraic $k$-groups), and fits in a complex
        \begin{equation} \label{eq:DescentTypeAbelian}
            \H^1(k,G) \to \H^1(X,G) \xrightarrow{\dtype} \DType(X,G) \xrightarrow{\del} \H^2(k,G) \to \H^2(X,G).
        \end{equation}
    If $X$ is quasi-projective, reduced, and geometrically connected, then \eqref{eq:DescentTypeAbelian} is exact. 
\end{lemm}
\begin{proof}
    See \cite[Construction 2.6, Lemma 2.8, Theorem 3.8]{Linh2025Type}.
\end{proof}

We note that, if $G$ is a $k$-group of multiplicative type and $X$ is a smooth quasi-projective geometrically integral $k$-variety, then \eqref{eq:DescentTypeAbelian} is equivalent to the ``fundamental exact sequence'' for {\em extended types} by Harari--Skorobogatov \cite[Proposition 8.1]{HS2013Descent} (see \cite[Theorem 4.1]{Linh2025Type}).

\label{FiniteDescentTYpe} If $X$ is geometrically connected, a {\em finite descent type} on $X$ in the sense of Harpaz--Wittenberg \cite[Definition 2.1]{HW2024Supersolvable} is a connected \'etale cover $\bar{Y}$ of $\bar{X}$ such that the composite $\bar{Y} \to \bar{X} \to X$ is Galois ({\em i.e.}, the field extension $\bar{k}(\bar{Y})/k(X)$ is Galois). In this case, the homotopy exact sequence
    \begin{equation*}
        1 \to \Aut(\bar{Y}/\bar{X}) \to \Aut(\bar{Y}/X) \to \Gamma_k \to 1
    \end{equation*}
yields an outer action $\kappa$ of $\Gamma_k$ on $\bar{G}:=\Aut(\bar{Y}/\bar{X})$ (viewed as a finite constant algebraic $\bar{k}$-group), hence a $k$-kernel $L = (\bar{G},\kappa)$, as well as a descent type $[\bar{Y},\Aut(\bar{Y}/X)]$ on $X$ bound by $L$.

\subsection{Abelianization of nonabelian cohomology} \label{subsec:Abelianization}

The proof of Theorem \ref{customthm:DescentConnected} ({\em i.e.}, Conjecture \ref{conj:Descent} for torsors under connected linear algebraic groups) requires the machinery of abelianization of nonabelian cohomology. Fix a base scheme $X$, which, for simplicity, is assumed to have residue characteristics $0$. Let us start with the following remark on $2$-term complexes of tori. Let $C = [T \xrightarrow{\rho} S]$ be a complex of $X$-tori, with $T$ in degree $-1$. Let $\wh{T}:=\iHom_X(T,\Gbb_m)$, $\wh{S}:=\iHom_X(S,\Gbb_m)$, and let $\wh{C} = [\wh{S} \xrightarrow{\rho^\ast} \wh{T}]$ (with $\wh{S}$ in degree $-1$) be the ``Cartier dual'' of $C$ . Following Demarche's construction \cite[p. 4]{Demarche2011PT}, we have a pairing
	\begin{equation} \label{eq:DescentPairingComplexTori1}
		C \otimes^{\Lbb} \wh{C} \to \Gbb_m[1]
	\end{equation}
in $\Dcal^+(X)$, as follows. Represent the derived tensor product $C \otimes^{\Lbb} \wh{C}$ by the ``total complex'' 
	\begin{equation*}
		\Tot(C \otimes \wh{C}):=[T \otimes \wh{S} \to (T \otimes \wh{T}) \oplus (S \otimes \wh{S}) \to S \otimes \wh{T}],
	\end{equation*}
concentrated in degrees $-2,-1,0$, where the first differential is $(\id_T \otimes \rho^\ast, -\rho \otimes \id_{\wh{S}})$ and the second one is $\rho \otimes \id_{\wh{T}} + \id_S \otimes \rho^\ast$. The composite
	\begin{equation*}
		(T \otimes \wh{T}) \oplus (S \otimes \wh{S}) \to \Gbb_m \oplus \Gbb_m \to \Gbb_m,
	\end{equation*}
where the first arrow is just the direct sum of the canonical pairings and the second one is the product map, gives rise to a morphism of complex $\Tot(C \otimes \wh{C}) \to \Gbb_m[1]$, hence a pairing as in  \eqref{eq:DescentPairingComplexTori1}. From its construction, this pairing is compatible with the canonical pairings between the tori $T$, $S$ and their respective Cartier duals. Since the pairing $T \otimes \wh{T} \to \Gbb_m$ induces an isomorphism $T \cong \iHom_X(\wh{T},\Gbb_m)$, and $\iExt_X^i(\wh{T},\Gbb_m) = 0$ for all $i > 0$ \cite[Lemme 1.3.3 (ii)]{CTS1987Descente} (and similarly for $S$), we deduce, using a d\'evissage argument, that \eqref{eq:DescentPairingComplexTori1} induces a quasi-isomorphism
	\begin{equation} \label{eq:DescentPairingComplexTori2}
		C \xrightarrow{\simeq} \Rbb\iHom_X(\wh{C},\Gbb_m)[1].
	\end{equation}
In particular, taking hypercohomology yields $\H^i(X,C) \cong \Ext^{i+1}_X(\wh{C},\Gbb_m)$ for all $i \ge -1$.

A {\em reductive $X$-group scheme} is a smooth affine group scheme $G \to X$ with connected reductive fibres \cite[D\'efinition 2.7]{Demazure1970Reductive1}. We use $G^{\ssimple}$ to denote its derived subgroup (whose fibres are connected semisimple), $G^{\tor} = G/G^{\ssimple}$ to denote its coradical torus \cite[D\'efinition \S{6.2}]{Demazure1970Reductive2}, and $G^{\sconnect}$ the universal covering of $G^{\ssimple}$ (whose fibres are simply connected semisimple) \cite[Exercise 1.6.13 (ii)]{Conrad2011Reductive}. The morphism is $G^{\sconnect} \to G^{\ssimple}$ is a central isogeny, whose kernel $\mu:=\pi_1(G^{\ssimple})$ is a finite \'etale commutative $X$-group scheme (the {\em algebraic fundamental group of $G^{\ssimple}$}). 

Associated with $G$ is a $2$-term complex of $X$-tori $C$ concentrated in degrees $-1$ and $0$, which is canonically defined up to a quasi-isomorphism, called the {\em algebraic fundamental complex of $G$}. It can be defined in terms of {\em $t$-resolutions} \cite[Definition 2.1, Propositions 2.2 and 4.1]{BGA2014Fundamental} as follows. A $t$-resolution of $G$ is a central extension
    \begin{equation*}
	1 \to T \to H \to G \to 1,
    \end{equation*}
where $T$ is an $X$-torus and $H$ is a reductive $X$-group scheme such that $H^{\ssimple} = H^{\sconnect}$. The complex $C$ is then represented by the cone of the composite $T \hookrightarrow H \twoheadrightarrow S:=H^{\tor}$, that is, the complex $[T \to S]$ (with $T$ in degree $-1$). The cohomology of $C$ is $\Hscr^{-1}(C) = \mu$ and $\Hscr^0(C) = S/T = G^{\tor}$. 

The {\em dual algebraic fundamental complex} $\pi_1^D(G)$ of $G$ is the complex $\wh{C}:=[\wh{S} \xrightarrow{\rho^\ast} \wh{T}]$, with $\wh{S}$ in degree $-1$. Its cohomology is $\Hscr^{-1}(\wh{C}) = \wh{G^{\tor}} =: \wh{G}$ and $\Hscr^0(\wh{C}) = \wh{\mu}$. If $G$ has a maximal torus $T^{\max}$ ({\em e.g.}, when $X = \Spec(k)$, where $k$ is a field of characteristic $0$), let $T: = T^{\max} \times_G G^{\sconnect}$ (which is maximal torus of $G^{\sconnect}$). Then, $C$ is represented by the complex $[T \to T^{\max}]$, with $T$ in degree $-1$ (see \cite[Lemma 3.9]{BGA2014Fundamental}). There are also the {\em abelianization maps}
\begin{equation*}
	\ab^i: \H^i(X,G) \to \H^i(X,C),
\end{equation*}
defined for $i = 0,1$. When $X = \Spec(k)$, they were introduced by Borovoi in \cite{Borovoi1998Reductive}. This was later generalized to arbitrary base schemes by Gonz\'alez-Avil\'es (combine \cite[(3.12)]{GA2012Abelianization} and \cite[Proposition 4.1]{BGA2014Fundamental}). There are two extreme cases where the maps $\ab^i$ have a simple description.
\begin{itemize}
	\item If $G^{\ssimple} = G^{\sconnect}$, then $C \simeq G^{\tor}$, and $\ab^i: \H^i(X,G) \to \H^i(X,G^{\tor})$ are the pushforward maps.
	
	\item If $G = G^{\ssimple}$, then $C \simeq \mu[1]$, and $\ab^i: \H^i(X,G) \to \H^{i+1}(X,\mu)$ are the connecting maps induced by the central extension $1 \to \mu \to G^{\sconnect} \to G \to 1$.
\end{itemize}

\begin{prop} \label{prop:DescentAbelianizationMaps}
    With the above notations, the following are true.
    \begin{enumerate}
	\item \label{prop:DescentAbelianizationMaps1} {\em (Breen's exact sequence)} We have a functorial exact sequence
		\begin{equation*}
			1 \to \mu(X) \to G^{\sconnect}(X) \xrightarrow{\rho} G(X) \xrightarrow{\ab^0} \H^0(X,C) \to \H^1(X,G^{\sconnect}) \xrightarrow{\rho_\ast} \H^1(X,G) \xrightarrow{\ab^1} \H^1(X,C)
		\end{equation*}
        of pointed sets. Furthermore, $\ab^0$ is a group homomorphism which induces an injection
		\begin{equation*}
			\H^0(X,C)/\ab^0(G(X)) \hookrightarrow \H^1(X,G^{\sconnect}).
		\end{equation*}
		
	\item \label{prop:DescentAbelianizationMaps2} If $X = \Spec(k)$, where $k$ is a local field of characteristic $0$, then $\ab^1: \H^1(k,G) \to \H^1(k,C)$ is surjective. It is bijective if $k$ is nonarchimedean.
		
	\item \label{prop:DescentAbelianizationMaps3} If $X = \Spec(k)$, where $k$ is a number field, then there is a Cartesian square
		\begin{equation*}
			\xymatrix{
				\H^1(k,G) \ar[d] \ar[r]^{\ab^1} & \H^1(k,C) \ar[d] \\
				\prod\limits_{v \in \Omega_k^{\infty}} \H^1(k_v,G) \ar[r]^{\ab^1} & \prod\limits_{v \in \Omega_k^{\infty}} \H^1(k_v,C)
			}
		\end{equation*}
		with surjective arrows. In particular, $\ab^1$ induces a bijection $\Sha(G) \xrightarrow{\cong} \Sha^1(k,C)$.
    \end{enumerate}
\end{prop}
\begin{proof}
    For \ref{prop:DescentAbelianizationMaps1}, see \cite[Theorem 1.1, Proposition 3.14 (a)]{GA2012Abelianization}. For \ref{prop:DescentAbelianizationMaps2} (resp. \ref{prop:DescentAbelianizationMaps3}), we refer to Theorem 5.4 and Corollary 5.4.1 (resp. Theorems 5.12 and 5.13) in \cite{Borovoi1998Reductive}.
\end{proof}

\begin{prop} \label{prop:DescentAbelianizationMapsTwisting}
    Let $G$ be a reductive $X$-group scheme and $f: Y \to X$ a right torsor under $G$. Let $P \to X$ be a left torsor under $G$ and let $G^P$ be the corresponding twisted $X$-group ({\em cf.} \cite[Chapitre III, \S{2.3}]{Giraud1971Cohomology}; in particular, $P \to X$ is also a right torsor under $G^P$). Let $f^P: Y^P \to X$ be the twist of $Y$ by $P$ (that is, $Y^P$ is the quotient of $Y \times_X P$ by the left action $(g,y,p) \mapsto (y \cdot g^{-1}, g \cdot p)$ of $G$), which is a right torsor under $G^P$. Denote by $C$ the algebraic fundamental complex of $G$, and by $\wh{C}$ the dual algebraic fundamental complex $\pi_1^D(G)$. 
	\begin{enumerate}
		\item \label{prop:DescentAbelianizationMapsTwisting1} The constructions $G \mapsto C$ and $G \mapsto \wh{C}$ commute with base change.
        \end{enumerate}
        Assume that $G$ admits a maximal torus $T^{\max}$. Then, the following are true.
        \begin{enumerate} \setcounter{enumi}{1}
		\item \label{prop:DescentAbelianizationMapsTwisting2}  Then, the algebraic fundamental complex (resp. the dual algebraic fundamental complex) of $G^P$ is canonically isomorphic to $C$ (resp. $\wh{C}$).
		
		\item \label{prop:DescentAbelianizationMapsTwisting3} We have $\ab^1([Y^P]) = \ab^1([Y]) - \ab^1([P])$ in $\H^1(X,C)$.
	\end{enumerate}
\end{prop}
\begin{proof}
    \begin{enumerate}
	\item This follows from the fact that $t$-resolutions are preserved by base change.
		
	\item It suffices to combine \cite[Definition 3.1]{GA2012Abelianization} and \cite[Proposition 4.1]{BGA2014Fundamental}, noting that the centers of $G$ and $G^{\sconnect}$ remains unchanged under twisting.
    
	\item See \cite[Proposition 3.11]{GA2012Abelianization}.
	\end{enumerate}
\end{proof}

\subsection{The relative units-Picard complex} \label{subsec:UPic}

Let $f: Y \to X$ be a morphism of schemes. Then we have a canonical morphism $\Gbb_{m,X}\to f_\ast \Gbb_{m,Y}$ of sheaves on $X_{\et}$. We define the {\em sheaf of relative units} of $Y/X$ to be 
    \begin{equation*}
        \Unit_{Y/X}:=\Coker(\Gbb_{m,X}\to f_\ast \Gbb_{m,Y}).
    \end{equation*}
We also recall the {\em relative Picard sheaf} $\Pic_{Y/X}:=\Rbb^1 f_\ast \Gbb_{m,Y}$. Following Gonz\'alez-Avil\'es \cite[(3.15)]{GA2018UPic}, we define the {\em relative units-Picard complex} $\UPic_{Y/X}$ to be the cone of the morphism 
    \begin{equation*}
        \Gbb_{m,X}[1] \to (\tau_{\le 1} \Rbb f_\ast \Gbb_{m,Y})[1].
    \end{equation*}
It is represented by a complex concentrated in degrees $-1$ and $0$, with cohomology $\Hscr^{-1}(\UPic_{Y/X}) = \Unit_{Y/X}$ and $\Hscr^0(\UPic_{Y/X}) = \Pic_{Y/X}$. By definition, we have a distinguished triangle
\begin{equation} \label{eq:DescentTriangleUPic}
	\Gbb_{m,X} \to \tau_{\le 1} \Rbb f_\ast \Gbb_{m,Y} \to \UPic_{Y/X}[-1] \to \Gbb_{m,X}[1],
\end{equation}
which splits if $f$ has a section \cite[Lemma 3.8]{GA2018UPic}. The class in $\Ext^2_X(\UPic_{Y/X},\Gbb_m)$ of the morphism $\UPic_{Y/X} \to \Gbb_{m,X}[2]$ in $\Dcal^+(X)$ associated with \eqref{eq:DescentTriangleUPic} is then an obstruction to the existence of sections of $f$. 

\begin{prop} \label{prop:DescentUPicOfTorsor}
    Let $X$ be a locally Noetherian regular scheme, $G$ a reductive $X$-group scheme, and $Y \to X$ a torsor under $G$. There exists an isomorphism $\varphi_Y: \UPic_{Y/X} \xrightarrow{\simeq} \pi_1^D(G)$ in $\Dcal^+(X)$, functorial in $Y$ and $G$. In particular, $\Unit_{Y/X} \cong \wh{G}$ and $\Pic_{Y/X} \cong \wh{\mu}$, where $\mu:=\Ker(G^{\sconnect} \to G^{\ssimple})$. 
\end{prop}
\begin{proof}
    The case where $Y = G$ is the trivial torsor was dealt with by Gonz\'alez-Avil\'es in \cite[Theorem 1.1]{GA2019Reductive}. The case where $X = \Spec(k)$ is due to Borovoi--van Hamel \cite[Lemma 5.2]{BvH2009Picard}. The proof in the general case is similar; we shall deploy Sansuc's argument, using the additivity theorem for $\UPic$. More precisely, the projections $p_Y: Y \times_X G \to Y$ and $p_G: Y \times_X G \to G$ induce a canonical morphism
	\begin{equation*}
		\varsigma_Y: \UPic_{Y/X} \oplus \UPic_{G/X} \to \UPic_{(Y \times_X G) / X}
	\end{equation*}
	in $\Dcal^+(X)$. By \cite[Proposition 4.4, Remark 5.2 (a)]{GA2018UPic}, $\varsigma_Y$ is an isomorphism. Note that the existence of an ``\'etale quasi-section'' of $Y \times_X G \to X$ as in (i) of {\em loc. cit.} is guaranteed by \cite[Corollaires 17.16.2 and 17.16.3 (ii)]{EGAIV4}, since this morphism is smooth \cite[Chapter III, Proposition 4.2]{Milne1980Etale} and surjective, and that for every generic point $\eta \in X$, the ``\'etale index'' \cite[Definition 2.3]{GA2018UPic} of $G_{\eta}$ is $1$ since it has a $\kappa(\eta)$-rational point. Let $\varphi_Y$ denote the composite
	\begin{equation*}
		\UPic_{Y/X} \xrightarrow{\alpha^\ast} \UPic_{(Y \times_X G) / X} \xrightarrow{\varsigma_Y^{-1}} \UPic_{Y/X} \oplus \UPic_{G/X} \xrightarrow{\pi_G} \UPic_{G/X}.
	\end{equation*}
	where $\alpha: Y \times_X G \to Y$ is the $X$-action of $G$ on $Y$, and $\pi_G$ is the second projection. According to its construction, $\varphi_Y$ is functorial in $Y$ and $G$. Since the result is known for the trivial torsor $G \to X$, it remains to show that $\varphi_Y$ is a quasi-isomorphism, that is, $\Hscr^{-1}(\varphi_Y): \Unit_{Y/X} \to \Unit_{G/X}$ and $\Hscr^0(\varphi_Y): \Pic_{Y/X} \to \Pic_{G/X}$ are isomorphisms of \'etale sheaves over $X$. This follows from the units-Picard-Brauer sequence, see \cite[Lemmata 2.5 and 3.2, Proposition 3.6, Corollary 3.8]{GA2019Reductive} and their proofs.
\end{proof}

Let $X$ be a locally Noetherian regular scheme, $G$ a reductive $X$-group scheme, and $f: Y \to X$ a torsor under $G$. Let $C$ (resp. $\wh{C}$) denote the fundamental algebraic complex (resp. the dual fundamental algebraic complex) of $G$. In view of Proposition \ref{prop:DescentUPicOfTorsor}, triangle \eqref{eq:DescentTriangleUPic} becomes
\begin{equation} \label{eq:DescentTriangleCHat}
	\Gbb_{m,X} \to \tau_{\le 1} \Rbb f_\ast \Gbb_{m,Y} \to \wh{C}[-1] \to \Gbb_{m,X}[1].
\end{equation}
Our aim is to relate the class in $\Ext^2_X(\wh{C},\Gbb_m)$ of the corresponding morphism $\wh{C} \to \Gbb_{m,X}[2]$ in $\Dcal^+(X)$ and the abelianized class $\ab^1([Y]) \in \H^1(X,C)$ under the identification \eqref{eq:DescentPairingComplexTori2}.

\begin{prop} \label{prop:DescentAb1Y}
    With the above notations, the isomorphism $\H^1(X,C) \xrightarrow{\cong} \Ext^2_X(\wh{C},\Gbb_m)$ induced by \eqref{eq:DescentPairingComplexTori2} takes $\ab^1([Y])$ to the inverse class of the morphism $\wh{C} \to \Gbb_{m,X}[2]$ from \eqref{eq:DescentTriangleCHat}, in the following cases.
    \begin{enumerate}
        \item \label{prop:DescentAb1Y1} $G^{\ssimple} = G^{\sconnect}$ (for example, when $G$ is a torus).
	
        \item \label{prop:DescentAb1Y2} The torsor $Y \to X$ is {\em isotrivial}, that is, there exists a finite \'etale surjective morphism $U \to X$ such that $Y_U \to U$ is a trivial torsor under $G_U$ (for example, if $X = \Spec(k)$, where $k$ is a field of characteristic $0$).
    \end{enumerate}
\end{prop}
\begin{proof}
	\begin{enumerate}
		\item When $G^{\ssimple} = G^{\sconnect}$, one has $C \simeq G^{\tor}$ and $\wh{C} \simeq \wh{G^{\tor}}[1] = \wh{G}[1]$, and the abelianization map $\ab^1: \H^1(X,G) \to \H^1(X,C)$ is simply the pushfoward along the projection $G \to G^{\tor}$. By functoriality, it is enough to prove the claim when $G = T$ is a torus. In this case, $\Pic_{Y/X} = 0$, and it was shown by Colliot-Th\'el\`ene and Sansuc that triangle \eqref{eq:DescentTriangleCHat} is quasi-isomorphic to the short exact sequence
		\begin{equation*} 
			0 \to \Gbb_{m,X} \to f_\ast \Gbb_{m,Y} \to \wh{T} \to 0,
		\end{equation*}
		provided by the relative version of Rosenlicht's lemma, whose class in $\Ext^1_X(\wh{T},\Gbb_m)$ turns out to be precisely $-[Y]$ (see \cite[Propositions 1.4.2 and 1.4.3]{CTS1987Descente}).
		
		\item We follow Kottwitz's argument in \cite[p. 369]{Kottwitz1986Trace}. Suppose that $\pi: U \to X$ is a surjective finite \'etale morphism trivializing $Y$. Take any $t$-resolution 
		\begin{equation*}
			1 \to T \to H \to G \to 1
		\end{equation*}
	of $G$, and consider the sheaf $T_1:=\pi_\ast T_U$ (which is representable by an $X$-torus). Let $\del: \H^1(X,G) \to \H^2(X,T)$ denote the connecting map. Then the restriction of $\del([Y])$ to $\H^2(U,T)$ is $0$. Now, by the global version of Shapiro's lemma (which can be obtained by combining \cite[Corollaire 5.6]{Grothendieck1973Foncteurs} and the Leray spectral sequence), one has $\H^2(U,T) \cong \H^2(X,T_1)$. Hence, the map $i_\ast: \H^2(X,T) \to \H^2(X,T_1)$ (induced by the inclusion $i: T \hookrightarrow T_1$) takes $\del([Y])$ to $0$. Let $H_1$ denote the pushout of $i$ and the inclusion $T \hookrightarrow H$. Then $H_1^{\ssimple} = H^{\ssimple} = H^{\sconnect}$ is simply connected, hence $1 \to T_1 \to H_1 \to G \to 1$ is also a $t$-resolution of $G$, and the corresponding connecting map $\delta: \H^1(X,G) \to \H^2(X,T_1)$ takes $[Y]$ to $0$. It follows that there is a torsor $Z \to X$ under $H_1$ such that the map $\H^1(X,H_1) \to \H^1(X,G)$ takes $[Z]$ to $[Y]$. We know from \ref{prop:DescentAb1Y1} that the statement of the proposition holds for $Z$. By functoriality, it also holds for $Y$.
	\end{enumerate}
\end{proof}

\begin{remk}
    When $X = \Spec(k)$, the statement of Proposition \ref{prop:DescentAb1Y} {\em without any restriction on $G$ or $Y$} was proved by Borovoi and van Hamel in \cite[Theorem 5.5]{BvH2009Picard}. The author, unfortunately, does not know if this result holds true in full generality. Should it be the case, a possible approach would be to adapt the proof from the toric case by Colliot-Th\'el\`ene and Sansuc \cite{CTS1987Descente} [Propositions 1.4.2 and 1.4.3] using $\check{\text{C}}$ech cocycles for \'etale hypercohomology. A catch there is that the quasi-isomorphism $\operatorname{UPic}_{Y/X} \simeq \pi_1^D(G)$ from Proposition \ref{prop:DescentUPicOfTorsor} is not explicit (it is constructed using the axioms of triangulated categories, {\em cf.} the proof of \cite[Theorem 1.1]{GA2019Reductive}). This would make the computation rather cumbersome, and, since the result in the isotrivial case is sufficient for further uses in the article, the author decides to refrain from including a proof for the general case (if it exists at all).
\end{remk}

\begin{coro} \label{coro:DescentAb1Y}
	Let $k$ be a field of characteristic $0$, $p: X \to \Spec(k)$ a smooth geometrically integral $k$-variety, $G$ a connected reductive linear algebraic $k$-group, and $f: Y \to X$ a torsor under $G$. Let $C$ (resp. $\wh{C}$) denote the algebraic fundamental complex (resp. dual algebraic fundamental complex) of $G$. Then, for every Galois cocycle $\sigma \in \Z^1(k,G)$, we have a distinguished triangle
	\begin{equation*}
		\Gbb_{m,X} \to \tau_{\le 1} \Rbb f_\ast^{\sigma} \Gbb_{m,Y^\sigma} \to \wh{C}_X[-1] \to \Gbb_{m,X}[1]
	\end{equation*}
	in $\Dcal^+(X)$. If either $G^{\ssimple} = G^{\sconnect}$ or the torsor $f: Y \to X$ is isotrivial ({\em cf.} Proposition \ref{prop:DescentAb1Y} \ref{prop:DescentAb1Y2}), then the map $\H^1(X,\wh{C}) \to \H^3(X,\Gbb_m)$ associated with the above triangle is given by $\wh{c} \mapsto (-\ab^1([Y]) + p^\ast[\sigma])\cup \wh{c}$, where the cup-product is induced by the pairing \eqref{eq:DescentPairingComplexTori1}.
\end{coro}
\begin{proof}
	This follows directly from Propositions \ref{prop:DescentAbelianizationMapsTwisting} and \ref{prop:DescentAb1Y}. We note that in the case where $f: Y \to X$ is isotrivial, the  twisted torsor $f^{\sigma}: Y^{\sigma} \to X$ is also isotrivial for any $\sigma \in \Z^1(k,G)$. This is because the class $[\sigma] \in \H^1(k,G)$ is trivialized by a finite extension of $k$.
\end{proof}

For further use in Section \ref{sec:DescentConnected}, we also require the following technical results.

\begin{prop} \label{prop:DescentTResolution}
    Let $X$ be an integral Noetherian regular scheme and let $1 \to T \to H \xrightarrow{p} G \to 1$ be a $t$-resolution of a reductive $X$-group scheme $G$. Let $Z \to X$ (resp. $Y \to X$) be a torsor under $H$ (resp. $G$) such that $p_\ast[Z] = [Y]$ ({\em i.e.} there exists a $p$-equivariant $X$-morphism $\pi: Z \to Y$ making $Z$ a $Y$-torsor under $T$). Then, the following are true.

    \begin{enumerate}
        \item \label{prop:DescentTResolution1} There is a distinguished triangle
            \begin{equation} \label{eq:DescentTResolution1}
                \UPic_{Y/X} \to \UPic_{Z/X} \to \wh{T}[1] \to \UPic_{Y/X}[1].
            \end{equation}

        \item \label{prop:DescentTResolution2} Triangle \eqref{eq:DescentTResolution1} is functorial in the following sense. Let $1 \to T' \to H' \xrightarrow{p'} G' \to 1$ be a $t$-resolution of another reductive $X$-group scheme $G'$. Let $Z' \to X$ (resp. $Y' \to X$) be a torsor under $H'$ (resp. $G'$) such that there exists a $p'$-equivariant $X$-morphism $\pi': Z' \to Y'$. Let $u: T' \to T$ be a morphism of $X$-tori and let $v: Z' \to Z$, $w: Y' \to Y$ be morphisms of $X$-schemes such that the diagram
            \begin{equation*}
                \xymatrix{
                    Z' \times_X T' \ar[rr]^-{(z,t) \mapsto (z \cdot t)} \ar[d]^{v \times_X u} && Z' \ar[d]^{v} \ar[rr]^{\pi'} && Y' \ar[d]^{w} \\
                    Z \times_X T \ar[rr]^-{(z',t') \mapsto (z' \cdot t')} && Z \ar[rr]^{\pi} && Y 
                }
            \end{equation*}
        commutes. Then we have morphism of distinguished triangles
            \begin{equation*}
                \xymatrix{
                    \UPic_{Y/X} \ar[r] \ar[d]^{w^\ast} & \UPic_{Z/X} \ar[r] \ar[d]^{v^\ast} & \wh{T}[1] \ar[r] \ar[d]^{u^\ast} & \UPic_{Y/X}[1] \ar[d]^{w^\ast[1]} \\
                    \UPic_{Y'/X} \ar[r] & \UPic_{Z'/X} \ar[r] & \wh{T'}[1] \ar[r] & \UPic_{Y'/X}[1]
                }
            \end{equation*}
    \end{enumerate}
\end{prop}
\begin{proof}
\begin{enumerate}
    \item We construct a commutative diagram 
        \begin{equation} \label{eq:DescentTResolution2}
                \xymatrix{
                    \UPic_{Y/X} \ar[r]^-{\pi^\ast} \ar[d]^{\varphi_Y} & \UPic_{Z/X} \ar[r] \ar[d]^{\varphi_Z} & \wh{T}[1] \ar[r]^-{\lambda_{\pi}[1]} \ar[d]^{\simeq} & \UPic_{Y/X}[1] \ar[d]^{\varphi_Y[1]} \\
                    \pi_1^D(G) \ar[r] & \pi_1^D(H) \ar[r] & \pi_1^D(T) \ar[r] & \pi_1^D(G)[1],
                }
            \end{equation}
        where the bottom row is the distinguished triangle associated with the given $t$-resolution of $G$ (see \cite[Proposition 4.12]{GA2019Reductive}). The quasi-isomorphisms $\varphi_Y$ and $\varphi_Z$ are given by Proposition \ref{prop:DescentUPicOfTorsor}. The unlabeled arrow on the top row of \eqref{eq:DescentTResolution2} is the composite
            \begin{equation*}
                \UPic_{Z/X} \xleftarrow{\simeq} \Unit_{Z/X}[1] \xrightarrow{\chi_\pi[1]} \wh{T},
            \end{equation*}
        where $\chi_\pi: \Unit_{Z/X} \to \wh{T}$ is the map from the units-Picard-Brauer sequence ({\em cf.} \cite[Corollary 3.8 (i)]{GA2019Reductive}), and where the inclusion $\Unit_{Z/X}[1] \hookrightarrow \UPic_{Z/X}$ is a quasi-isomorphism because $\Pic_{Z/X} = 0$ (note that $H^{\sconnect} = H^{\ssimple}$). The map $\lambda_\pi: \wh{T} \to \UPic_{Y/X}[1]$ is the {\em extended type} of the torsor $Z \xrightarrow{\pi} Y$ under $T$, constructed in \cite[Definition 8.2]{HS2013Descent}. Let us verify that diagram \eqref{eq:DescentTResolution2} indeed commutes. 
        
        To this end, we recall the construction of $\varphi_Y$ and $\varphi_Z$ from the proof of Proposition \ref{prop:DescentUPicOfTorsor}. The left square of \eqref{eq:DescentTResolution2} commutes because we have a commutative diagram
        \begin{equation*}
            \xymatrix{
                \UPic_{Y/X} \ar[d]^{\pi^\ast} \ar[r]^-{\alpha_Y^\ast} & \UPic_{(Y \times_X G)/X} \ar[d]^{(\pi \times_X p)^\ast} & \UPic_{Y/X} \oplus \UPic_{G/X} \ar[l]_-{\simeq} \ar[d]^{\pi^\ast \oplus p^\ast} \ar[r] & \UPic_{G/X} \ar[d]^{p^\ast} \\
                \UPic_{Z/X} \ar[r]^-{\alpha_Z^\ast} & \UPic_{(Z \times_X H)/X} & \UPic_{Z/X} \oplus \UPic_{H/X} \ar[l]_-{\simeq} \ar[r] & \UPic_{H/X},
            }
        \end{equation*}
    where $\alpha_Y: Y \times_X G \to Y$ and $\alpha_Z: Z \times_X H \to Z$ denote the respective action morphisms, and where right horizontal arrows are the second projections. Indeed, the left square commutes because $\alpha_Y \circ (\pi \times_X p) = \pi \circ \alpha_Z$ (the morphism $\pi$ being $p$-equivariant), and the other two squares obviously commute. 

    The middle square of \eqref{eq:DescentTResolution2} commutes because it corresponds to a commutative square
        \begin{equation*}
            \xymatrix{
                \Unit_{Z/X} \ar[d]^{\Hscr^{-1}(\varphi_Z)} \ar[rr]^-{\chi_\pi} && \wh{T} \ar@{=}[d]\\
                \wh{H} \ar[rr] && \wh{T}.
            }
        \end{equation*}
    Indeed, $\Hscr^{-1}(\varphi_Z)$ is just the map from the units-Picard-Brauer sequence for $Z$ ({\em cf.} \cite[Corollary 3.8 (i)]{GA2019Reductive}, and hence the commutativity of the above square follows from the functoriality of the same sequence (see {\em loc. cit.}, Remarks 3.7 (a)).

    As for the right square of \eqref{eq:DescentTResolution2}, since $T$ is central in $H$, we have a commutative diagram
        \begin{equation*}
                \xymatrix{
                    Z \times_X H \times_X T \times_X T \ar[rrr]^-{(z,h,t_1,h_2) \mapsto (z \cdot t_1, ht_2)} \ar[d]^{\alpha_Z \times_X \nabla} &&& Z \times_X H \ar[d]^{\alpha_Z} \ar[rrr]^{\pi \times_X p} &&& Y \times_X G \ar[d]^{\alpha_Y} \\
                    Z \times_X T \ar[rrr]^-{(z,t) \mapsto (z \cdot t)} &&& Z \ar[rrr]^{\pi} &&& Y, 
                }
            \end{equation*}
        where $\nabla: T \times_X T \to T$ denotes the multiplication morphism. By functoriality of the extended type construction, we obtain the commutativity of the left square the following diagram:
            \begin{equation*}
            \xymatrix{
                \wh{T} \ar[d]^{\lambda_\pi[-1]} \ar[r]^-{\nabla^\ast} & \wh{T \times_X T} \ar[d]^{\lambda_{\pi \times_X p}[-1]} & \wh{T} \oplus \wh{T} \ar[l]_-{\simeq} \ar[d]^{(\lambda_\pi \oplus \lambda_p)[-1]} \ar[r] & \wh{T} \ar[d]^{\lambda_p[-1]} \\
                \UPic_{Y/X} \ar[r]^-{\alpha_Y^\ast} & \UPic_{(Y \times_X G)/X} & \UPic_{Y/X} \oplus \UPic_{G/X} \ar[l]_-{\simeq} \ar[r] & \UPic_{G/X},
            }
        \end{equation*}
        The other squares in the above diagram obviously commute. The composite of the top row is $\id_{\wh{T}}$ (composing the first two arrows yields the diagonal map $\wh{T} \to \wh{T} \oplus \wh{T}$), and that of the bottom row is $\varphi_Y$ (when $\UPic_{G/X}$ is identified to $\pi_1^D(G)$). Furthermore, by unwinding the definition of the extend type in \cite[Proposition 8.1, Definition 8.2]{HS2013Descent} (see also \cite[\S {4.1}, pp. 415---416]{Linh2025Type}, one checks that the composite $\wh{T} \xrightarrow{\simeq} \UPic_{T/X} \to \UPic_{G/X}[1]$ is precisely the extended type $\lambda_p$ of the torsor $p: H \to G$ under $T$. This shows that the right square of \eqref{eq:DescentTResolution2} commutes.

        Since the bottom row of \eqref{eq:DescentTResolution2} is a distinguished triangle and since the vertical arrows there are quasi-isomorphisms, the top row is also a distinguished triangle.

        \item The catch here is that the morphism $v: Z' \to Z$ needs not correspond to a morphism $H' \to H$ of $X$-group schemes. Nevertheless, once the maps in \eqref{eq:DescentTResolution1} are explicitly constructed, our claims follow easily from the functoriality of the maps from the units-Picard-Brauer sequence as well as the extended types.
    \end{enumerate}
\end{proof}

\begin{prop} \label{prop:DescentContractedProduct}
    Let $X$ be an integral Noetherian regular scheme, $G$ a reductive $X$-group scheme, and $f: Y \to X$ a right torsor under $G$. Let $P \to X$ be a left torsor under $G$ and let $G^P$ be the corresponding twisted $X$-group. Let $f^P: Y^P \to X$ be the twist of $Y$ by $P$ ({\em cf.} Proposition \ref{prop:DescentAbelianizationMapsTwisting}). Let $C$ (resp. $\wh{C}$) denote the fundamental algebraic complex (resp. the dual fundamental algebraic complex) of $G$ and assume either of the following.
    \begin{enumerate}
        \item $G^{\ssimple} = G^{\sconnect}$ (for example, when $G$ is a torus).

        \item The $X$-torsors $Y$ and $P$ are isotrivial ({\em cf.} Proposition \ref{prop:DescentAb1Y} \ref{prop:DescentAb1Y2}).
    \end{enumerate}
    Then, we have a commutative square
    \begin{equation*}
        \xymatrix{
            \UPic_{Y^P/X} \ar[r]^-{\simeq} \ar[d]^{\pi^\ast} & \wh{C} \ar[d]^{\Delta}\\
            \UPic_{(Y \times_X P)/X} \ar[r]^-{\simeq} & \wh{C} \oplus \wh{C}.
        }
    \end{equation*}
    in $\Dcal^+(X)$, where the horizontal maps are the quasi-isomorphisms from Proposition \ref{prop:DescentUPicOfTorsor} (see also Proposition \ref{prop:DescentAbelianizationMapsTwisting} \ref{prop:DescentAbelianizationMapsTwisting2}), where $\pi: Y \times_X P \to Y^P$ is the canonical projection, and where $\Delta$ is the diagonal map.
\end{prop}
\begin{proof}
    The main difficulty here is that if $G$ is not commutative, the morphism $\pi$ is not equivariant with respect to any morphism $G \times_X G^P \to G^P$ of $X$-group schemes. We start with the case where $G = T$ is a torus. In this case, one has $G^P = T$, $\wh{C} \simeq \wh{T}[1]$ and $\pi$ is $\nabla$-equivariant, where $\nabla: T \times_X T \to T$ is the multiplication morphism. Since $\nabla$ induces the diagonal map $\Delta: \wh{T} \to \wh{T} \oplus \wh{T}$ on character modules, the claim in the toric case follows.

    Next, we consider the case where $G^{\ssimple} = G^{\sconnect}$. Then $\wh{C} \simeq \wh{G^{\tor}}[1]$. Let $Z$ (resp. $Q$) denote the pushforward of $Y$ (resp. $P$) along the projection $u: G \to G^{\tor}$, which are torsors under the $X$-torus $G^{\tor}$. Then, the projections $u^P: G^P \to (G^P)^{\tor} = G^{\tor}$ and $u \times_X u^P: G \times_X G^P \to (G \times_X G^P)^{\tor} = G^{\tor} \times_X G^{\tor}$ have the respective properties that $\pi^P_\ast Y^P = Z^Q$ and $(\pi \times_X \pi^P)_\ast (Y \times_X P) = Z \times_X Q$. Let $\varpi: Z \times_X Q \to Z^Q$ be the canonical projection. The commutativity of the square
    \begin{equation*}
        \xymatrix{
            Y \times_X P \ar[r] \ar[d]^{\pi} & Z \times_X Q \ar[d]^{\varpi} \\
            Y^P \ar[r] & Z^Q
        }
    \end{equation*}
    yields the commutativity of the left square in the following diagram:
        \begin{equation*}
        \xymatrix{
            \UPic_{Y^P/X} \ar[r] \ar[d]^{\pi^\ast} & \UPic_{Z^Q/X} \ar[r]^-{\simeq} \ar[d]^{\varpi^\ast} & \wh{G^{\tor}}[1] \ar[d]^{\Delta}\\
             \UPic_{(Y \times_X P)/X} \ar[r] & \UPic_{(Z \times_X Q)/X} \ar[r]^-{\simeq} & (\wh{G^{\tor}} \oplus \wh{G^{\tor}})[1].
        }
    \end{equation*}
    The right square commutes thanks to the toric case which was dealt with above.

    Now we consider the case where the torsors $Y \to X$ and $P \to X$ are isotrivial. Take a surjective finite \'etale morphism $U \to X$ trivializing both $Y$ and $P$. Then, by repeating the argument from the proof of Proposition Proposition \ref{prop:DescentAb1Y} \ref{prop:DescentAb1Y2}, we obtain a $t$-resolution
        \begin{equation*}
            1 \to T \to H \xrightarrow{u} G \to 1
        \end{equation*}
    such that $Y = u_\ast Z$ and $P = u_\ast Q$ for some right torsor $Z \to X$ and some left torsor $Q \to X$ under $H$. Note that 
        \begin{equation*}
            1 \to T \to H^Q \xrightarrow{u^Q} G^P \to 1 \quad \text{and} \quad 1 \to T \times_X T \to H \times_X H^Q \xrightarrow{u \times_X u^Q} G \times_X G^P \to 1
        \end{equation*}
    are also $t$-resolutions (recall that $T$ is central in $H$ hence remains unchanged under twisting), and that we have $u^Q_\ast Z^Q = Y^P$ and $(u \times_X u^Q)_\ast (Z \times_X Q) = Y \times_X P$. Let $\varpi: Z \times_X Q \to Z^Q$ (resp. $\nabla: T \times_X T \to T$) denote the canonical projection (resp. the multiplication morphism). One easily checks that the diagram    
        \begin{equation*}
            \xymatrix{
                Z \times_X Q \times_X T \times_X T \ar[rrr]^-{(z,q,t_1,h_2) \mapsto (z \cdot t_1, q \cdot t_2)} \ar[d]^{\varpi \times_X \nabla} &&& Z \times_X Q \ar[d]^{\varpi} \ar[rrr] &&& Y \times_X P \ar[d]^{\pi} \\
                Z^Q \times_X T \ar[rrr]^-{(z',t) \mapsto (z' \cdot t)} &&& Z^Q \ar[rrr] &&& Y, 
            }
        \end{equation*}
    commutes. By Proposition \ref{prop:DescentTResolution}, we have a morphism of distinguished triangles
        \begin{equation*}
                \xymatrix{
                    \UPic_{Y^P/X} \ar[r] \ar[d]^{\pi^\ast} & \UPic_{Z^Q/X} \ar[r] \ar[d]^{\varpi^\ast} & \wh{T}[1] \ar[r] \ar[d]^{\Delta} & \UPic_{Y^P/X}[1] \ar[d]^{\pi^\ast[1]} \\
                    \UPic_{(Y \times_X P)/X} \ar[r] & \UPic_{(Z \times_X Q)/X} \ar[r] & (\wh{T} \oplus \wh{T})[1] \ar[r] & \UPic_{(Y \times_X P)/X}[1].
                }
            \end{equation*}
   Now, since $H^{\ssimple} = H^{\sconnect}$, we deduce from the result in the previous case that the second column in the above diagram can be replaced by the map $\Delta[1]: \wh{S}[1] \to (\wh{S} \oplus \wh{S})[1]$, where $S:=H^{\tor}$. Since $\wh{C}$ is by definition the cone of of the morphism $\wh{S} \to \wh{T}$, we obtain the claimed result using the functoriality of the quasi-isomorphisms constructed from Proposition \ref{prop:DescentUPicOfTorsor}.
\end{proof}
\section{Torsors under connected linear algebraic groups} \label{sec:DescentConnected}

This section is devoted to the proof of Theorems \ref{customthm:DescentConnected} and \ref{customthm:DescentFinite}.

\subsection{Some calculations} \label{subsec:DescentConnectedPreliminaries}

In this subsection, we fix a field $k$ of characteristic $0$, a smooth geometrically integral $k$-variety $p: X \to \Spec(k)$, and a connected reductive linear algebraic group $G$ over $k$. Let $C$ (resp. $\wh{C}$) denote the algebraic fundamental complex (resp. dual algebraic fundamental complex) of $G$. Let $f: Y \to X$ be a torsor under $G$. Taking hypercohomology of the distinguished triangle \eqref{eq:DescentTriangleCHat}, we obtain an exact sequence
\begin{equation} \label{eq:DescentExactSequence1}
	\Br(X) \to \H^2(X,\tau_{\le 1} \Rbb f_\ast \Gbb_{m,Y}) \to \H^1(X,\wh{C}) \to \H^3(X,\Gbb_m),
\end{equation}
where the last arrow is given by $\wh{c} \mapsto -\ab^1([Y]) \cup \wh{c}$ if either $G^{\ssimple} = G^{\sconnect}$ or the torsor $Y \to X$ is isotrivial (Corollary \ref{coro:DescentAb1Y}). Let us investigate the two middle terms of \eqref{eq:DescentExactSequence1}. First, we need some lemmata.

\begin{lemm} \label{lemm:DescentUnramifiedBrauer}
    If $G$ is split, then $\Br_{\nr}(Y) \subseteq f^\ast(\Br(X))$.
\end{lemm}
\begin{proof}
    Take a split maximal torus $S$ of $G$, and let $T$ be its preimage in $G^{\sconnect}$. Then $[\wh{S} \to \wh{T}]$ (where $\wh{S}$ is in degree $-1$) is a complex of split tori representing $\wh{C}$. Let $K:=k(X)$ and let $V \to \Spec(K)$ be the generic fibre of $f$. By \cite[Corollary B]{BK2000Brauer}, the group $\Br_{\nr}(V)/ f^\ast(\Br(K))$ injects into the subgroup $\Sha^1_\omega(K,\wh{C})$ of $\H^1(K,\wh{C})$ consisting of elements that become trivial after restriction to every procyclic subgroup of $\Gamma_K$. Let $r:=\dim(T) = \dim(S)$, then $\wh{T} \cong \wh{S} \cong \Zbb^r$ as $\Gamma_K$-modules. In particular, one has $\H^1(K,\wh{T}) = 0$. It follows that $\H^1(K,\wh{C})$ injects into $\H^2(K,\wh{S})$. On the other hand, there is an exact sequence
	\begin{equation}
		0 \to \wh{S} \to \Qbb^r \to (\Qbb/\Zbb)^r \to 0
	\end{equation}
    of abelian groups equipped with the trivial actions of $\Gamma_K$. The term $\Qbb^r$ is uniquely divisible, hence {\em cohomologically trivial}. It follows that $\H^2(K,\wh{S}) = \H^1(K,(\Qbb/\Zbb)^r) = \Hom_{\cts}(\Gamma_K, (\Qbb/\Zbb)^r)$, and similarly for any procyclic subgroup of $\Gamma_K$. It follows that $\Sha^2_{\omega}(K,\wh{S}) = \Sha^1_{\omega}(K,(\Qbb/\Zbb)^r) = 0$, {\em a fortiori} $\Sha^1_{\omega}(K,\wh{C}) = 0$, and hence the map $f^\ast: \Br(K) \to \Br_{\nr}(V)$ is surjective.
	
	Viewing that $k(Y) = K(V)$, we have $\Br(Y) \subseteq \Br(V) \subseteq \Br(k(Y))$. Let $\alpha \in \Br_{\nr}(Y) \subseteq \Br_{\nr}(V)$. From what has been established above, $\alpha = f^\ast \beta$ for some $\beta \in \Br(K)$. It remains to show that $\beta \in \Br(X)$. Indeed, let $\xi \in X^{(1)}$ be any codimension $1$ point. Let $\eta$ be the generic point of the fibre $f^{-1}(\xi) \to \Spec(k(\xi))$, which is a codimension 1 point of $Y$. Since $\alpha \in \Br(Y)$, the residue $\del_\eta(\alpha) \in \H^1(k(\eta),\Qbb/\Zbb)$ vanishes ({\em cf.} \eqref{eq:GrothendieckPurity}). Since $f^{-1}(\xi)$ is geometrically integral, the map $\H^1(k(\xi),\Qbb/\Zbb) \to \H^1(k(\eta),\Qbb/\Zbb)$ is injective. By functoriality of the residue map \cite[Theorem 3.7.5]{CTS2021Brauer}, we have $\del_\xi(\beta) = 0$. It follows that $\beta \in \Br(X)$, which concludes the proof.
\end{proof}

Following Harpaz and Wittenberg \cite[\S{2.2}]{HW2020Galois}, we define the {\em relative algebraic Brauer group} 
\begin{equation*}
	\Br_1(Y/X):=\Ker(\Br(Y) \to \Br(\bar{Y}) / f^\ast (\Br(\bar{X}))).
\end{equation*}
It is obvious that $f^\ast(\Br(X)) \subseteq \Br_1(Y/X)$. By Lemma \ref{lemm:DescentUnramifiedBrauer}, one has $\Br_{\nr}(Y) \subseteq \Br_1(Y/X)$. 

\begin{lemm} \label{lemm:DescentRelativeAlgebraicBrauer}
	The group $\H^2(X, \tau_{\le 1}\Rbb f_\ast \Gbb_{m,Y})$ is identified to 
	\begin{equation*}
		\Ker(\Br(Y) \to \H^0(X, \Rbb^2 f_\ast \Gbb_{m,Y})),
	\end{equation*}
	and it contains $\Br_1(Y/X)$.
\end{lemm}
\begin{proof}
    Taking hypercohomology of the distinguished triangle
	\begin{equation*}
            \tau_{\le 1}\Rbb f_\ast \Gbb_{m,Y} \to \tau_{\le 2}\Rbb f_\ast \Gbb_{m,Y} \to \Rbb^2 f_\ast \Gbb_{m,Y}[-2] \to (\tau_{\le 1}\Rbb f_\ast \Gbb_{m,Y})[1]
	\end{equation*}
    in $\Dcal^+(X)$ yields an exact sequence
	\begin{equation} \label{eq:DescentRelativeAlgebraicBrauer1}
            0 \to \H^2(X,\tau_{\le 1} \Rbb f_\ast \Gbb_{m,Y}) \to \H^2(X,\tau_{\le 2} \Rbb f_\ast \Gbb_{m,Y}) \to \H^0(X,\Rbb^2 f_\ast \Gbb_{m,Y}).
	\end{equation}
    On the other hand, by taking cohomology of the distinguished triangle
	\begin{equation*}
		\tau_{\le 2}\Rbb f_\ast \Gbb_{m,Y} \to \Rbb f_\ast \Gbb_{m,Y} \to \tau_{\ge 3}\Rbb f_\ast \Gbb_{m,Y} \to (\tau_{\le 2}\Rbb f_\ast \Gbb_{m,Y})[1],
	\end{equation*}
    we have $\H^2(X,\tau_{\le 2} \Rbb f_\ast \Gbb_{m,Y}) = \H^2(X,\Rbb f_\ast \Gbb_{m,Y}) = \H^2(Y,\Gbb_m) = \Br(Y)$, noting that $\tau_{\ge 3}\Rbb f_\ast \Gbb_{m,Y}$ is acyclic in degrees $0$, $1$, $2$, and that $\Hbb(X,-) \circ \Rbb f_\ast \simeq \Hbb(Y,-)$ \cite[Corollary 10.8.3]{Weibel1994Homological}. It follows that \eqref{eq:DescentRelativeAlgebraicBrauer1} yields an identification 
	\begin{equation*}
		\H^2(X,\tau_{\le 1} \Rbb f_\ast \Gbb_{m,Y}) = \Ker(\Br(Y) \to \H^0(X,\Rbb^2 f_\ast \Gbb_{m,Y})).
	\end{equation*}
	Let us check that this group contains $\Br_1(Y/X)$. First, note that the composite 
        \begin{equation*}
            \Br(\bar{X}) \xrightarrow{f^\ast} \Br(\bar{Y}) \to \H^0(\bar{X},\Rbb^2 f_\ast \Gbb_{m,\bar{Y}})
        \end{equation*}
    is $0$ because $f^\ast$ factors through the map 
        \begin{equation*}
            \Br(\bar{X}) \to \H^2(\bar{X},\tau_{\le 1} \Rbb f_\ast \Gbb_{m,\bar{Y}}) = \Ker(\Br(\bar{Y}) \to \H^0(\bar{X},\Rbb^2 f_\ast \Gbb_{m,\bar{Y}}))
        \end{equation*}
    from \eqref{eq:DescentRelativeAlgebraicBrauer1}. In view of the commutative diagram with exact rows
	\begin{equation*}
		\xymatrix{
			0 \ar[r] & \H^2(X,\tau_{\le 1} \Rbb f_\ast \Gbb_{m,Y}) \ar[r] \ar[d] & \Br(Y) \ar[r] \ar[d] & \H^0(X,\Rbb^2 f_\ast \Gbb_{m,Y}) \ar[d] \\
			0 \ar[r] & \H^2(\bar{X},\tau_{\le 1} \Rbb f_\ast \Gbb_{m,\bar{Y}}) \ar[r] & \Br(\bar{Y}) \ar[r] & \H^0(\bar{X},\Rbb^2 f_\ast \Gbb_{m,\bar{Y}}),
		}
	\end{equation*} 
	we see that the inclusion $\Br_1(Y/X) \subseteq \H^2(X,\tau_{\le 1} \Rbb f_\ast \Gbb_{m,Y})$ would follow from the injectivity of the canonical map $\H^0(X,\Rbb^2 f_\ast \Gbb_{m,Y}) \to \H^0(\bar{X},\Rbb^2 f_\ast \Gbb_{m,\bar{Y}})$. Indeed, since $\Rbb^2 f_\ast \Gbb_{m,Y}$ is the sheaf on $X_{\et}$ associated to the presheaf $(X' \to X) \mapsto \Br(Y')$ (where $Y' := Y \times_X X'$), any element $a \in \H^0(X,\Rbb^2 f_\ast \Gbb_{m,Y})$ comes from some element $\alpha' \in \Br(Y')$ for some \'etale cover $X' \to X$. If the image of $a$ in $\H^0(\bar{X},\Rbb^2 f_\ast \Gbb_{m,\bar{Y}})$ is $0$, then there is a finite field extension $K/k$ and an \'etale cover $X'' \to X'_K$ such that the image of $\alpha'$ in $\Br(Y'')$ is $0$, where $Y'' := Y_K \times_{X_K} X'' = Y \times_X X''$. This means $\alpha'$ is trivialized by the \'etale cover $X'' \to X$ refining $X'$, hence $a = 0$.
\end{proof}

\begin{lemm} \label{lemm:DescentH1KC}
	We have $\H^1(k,\wh{C}) = \Ker(\H^1(X,\wh{C}) \to \H^1(\bar{X},\wh{C}))$.
\end{lemm}
\begin{proof}
    Represent $C$ by a complex $[T \to S]$ of $k$-tori (with $T$ in degree $-1$), then $\wh{C} \simeq [\wh{S} \to \wh{T}]$ (with $\wh{S}$ in degree $-1$). We have a commutative diagram
	\begin{equation*}
		\xymatrix{
			\H^1(k,\wh{S}) \ar[r] \ar[d] & \H^1(k,\wh{T}) \ar[r] \ar[d] & \H^1(k,\wh{C}) \ar[r] \ar[d] & \H^2(k,\wh{S}) \ar[r] \ar[d] & \H^2(k,\wh{T}) \ar[d] \\
			\H^1(X,\wh{S}) \ar[r] \ar[d] & \H^1(X,\wh{T}) \ar[r] \ar[d] & \H^1(X,\wh{C}) \ar[r] \ar[d] & \H^2(X,\wh{S}) \ar[r] \ar[d] & \H^2(X,\wh{T}) \ar[d] \\      \H^1(\bar{X},\wh{S}) \ar[r] & \H^1(\bar{X},\wh{T}) \ar[r] & \H^1(\bar{X},\wh{C}) \ar[r] & \H^2(\bar{X},\wh{S}) \ar[r] & \H^2(\bar{X},\wh{T})
		}
	\end{equation*}
	with exact rows. Since $\wh{S} \cong \Zbb^{\dim(S)}$ over $\bar{k}$, we have $\H^1(\bar{X},\wh{S}) = 0$. By the Hochschild--Serre spectral sequence $\H^p(k,\H^q(\bar{X},\wh{S})) \Rightarrow \H^{p+q}(X,\wh{S})$, one has $\H^1(k,\wh{S}) = \H^1(X,\wh{S})$ and $\H^2(k,\wh{S}) = \Ker(\H^2(X,\wh{S}) \to \H^2(\bar{X},\wh{S}))$. Similarly, we have $\H^1(\bar{X},\wh{T}) = 0$, $\H^1(k,\wh{T}) = \H^1(X,\wh{T})$, and $\H^2(k,\wh{T}) =
	\Ker(\H^2(X,\wh{T}) \to \H^2(\bar{X},\wh{T}))$. By a diagram chasing, we check that the map $\H^1(k,\wh{C}) \to \H^1(X,\wh{C})$ is injective, and that its image is precisely the kernel of the map $\H^1(X,\wh{C}) \to \H^1(\bar{X},\wh{C})$.
\end{proof}

In view of Lemmata \ref{lemm:DescentRelativeAlgebraicBrauer} and \ref{lemm:DescentH1KC}, we extract from \eqref{eq:DescentExactSequence1} an exact sequence
\begin{equation} \label{eq:DescentExactSequence2}
	\Br(X) \xrightarrow{f^\ast} \Br_1(Y/X) \xrightarrow{\varphi^1} \H^1(k,\wh{C}) \xrightarrow{\delta} \H^3(X,\Gbb_m),
\end{equation}
which lies in the heart of the proof of Theorem \ref{customthm:DescentConnected}. By Corollary \ref{coro:DescentAb1Y}, if either $G^{\ssimple} = G^{\sconnect}$ or the torsor $Y \to X$ is isotrivial, then the map $\delta$ in \eqref{eq:DescentExactSequence2} is given by $\delta(\wh{c}) = -\ab^1([Y]) \cup p^\ast \wh{c}$ (the cup-product being induced by the pairing \eqref{eq:DescentPairingComplexTori1}).

\subsection{The method of Harpaz and Wittenberg} \label{subsec:DescentConnectedProof}

Having established the exact sequence \eqref{eq:DescentExactSequence2}, we shall now prove in this subsection the following generalization of \cite[Th\'eor\`eme 2.1]{HW2020Galois}, which would imply Theorem \ref{customthm:DescentConnected}.

\begin{thm} \label{thm:DescentReductive} 
	Let $G$ be a connected reductive linear algebraic group over a number field $k$. Let $X$ be a smooth geometrically integral $k$-variety and $f: Y \to X$ a torsor under $G$. Let $A \subseteq \Br(X)$ denote the preimage of $\Br_{\nr}(Y)$ by the map $f^\ast: \Br(X) \to \Br(Y)$. Then
	\begin{equation*}
		X(\Abb_k)^A = \bigcup_{[\sigma] \in \H^1(k,G)} f^{\sigma}(Y^{\sigma}(\Abb_k)^{\Br_{\nr} (Y^{\sigma})}).
	\end{equation*}
\end{thm}

We keep the notations $k, G, X, Y, p, f, C$, and $\wh{C}$ from the previous subsection.

\begin{lemm} \label{lemm:DescentReductive1}
	Let $Z \to \Spec(k)$ be a torsor under $G$. Consider the exact sequence \eqref{eq:DescentExactSequence2} for $Z$:
	\begin{equation*}
		\Br(k) \to \Br_1(Z) \xrightarrow{\varphi} \H^1(k,\wh{C}) \to \H^3(k,\Gbb_m).
	\end{equation*}
	Then, for any field extension $K/k$ and any $z \in Z(K)$, $g \in G(K)$, $\alpha \in \Br_1(Z)$, one has 
	\begin{equation*}
		\alpha(z \cdot g) = \alpha(z) + \ab^0(g) \cup \res_{K/k}(\varphi(\alpha))
	\end{equation*} 
	in $\H^2(K,\Gbb_m) = \Br(K)$. Here, the cup-product is induced by the pairing \eqref{eq:DescentPairingComplexTori1}.
\end{lemm}
\begin{proof}
    If $G$ is a torus, this is \cite[Lemme 2.4]{HW2020Galois}. When $G^{\ssimple} = G^{\sconnect}$, we have $C \simeq G^{\tor}$, and the map $\ab^0$ is just the homomorphism $G(K) \to G^{\tor}(K)$. By functoriality and the toric case, the result also holds for the case where $G^{\ssimple} = G^{\sconnect}$. 
	
    In the general case, it suffices to repeat the argument in \cite[p. 369]{Kottwitz1986Trace}. Take any $t$-resolution
	\begin{equation*}
		1 \to T \to H \to G \to 1
	\end{equation*}
	of $G$ and let $\del: \H^1(k,G) \to \H^2(k,T)$ denote the connecting map. Let $K/k$ be a finite extension trivializing the class $\del([Z]) \in \H^2(k,T)$ and let $T_1:=\Res_{K/k}(T_K)$ (where $\Res$ denotes the restriction of scalars {\em \`a la} Weil). Let $i: T \hookrightarrow T_1$ denote the canonical inclusion. By Shapiro's lemma and \cite[Proposition 1.6.5]{NSW2008Cohomology}, the map $i_\ast: \H^2(k,T) \to \H^2(k,T_1) \cong \H^2(K,T)$  takes $\del([Z])$ to $0$. Thus, we obtain a $t$-resolution
	\begin{equation*}
		1 \to T_1 \to H_1 \xrightarrow{\pi} G \to 1,
	\end{equation*}
	where $H_1$ is the pushforward of $H$ along $i$, with the property that the connecting map $\H^1(k,G) \to \H^2(k,T_1)$ takes $[Z]$ to $0$, that is, there exists a $k$-torsor $W$ under $H_1$ and a $\pi$-equivariant morphism $W \to Z$. 
	
	Now, we consider the torsors $H_1 \to G$ and $W \to Z$ under $T_1$, and the fibres $(H_1)_g \to \Spec(K)$, $W_z \to \Spec(K)$. Take a finite extension $K'/K$ trivializing the classes $[(H_1)_g], [W_z] \in \H^1(K,T_1)$. Put $T_2:=\Res_{K'/K}((T_1)_{K'})$ and $T_3:=\Res_{K/k}(T_2) = \Res_{K'/k}(T_{K'})$, then $(T_3)_K \cong T_2^{[K:k]}$. Let $j: T_1 \to T_3$ denote the canonical inclusion. By Shapiro's lemma and \cite[Proposition 1.6.5]{NSW2008Cohomology}, the map $j_\ast: \H^1(K,T_1) \hookrightarrow \H^1(K,T_3) \cong \H^1(K,T_2)^{\oplus [K:k]} \cong \H^1(K',T_1)^{\oplus [K:k]}$ takes $[(H_1)_g]$ and $[W_z]$ to $0$. Thus, we obtain a $t$-resolution
	\begin{equation*}
		1 \to T_3 \to H_3 \to G \to 1,
	\end{equation*}
    where $H_3$ denotes the pushforward of $H_1$ along $j$, with the following properties.
	\begin{itemize}
		\item There exists a $k$-torsor $W' \to Z$ under $H_3$ such that the map $\H^1(k,H_3) \to \H^1(k,G)$ takes $[W']$ to $[Z]$.
		
		\item The point $g \in G(K)$ lifts to a point $h \in H_3(K)$.
		
		\item The point $z \in Z(K)$ lifts to a point $y \in W'(K)$.
	\end{itemize}
    Since the claim in the lemma holds for $W'$ (because $H_3^{\ssimple} = H_3^{\sconnect}$), it also holds for $Z$ by functoriality.
\end{proof}

Let $\sigma \in \Z^1(k,G)$ be a Galois cycle. The complexes $C$ and $\wh{C}$ remain unaffected after twisting by $G$ (see Proposition \ref{prop:DescentAbelianizationMapsTwisting} \ref{prop:DescentAbelianizationMapsTwisting2}). By Lemma \ref{lemm:DescentUnramifiedBrauer}, we have 
    \begin{equation*}
        \Br_{\nr}(Y^\sigma) \subseteq \Br_1(Y^{\sigma}/X):=\Ker(\Br(Y^\sigma) \to \Br(\bar{Y})/ f^\ast (\Br(\bar{X}))),
    \end{equation*}
and the exact sequence \eqref{eq:DescentExactSequence2} for the twisted torsor $f^\sigma: Y^\sigma \to X$ is
\begin{equation} \label{eq:DescentExactSequence3}
	\Br(X) \xrightarrow{(f^\sigma)^\ast} \Br_1(Y^\sigma / X) \xrightarrow{\varphi^{\sigma}} \H^1(k,\wh{C}) \xrightarrow{\delta} \H^3(X,\Gbb_m).
\end{equation}

The following is an analogue of \cite[Proposition 2.3]{HW2020Galois}.

\begin{prop} \label{prop:DescentReductive}
    Let $Z := G_\sigma$ be the left torsor under $G$ defined by $\sigma$ (which is also a right torsor under $G^{\sigma}$), so that $Y^{\sigma} = Y \times_k^G Z$ (see Subsection \ref{subsec:Notation}). Then, we have an exact sequence {(\em cf. \eqref{eq:DescentExactSequence2})}
	\begin{equation} \label{eq:DescentExactSequence4}
		\Br(k) \to \Br_1(Z) \xrightarrow{\varphi} \H^1(k,\wh{C}) \to \H^3(k,\Gbb_m).
	\end{equation}
	\begin{enumerate}
		\item \label{prop:DescentReductive1} Suppose that $\H^3(k,\Gbb_m) = 0$.
		Let $\pi: Y \times_k Z \to Y^{\sigma}$, $p_1: Y \times_k Z \to Y$, and $p_2: Y \times_k Z \to Z$ be the canonical projections. Then, for every $\alpha^{\sigma} \in \Br_{\nr}(Y^\sigma)$, there exist $\alpha \in \Br_1(Z)$ and $\alpha^1 \in \Br_{\nr}(Y)$ such that $\varphi^\sigma(\alpha^\sigma) = \varphi^1(\alpha^1) = \varphi(\alpha)$ in $\H^1(k,\wh{C})$ and $p_1^\ast \alpha^1 + p_2^\ast \alpha = \pi^\ast \alpha^{\sigma}$ in $\Br_1((Y \times_k Z) / X)$. 
		
		Also, the subgroup $\varphi^{\sigma}(\Br_{\nr}(Y^\sigma)) \subseteq \H^1(k,\wh{C})$ does not depend on $\sigma$.
		
		\item \label{prop:DescentReductive2} The preimage of $\Br_{\nr}(Y^{\sigma})$ by $(f^\sigma)^\ast: \Br(X) \to \Br(Y^\sigma)$ does not depend on $\sigma$.
	\end{enumerate}
\end{prop} 
\begin{proof}
	\begin{enumerate}
		\item Since $\H^3(k,\Gbb_m) = 0$, the map $\varphi$ from  \eqref{eq:DescentExactSequence4} is surjective, hence $\varphi^{\sigma}(\alpha^{\sigma}) = \varphi(\alpha)$ for some $\alpha \in \Br_1(Z)$. Take a nonempty open subset $U \subseteq X$ such that the restriction $V:=f^{-1}(U) \xrightarrow{f|_V} U$ is an isotrivial torsor under $G$. When we replace $f^{\sigma}$ by $f^{\sigma}|_{V^{\sigma}}$ in the sequence \eqref{eq:DescentExactSequence3}, the map $\delta: \H^1(k,\wh{C}) \to \H^3(U,\Gbb_m)$ is given by
		\begin{equation*}
			\delta(\wh{c}) = -\ab^1([V]) \cup p^\ast \wh{c} + p^\ast([\sigma] \cup \wh{c}) = -\ab^1([V]) \cup p^\ast \wh{c},
		\end{equation*}
		thanks to Corollary \ref{coro:DescentAb1Y} and the fact that $[\sigma] \cup \wh{c} \in \H^3(k,\Gbb_m) = 0$ (the cup-product being induced by the pairing \eqref{eq:DescentPairingComplexTori1}). Thus, this map does not depend on $\sigma$. It follows that the subgroup $\varphi^{\sigma}(\Br_1(V^{\sigma}/U)) \subseteq \H^1(k,\wh{C})$ does not depend on $\sigma$. Hence, we may write $\varphi^{\sigma}(\alpha^{\sigma}) = \varphi^1(\alpha^1)$ for some $\alpha^1 \in \Br_1(V/U)$. Let us now make the necessary modifications to $\alpha^1$ and $\alpha$ so that the conditions in the statement are fulfilled. From now on, the argument is the same as in \cite[Proposition 2.3 (i) and (ii)]{HW2020Galois}. Let $q: Z \to \Spec(k)$ denote the structure morphism and let $h:= f \times_k q = f \circ p_1 = f^{\sigma} \circ \pi$. Then $h: Y \times_k Z \to X$ is a torsor under $G \times_k G^{\sigma}$ (a connected reductive $k$-group with dual algebraic fundamental complex $\wh{C} \oplus \wh{C}$).  Let us now consider the distinguished triangles \eqref{eq:DescentTriangleCHat} (deduced from \eqref{eq:DescentTriangleUPic} and Proposition \ref{prop:DescentUPicOfTorsor}) associated with the torsors $f^{\sigma}|_{V^{\sigma}}: V^{\sigma} \to U$, $h|_{V \times_k Z}: V \times_k Z \to U$, $f|_{V}: V \to U$, and $q: Z \to \Spec(k)$. By functoriality and by Proposition \ref{prop:DescentContractedProduct}, we have a commutative diagram
        \begin{equation*}
            \xymatrix{
                \Gbb_{m,U} \ar[r] \ar@{=}[d] & \tau_{\le 1} \Rbb f^\sigma_\ast \Gbb_{m,V^{\sigma}} \ar[r] \ar[d]^{\pi^\ast} & \wh{C}_U[-1] \ar[r] \ar[d]^{\Delta[-1]} & \Gbb_{m,U}[1] \ar@{=}[d] \\
                \Gbb_{m,U} \ar[r]  & \tau_{\le 1} \Rbb h_\ast \Gbb_{m,V^{\sigma}} \ar[r] & (\wh{C}_U \oplus \wh{C}_U) [-1] \ar[r] & \Gbb_{m,U}[1]\\
                \Gbb_{m,U} \oplus \Gbb_{m,U} \ar[r] \ar[u]_{\nabla} & (\tau_{\le 1} \Rbb f_\ast \Gbb_{m,V}) \oplus p^\ast(\tau_{\le 1} \Rbb q_\ast \Gbb_{m,Z}) \ar[r] \ar[u]_{p_1^\ast + p_2^\ast} & (\wh{C}_U \oplus \wh{C}_U)[-1] \ar[r] \ar@{=}[u] & (\Gbb_{m,U} \oplus \Gbb_{m,U})[1], \ar[u]_{\nabla}
            }
        \end{equation*}
        where $\Delta: \wh{C}_U \to \wh{C}_U \oplus \wh{C}_U$ denotes the diagonal map and where $\nabla: \Gbb_{m,U} \oplus \Gbb_{m,U} \to \Gbb_{m,U}$ is the addition map. Now, the exact sequences of the shape \eqref{eq:DescentExactSequence2} associated with the rows of the above diagram are given by	
		\begin{equation*}
			\xymatrix{
				\Br(U) \ar[r]^-{(f^\sigma)^\ast} \ar@{=}[d] & \Br_1(V^{\sigma}/U) \ar[r]^-{\varphi^\sigma} \ar[d]^{\pi^\ast} & \H^1(k,\wh{C}) \ar[d]^{\Delta} \\
				\Br(U) \ar[r]^-{h^\ast} & \Br_1((V \times_k Z)/U) \ar[r] & \H^1(k,\wh{C}) \oplus \H^1(k,\wh{C}) \\ 
				\Br(U) \oplus \Br(k) \ar[u]_{\id_{\Br(U)} + p^\ast} \ar[r]^-{f^\ast \oplus q^\ast} & \Br_1(V/U) \oplus \Br_1 (Z) \ar[r]^{\varphi^1 \oplus \varphi} \ar[u]_{p_1^\ast + p_2^\ast} & \H^1(k,\wh{C}) \oplus \H^1(k,\wh{C}). \ar@{=}[u]
			}
		\end{equation*}
	Since $\Delta(\varphi^{\sigma}(\alpha^{\sigma})) = (\varphi^{\sigma}(\alpha^{\sigma}),\varphi^{\sigma}(\alpha^{\sigma})) = (\varphi^1(\alpha^1),\varphi(\alpha))$, the elements $p_1^\ast \alpha^1 + p_2^\ast \alpha$ and $\pi^\ast \alpha^{\sigma}$ have the same image in $\H^1(k,\wh{C}) \oplus \H^1(k,\wh{C})$. By exactness of the middle row, there exists $\beta \in \Br(U)$ such that $p_1^\ast \alpha^1 + p_2^\ast \alpha = \pi^\ast \alpha^{\sigma} + h^\ast \beta$, or 
		\begin{equation*}
			p_1^\ast(\alpha^1 - f^\ast \beta) + p_2^\ast \alpha = \pi^\ast \alpha^{\sigma}.
		\end{equation*} 
		On the other hand, $\varphi^1(\alpha^1-f^\ast \beta) = \varphi^1(\alpha^1)$ since $\varphi^1 \circ f^\ast = 0$. Replacing $\alpha^1$ by $\alpha^1 - f^\ast \beta$, we may assume that $\varphi^{\sigma}(\alpha^{\sigma}) + \varphi^1(\alpha^1) = \varphi(\alpha)$ and $p_1^\ast \alpha^1 + p_2^\ast \alpha = \pi^\ast \alpha^\sigma$. It remains to check that this implies $\alpha^1 \in \Br_{\nr}(Y)$. Indeed, since $\alpha^{\sigma} \in \Br_{\nr}(V^\sigma) = \Br_{\nr}(Y^{\sigma})$, we have $p_1^\ast \alpha^1 + p_2^\ast \alpha \in \Br_{\nr}(Y \times_k Z)$. Let $\xi$ be any codimension $1$ point on a smooth compactification of $Y$, and let $\eta$ denote the generic point $\xi \times_k Z$ (which is a $k(\xi)$-torsor under $G_{k(\xi)}$). Then $\eta$ is a codimension $1$ point on a smooth compactification of $Y \times_k Z$, hence $\del_{\eta}(p_1^\ast \alpha^1 + p_2^\ast \alpha) = 0$. It follows that $\del_{\eta}(p_1^\ast \alpha^1) = 0 \in \H^1(k(\eta),\Qbb/\Zbb)$ because $p_2^\ast \alpha \in \Br(Y \times_k Z)$. Since $\xi \times_k Z$ is geometrically integral, the restriction map $\H^1(k(\xi),\Qbb/\Zbb) \to \H^1(k(\eta),\Qbb/\Zbb)$ is injective. By functoriality of the residue map \cite[Theorem 3.7.5]{CTS2021Brauer}, we have $\del_\xi(\alpha^1) = 0 \in \H^1(k(\xi),\Qbb/\Zbb)$. Since $\xi$ is arbitrary, we conclude that $\alpha^1 \in \Br_{\nr}(Y)$.
		
		We have shown that $\varphi^{\sigma}(\Br_{\nr}(Y^{\sigma})) \subseteq \varphi^{1}(\Br_{\nr}(Y))$ for any torsor $Y \to X$ under any reductive group $G$ and any $\sigma \in \Z^1(k,G)$. If $\sigma,\tau \in \Z^1(k,G)$, then $\varphi^{\sigma}(\Br_{\nr}(Y^{\sigma})) \subseteq \varphi^{\tau}(\Br_{\nr}(Y^\tau))$ (resp. $\varphi^{\tau}(\Br_{\nr}(Y^\tau)) \subseteq \varphi^{\sigma}(\Br_{\nr}(Y^{\sigma}))$) thanks to a twisting argument by $\tau$ (resp. $\sigma$). It follows that $\varphi^{\tau}(\Br_{\nr}(Y^\tau)) = \varphi^{\sigma}(\Br_{\nr}(Y^{\sigma}))$.
		
	\item It suffices to repeat the proof of \cite[Proposition 2.3 (iii)]{HW2020Galois}.
	\end{enumerate}
\end{proof}

We can now proceed to prove Theorem \ref{thm:DescentReductive}. Keep the notations and assumptions therein.

\begin{proof} [Proof of Theorem \ref{thm:DescentReductive}]
    For any Galois cocycle $\sigma \in \Z^1(k,G)$, the subgroup $A \subseteq \Br(X)$ is the preimage of $\Br_{\nr}(Y^{\sigma})$ by $(f^\sigma)^\ast: \Br(X) \to \Br(Y^{\sigma})$, thanks to Proposition \ref{prop:DescentReductive} \ref{prop:DescentReductive2}. The inclusion
	\begin{equation*}
		\bigcup_{[\sigma] \in \H^1(k,G)} f^{\sigma}(Y^{\sigma}(\Abb_k)^{\Br_{\nr} (Y^{\sigma})}) \subseteq X(\Abb_k)^A
	\end{equation*}
    is then obvious, by functoriality of the Brauer--Manin pairing \eqref{eq:BrauerManinPairingAdelic}. Let us take a family $x_{\Abb} = (x_v)_{v \in \Omega_k} \in X(\Abb_k)^A$ and show that it can be lifted to $Y^{\sigma}(\Abb_k)^{\Br_{\nr} (Y^{\sigma})}$ for some $\sigma \in \Z^1(k,G)$. First, we try to lift $x_{\Abb}$ to $Y^{\sigma}(\Abb_k)$ by exploiting the orthogonality to $\Ker(f^\ast)$. This has already been done in \cite[Proof of Theorem 3.5]{CDX2019Compare}. We present here an alternative proof in the spirit of \cite[Proposition 2.5]{HW2020Galois} for the sake of uniformity.
	
    \begin{lemm} \label{lemm:DescentReductive2}
	There exists a Galois cocycle $\sigma \in \Z^1(k,G)$ such that $x_\Abb \in f^{\sigma}(Y^{\sigma}(\Abb_k))$.
    \end{lemm}
    \begin{proof}
	Since the torsor $Y \times_X Y \to Y$ under $G$ is trivial, the map $f^\ast: \H^1(X,G) \to \H^1(Y,G)$ takes $[Y]$ to $1$, hence $f^\ast(\ab^1([Y])) = 0$ in $\H^1(Y,C)$ by functoriality of $\ab^1$. It follows that  $f^\ast(\ab^1([Y]) \cup p^\ast \wh{c}) = 0$ for any $\wh{c} \in \H^0(k,\wh{C})$ (the cup-product being induced by the pairing \eqref{eq:DescentPairingComplexTori1}), or $\ab^1([Y]) \cup p^\ast \wh{c} \in \Ker (f^\ast) \subseteq A$. Since $x_\Abb \in X(\Abb_k)^A$, the family $\ab^1([Y](x_\Abb)) \in \Pbb^1(k,C)$ is orthogonal to $\H^0(k,\wh{C})$. By the Poitou--Tate exact sequence for $C$ (see \cite[Th\'eor\`eme 6.1]{Demarche2011PT}), there exists $c \in \H^1(k,C)$ such that $\ab^1([Y](x_v)) = \loc_v(c) \in \H^1(k_v,C)$ for all $v \in \Omega_k$.
		
    Let us now apply Proposition \ref{prop:DescentAbelianizationMaps} \ref{prop:DescentAbelianizationMaps3}. First, let $c = \ab^1([\sigma])$ for some cocycle $\sigma \in \Z^1(k,G)$, so that $\ab^1([Y](x_v)) = \ab^1(\loc_v([\sigma]))$ in $\H^1(k_v,C)$ for all $v \in \Omega_k$. By Proposition \ref{prop:DescentAbelianizationMapsTwisting} \ref{prop:DescentAbelianizationMapsTwisting3}, this is equivalent to $\ab^1([Y^{\sigma}](x_v)) = 0$. Twisting by $\sigma$, we may assume that $\ab^1([Y](x_v)) = 0$ for all $v \in \Omega_k$. If $v \in \Omega_k^{\f}$, this means $[Y](x_v) = 1 \in \H^1(k_v,G)$ by virtue of Proposition \ref{prop:DescentAbelianizationMaps} \ref{prop:DescentAbelianizationMaps2}. Since $\H^1(k,G)$ is the fibre product of $\prod_{v \in \Omega_k^{\infty}}\H^1(k_v,G)$ and $\H^1(k,C)$ over $\prod_{v \in \Omega_k^{\infty}} \H^1(k_v,C)$, there exists $\sigma \in \Z^1(k,G)$ such that $[Y](x_v) = \loc_v([\sigma])$ for all $v \in \Omega_k^{\infty}$ and $\ab^1([\sigma]) = 0$ (in particular, $\loc_v([\sigma]) = 1$ for $v \in \Omega_k^{\f}$). It follows that $[Y^\sigma](x_v) = 1$ in $\H^1(k_v,G^\sigma)$ for all $v \in \Omega_k$, {\em i.e.}, $x_\Abb \in f^\sigma(Y^\sigma(\Abb_k))$.
    \end{proof}
	
	Twisting by the cocycle provided by Lemma \ref{lemm:DescentReductive2}, we may assume that $x_\Abb \in f(Y(\Abb_k))$. Then $x_\Abb \in f^{\sigma}(Y^\sigma(\Abb_k))$ for any cocycle $\sigma$ such that $[\sigma] \in \Sha(G)$. For such a cocycle, we define a map
	\begin{equation} \label{eq:DescentReductive1}
		\varepsilon_\sigma: \Sha^1(k,\wh{C}) \cap \varphi^1(\Br_{\nr}(Y)) \to \Qbb/\Zbb, \quad \alpha' \mapsto \pair{y_\Abb,\alpha^{\sigma}}_{\BM},
	\end{equation}
	where $y_\Abb = (y_v)_{v \in \Omega_k} \in Y^{\sigma}(\Abb_k) = Y(\Abb_k)$ is any lifting of $x_\Abb$ and where $\alpha^{\sigma} \in \Br_{\nr} (Y^\sigma)$ is such that $\varphi^\sigma(\alpha^{\sigma}) = \alpha'$ (which exists by Proposition \ref{prop:DescentReductive} \ref{prop:DescentReductive1}). The map $\varepsilon_\sigma$ is well-defined:
	\begin{itemize}
		\item It does not depend on the choice of $\alpha^{\sigma}$. Indeed, upon replacing $\alpha^{\sigma}$ by $\alpha^{\sigma} + (f^{\sigma})^\ast \beta \in \Br_{\nr}(Y^{\sigma})$ (in view of the exact sequence \eqref{eq:DescentExactSequence3}) for some $\beta \in A$, one has
		\begin{equation*}
			\pair{y_\Abb, \alpha^{\sigma} + (f^{\sigma})^\ast \beta}_{\BM} = \pair{y_\Abb, \alpha^{\sigma}}_{\BM} + \pair{y_\Abb, (f^{\sigma})^\ast \beta}_{\BM} = \pair{y_\Abb,\alpha^{\sigma}}_{\BM} + \pair{f^{\sigma}(y_\Abb),\beta}_{\BM} = \pair{y_\Abb,\alpha^{\sigma}}_{\BM},
		\end{equation*}
		where the last equality is due to the fact that $f^{\sigma}(y_\Abb) = x_\Abb \in X(\Abb_k)^A$. 
		
		\item It does not depend on the choice of $y_{\Abb}$. Indeed, since $Y \to X$ is a torsor under $G$, upon replacing $y_{\Abb}$ by $y_\Abb \cdot g_\Abb:=(y_v \cdot g_v)_{v \in \Omega_k}$ for some $g_\Abb = (g_v)_{v \in \Omega_k} \in \Gbb(\Abb_k)$, one has
		\begin{equation*}
			\pair{y_\Abb \cdot g_\Abb, \alpha^{\sigma}}_{\BM} = \pair{y_\Abb, \alpha^{\sigma}}_{\BM} + \sum_{v \in \Omega_k} \ab^0(g_v) \cup \loc_v(\alpha'). 
		\end{equation*}
		by Lemma \ref{lemm:DescentReductive1} applied to each fibre $(f^\sigma)^{-1}(x_v) \to \Spec(k_v)$. But $\loc_v(\alpha') = 0$ for all $v \in \Omega_k$ since $\alpha' \in \Sha^1(k,\wh{C})$, hence $\pair{y_\Abb \cdot g_\Abb, \alpha^{\sigma}}_{\BM} = \pair{y_\Abb, \alpha^{\sigma}}_{\BM}$.
	\end{itemize} 
	In what follows, we denote by $\pair{-.-}_{\PT}: \Sha^1(k,C) \times \Sha^1(k,\wh{C}) \to \Qbb/\Zbb$ the global Poitou--Tate pairing induced by \eqref{eq:DescentPairingComplexTori1}, which is a perfect duality of finite groups (see \cite[Th\'eor\`eme 5.7]{Demarche2011PT}).
	
\begin{lemm} \label{lemm:DescentReductive3}
    Let $[\sigma] \in \Sha(G)$. Let $Z$ be the $k$-torsor under $G^{\sigma}$ in Proposition \ref{prop:DescentReductive} (hence $Z(\Abb_k) \neq \varnothing$) and $\varphi$ the map from \eqref{eq:DescentExactSequence4}. For any $\alpha \in \Be(Z)$ ({\em cf.} \eqref{eq:FirstObstruction}) and $z_\Omega \in Z(k_\Omega)$, one has 
        \begin{equation*}
            \pair{z_\Omega,\alpha}_{\BM} = -\pair{\ab^1([\sigma]),\varphi(\alpha)}_{\PT}.
        \end{equation*}
\end{lemm}
\begin{proof}
    Since $\alpha \in \Be(Z)$, it suffices to show the identity in the statement for any family $z_\Omega$, which can be chosen at our disposal. We deploy the same trick as in the proof of Lemma \ref{lemm:DescentReductive1}. If $G$ is a torus, the result is \cite[Lemme 8.4]{Sansuc1981Arithmetique}. When $G^{\ssimple} = G^{\sconnect}$, we have $C \simeq G^{\tor}$, and the map $\ab^1$ is just the natural pushforward $\H^1(k,G) \to \H^1(k,G^{\tor})$. Let $Z' \to \Spec(k)$ denote the $k$-torsor $Z \times_k^{G} G^{\tor}$ under $G_{\tor}$. Note that the map $\varphi$ induces an isomorphism $\Be(Z)/\Br(k) \xrightarrow{\cong} \Sha^1(k,\wh{C})$, and similarly, $\Be(Z')/\Br(k) \cong \Sha^1(k,\wh{C}) \cong \Be(Z)/\Br(k)$.	By functoriality and the toric case, the result also holds for the case $G^{\ssimple} = G^{\sconnect}$. In the general case, it suffices to show that there exist a surjection $\varpi: H \to G$, where $H$ is a connected reductive $k$-group with $H^{\ssimple} = H^{\sconnect}$ such that $[\sigma]$ lifts to $\Sha(H)$. Indeed, as in \cite[p. 369]{Kottwitz1986Trace}, there is a {\em $z$-extension} $1 \to T \to H \to G \to 1$, {\em i.e.}, a $t$-resolution with $T$ quasi-trivial, such that $[\sigma] \in \H^1(k,G)$ lifts to a class $[\tau] \in \H^1(k,H)$. Since $[\sigma] \in \Sha(G)$ and $\H^1(k_v,T) = 0$ for all $v \in \Omega_k$, one has $[\tau] \in \Sha(H)$. This concludes the proof of the lemma.
\end{proof}

Lemma \ref{lemm:DescentReductive3} allows us to relate the maps $\varepsilon_1$ and $\varepsilon_\sigma$ from \eqref{eq:DescentReductive1} as follows. Let $Z$ be as above, then $Y^{\sigma} = Y \times_k^G Z$. Let $\pi: Y \times_k Z \to Y^{\sigma}$, $p_1: Y \times_k Z \to Y$, and $p_2: Y \times_k Z \to Z$ be the canonical projections. For any $\alpha' \in \Sha^1(k,\wh{C}) \cap \varphi^1(\Br_{\nr}(Y))$, by Proposition \ref{prop:DescentReductive} \ref{prop:DescentReductive1}, there are $\alpha \in \Br_1(Z)$, $\alpha^1 \in \Br_{\nr}(Y)$, and $\alpha^{\sigma} \in \Br_{\nr}(Y^{\sigma})$ such that $\varphi(\alpha) = \varphi^1(\alpha^1) = \varphi^{\sigma}(\alpha^{\sigma}) = \alpha'$ and $p_1^\ast \alpha^1 + p_2^\ast \alpha = \pi^\ast \alpha^\sigma$. Take any $z_\Abb \in Z(\Abb_k)$, then
	\begin{align*}
		\varepsilon_\sigma(\alpha') & = \pair{y_\Abb,\alpha^{\sigma}}_{\BM}, & \text{by definition of } \varepsilon_\sigma,\\
		& = \pair{(y_\Abb,z_\Abb),\pi^\ast \alpha^{\sigma}}_{\BM}, & \text{by functoriality of the pairing \eqref{eq:BrauerManinPairingAdelic}},\\
		& = \pair{(y_\Abb,z_\Abb),p_1^\ast \alpha^1}_{\BM} + \pair{(y_\Abb,z_\Abb),p_2^\ast \alpha}_{\BM} \\
		& = \pair{y_\Abb,\alpha^1}_{\BM} + \pair{z_\Abb,\alpha}_{\BM}, & \text{by functoriality of the pairing \eqref{eq:BrauerManinPairingAdelic}}, \\
		& = \varepsilon_1(\alpha') - \pair{\ab^1([\sigma]), \alpha'}_{\BM}, & \text{by definition of $\varepsilon_1$ and Lemma \ref{lemm:DescentReductive3}}.
	\end{align*}
	Since the group $\Qbb/\Zbb$ is divisible, the homomorphism $\varepsilon_1$ can be extended to the whole group $\Sha^1(k,\wh{C})$. Since the Poitou--Tate pairing $\pair{-.-}_{\PT}: \Sha^1(k,C) \times \Sha^1(k,\wh{C}) \to \Qbb/\Zbb$ is a perfect duality and since the map $\ab^1: \H^1(k,G) \to \H^1(k,C)$ induces a bijection $\Sha(G) \xrightarrow{\cong} \Sha^1(k,C)$ by Proposition \ref{prop:DescentAbelianizationMaps} \ref{prop:DescentAbelianizationMaps3}, there exists $[\sigma] \in \Sha(G)$ such that $\varepsilon_1(\alpha') = \pair{\ab^1([\sigma]), \alpha'}_{\BM}$ for all $\alpha' \in \Sha^1(k,\wh{C}) \cap \varphi^1(\Br_{\nr}(Y))$. By the above calculation, we have $\varepsilon_\sigma = 0$. Twisting by the cocycle $\sigma$ as above, we may assume that $\varepsilon_1 = 0$. The rest of the proof is to show that this implies $x_\Abb \in f(Y(\Abb_k)^{\Br_{\nr}(Y)})$. As always, we fix an adelic point $y_\Abb \in Y(\Abb_k)$ lifting $x_\Abb$, and consider the homomorphism $\varepsilon: \Br_{\nr}(Y) \to \Qbb/\Zbb$ given by $\varepsilon(\alpha^1):=\pair{y_\Abb,\alpha^1}_{\BM}$. Then, one has
        \begin{equation*}
            \varepsilon|_{\Br_{\nr}(Y) \cap (\varphi^1)^{-1}(\Sha^1(k,\wh{C}))} = \varepsilon_1 \circ \varphi^1|_{\Br_{\nr}(Y) \cap (\varphi^1)^{-1}(\Sha^1(k,\wh{C}))} = 0.
        \end{equation*}
    Hence, $\varepsilon$ factors through a homomorphism $\bar{\varepsilon}: \frac{\Br_{\nr} (Y)}{\Br_{\nr} (Y) \cap (\varphi^1)^{-1}(\Sha^1(k,\wh{C}))} \to \Qbb/\Zbb$. On the other hand, the map $\varphi^1$ induces an injection $\bar{\varphi}^1: \frac{\Br_{\nr} (Y)}{\Br_{\nr} (Y) \cap (\varphi^1)^{-1}(\Sha^1(k,\wh{C}))} \hookrightarrow \frac{\H^1(k,\wh{C})}{\Sha^1(k,\wh{C})}$.
	Since $\Qbb/\Zbb$ is a divisible group, there exists a map $\bar{\varepsilon}': \frac{\H^1(k,\wh{C})}{\Sha^1(k,\wh{C})} \to \Qbb/\Zbb$ such that $\bar{\varepsilon}' \circ \bar{\varphi}^1 = \bar{\varepsilon}$, or equivalently, a map $\varepsilon': \H^1(k,\wh{C}) \to \Qbb/\Zbb$ such that $\varepsilon'|_{\Sha^1(k,\wh{C})} = 0$ and $\varepsilon' \circ \varphi^1 = \varepsilon$. By Proposition \ref{prop:DescentAbelianizationMaps} \ref{prop:DescentAbelianizationMaps1}, one has a commutative diagram with exact rows
	\begin{equation*}
		\xymatrix{
			& & & \Sha(G)  \ar[r]^-{\cong} \ar@{^{(}->}[d] & \Sha^1(k,C) \ar@{^{(}->}[d] \\
			& \H^0(k,C) \ar[r] \ar[d] & \H^1(k,G^{\sconnect}) \ar[r] \ar[d]^{\cong} & \H^1(k,G) \ar[r]^{\ab^1} \ar[d] & \H^1(k,C) \ar[d] \\
			G(\Abb_k) \ar[r]^-{\ab^0} & \Pbb^0(k,C) \ar[r] \ar[d] & \Pbb^1(k,G^{\sconnect})  \ar[r] & \Pbb^1(k,G) \ar[r]^{\ab^1} & \Pbb^1(k,C) \\
			& \H^1(k,\wh{C})^D  \ar[d] \\
			& \Sha^1(k,\wh{C})^D.
		}
	\end{equation*}
	Actually, the first column from the left is a complex, which is exact at the term $\H^1(k,\wh{C})^D $  \cite[Th\'eor\`eme 2.9]{Demarche2011Approximation}. The arrow on the second column is bijective (Kneser--Harder--Chernousov, see \cite[Theorems 6.4 and 6.6]{PR1994Group}). The other columns are obviously exact. Since $\varepsilon' \in \H^1(k,\wh{C})^D$ is mapped to $0 \in \Sha^1(k,\wh{C})^D$, it comes from an element $c_{\Abb} = (c_v)_{v \in \Omega_k} \in \Pbb^0(k,C)$. Its image $\gamma_{\Abb}$ in $\Pbb^1(k,G^{\sconnect})$ comes from a unique element $\gamma \in \H^1(k,G^{\sconnect})$, whose image $\delta$ in $\H^1(k,G)$ lies in $\Sha(G)$. Since $\ab^1(\delta) = 0$, we have $\delta = 1$ by Proposition \ref{prop:DescentAbelianizationMaps} \ref{prop:DescentAbelianizationMaps3}. It follows that $\gamma$ comes from an element $c \in \H^0(k,C)$. Its localization $\loc_\Omega(c) := (\loc_v(c))_{v \in \Omega_k} \in \Pbb^0(k,C)$ has the same image in $\Pbb^1(k,G^{\sconnect})$ as that of $c_\Abb$, hence there exists $g_\Abb = (g_v)_{v \in \Omega_k} \in G(\Abb_k)$ such that $\ab^0(g_v) = \loc_v(c) - c_v$ for all $v \in \Omega_k$. Since the first column is a complex, the image of $(\ab^0(g_v))_{v \in \Omega_k}$ in $\H^1(k,\wh{C})^D$ is the same as that of $-c_{\Abb}$, {\em i.e.}, $-\varepsilon'$. Let $y_\Abb \cdot g_\Abb:=(y_v \cdot g_v)_{v \in \Omega_k}$. By virtue Lemma \ref{lemm:DescentReductive1}, we have, for all $\alpha^1 \in \Br_{\nr}(Y)$ (noting that $\varepsilon = \varepsilon' \circ \varphi^1$),
	\begin{equation*}
		\pair{y_{\Abb} \cdot g_{\Abb}, \alpha^1}_{\BM} = \pair{y_{\Abb}, \alpha^1}_{\BM} + \sum_{v \in \Omega_k}\ab^0(g_v) \cup \varphi^1(\alpha^1) = \varepsilon(\alpha^1) - \varepsilon'(\varphi^1(\alpha^1)) = 0.
	\end{equation*}
	We have found a lifting $y_\Abb \cdot g_\Abb \in Y(\Abb_k)^{\Br_{\nr}(Y)}$ of $x_\Abb \in X(\Abb_k)^A$, and this concludes the proof of Theorem \ref{thm:DescentReductive}.
\end{proof}

\subsection{End of proof of the main theorem} \label{subsec:DescentConnectedFormal}

Let us now finish the proof of Theorem \ref{customthm:DescentConnected}, using Theorem \ref{thm:DescentReductive} from the previous subsection. The first step is to drop the reductive assumption.

\begin{thm} \label{thm:DescentNonReductive} 
    Let $G$ be a connected (not necessarily reductive) linear algebraic group over a number field $k$. Let $X$ be a smooth geometrically integral $k$-variety and $f: Y \to X$ a torsor under $G$. Let $A \subseteq \Br(X)$ denote the preimage of $\Br_{\nr}(Y)$ by the map $f^\ast: \Br(X) \to \Br(Y)$. Then
	\begin{equation*}
            X(\Abb_k)^A = \bigcup_{[\sigma] \in \H^1(k,G)} f^{\sigma}(Y^{\sigma}(\Abb_k)^{\Br_{\nr} (Y^{\sigma})}).
	\end{equation*}
\end{thm}
\begin{proof}
    Let $G^{\uni}$ be the unipotent radical of $G$ and $G^{\red} := G/G^{\uni}$ (which is connected reductive). Letting $Z:=Y/G^{\uni}$, we obtain a torsor $g: Z \to X$ under $G^{\red}$ and a torsor $h: Y \to Z$ under $G^{\uni}$. Now, \cite[Lemme 1.13]{Sansuc1981Arithmetique} tells us that $\H^1(K, (G^{\uni})^{\sigma}) = 1$ for any field extension $K/k$ and any cocycle $\sigma \in \Z^1(k,G)$, and that the map $\H^1(k,G) \to \H^1(k,G^{\red})$ is bijective. For every $[\tau] \in \H^1(k,G^{\red})$, there is a unique class $[\sigma] \in \H^1(k,G)$ lifting $[\tau]$, hence the torsor $f^{\sigma}: Y^{\sigma} \to X$ has a factorization $Y^{\sigma} \xrightarrow{h^\sigma} Z^\tau \xrightarrow{g^\tau} X$. Since $h^{\sigma}$ is a torsor under $(G^{\uni})^{\sigma}$ and $\H^1(k(Z^\tau), (G^{\uni})^{\sigma}) = 1$, it is a stable birational equivalence, hence induces an isomorphism $(h^{\sigma})^\ast: \Br_{\nr}(Z^\tau) \xrightarrow{\cong} \Br_{\nr}(Y^{\sigma})$. In particular, $A$ is equal to the preimage of $\Br_{\nr}(Z)$ by the map $g^\ast: \Br(X) \to \Br(Z)$. In view of Theorem \ref{thm:DescentReductive}, it remains to show that $h^{\sigma}(Y^{\sigma}(\Abb_k)^{\Br_{\nr}(Y^{\sigma})}) = Z^\tau(\Abb_k)^{\Br_{\nr}(Z^\tau)}$. The inclusion ``$\subseteq$'' is obvious. Conversely, let $z_\Abb \in Z^\tau(\Abb_k)^{\Br_{\nr}(Z^\tau)}$. Viewing that $\H^1(k_v,(G^{\uni})^{\sigma}) = 1$ for all $v \in \Omega_k$, we may lift $z_\Abb$ to a family $y_\Abb \in Y^{\sigma}(\Abb_k)$. Since $\Br_{\nr}(Z^\tau) \cong \Br_{\nr}(Y^{\sigma})$, we have $y_\Abb \in Y^{\sigma}(\Abb_k)^{\Br_{\nr}(Y^{\sigma})}$ by functoriality of the pairing \eqref{eq:BrauerManinPairingAdelic}. The theorem is proved.
\end{proof}

If $X$ is proper, then $\Br_{\nr}(X) = A = \Br(X)$, and the statement of Theorem \ref{thm:DescentNonReductive} becomes
    \begin{equation*}
        X(k_\Omega)^{\Br(X)} = \bigcup_{[\sigma] \in \H^1(k,G)} f^{\sigma}(Y^{\sigma}(\Abb_k)^{\Br_{\nr}(Y^{\sigma})}) = \bigcup_{[\sigma] \in \H^1(k,G)} f^{\sigma}(Y^{\sigma}(k_\Omega)^{\Br_{\nr}(Y^{\sigma})}),
    \end{equation*}
which implies Conjecture \ref{conj:Descent} for the torsor $f: Y \to X$. Under the rationally connected assumption, one can establish Theorem \ref{thm:DescentConnected} below ({\em i.e.}, Theorem \ref{customthm:DescentConnected}) without properness. To this end, it will be sufficient to deploy {\em Harari's formal lemma} \cite[Th\'eor\`eme 2.1.1, Corollaire 2.6.1]{Harari1994Fibration}.

\begin{lemm} \label{lemm:Formal}
    Let $X$ be a smooth geometrically integral variety over a field $k$ of characteristic $0$.
    \begin{enumerate}
        \item \label{lemm:Formal1} If $X$ is rationally connected, then $\Br_{\nr}(X)/\Br(k)$ is finite.
	
        \item \label{lemm:Formal2} Suppose that $k$ is a number field and $A \subseteq \Br(X)$ is a subgroup containing $\Br_{\nr}(X)$ such that $A/\Br(k)$ is finite. Then, one has $\ol{X(\Abb_k)^A} = X(k_\Omega)^{\Br_{\nr}(X)}$ in the topological space $X(k_\Omega)$ (both sides can be empty).
    \end{enumerate}
\end{lemm}
\begin{proof}
\begin{enumerate}
    \item Combine \cite[Lemma 1.3]{CTS2013Reduction} and \cite[Corollary 4.18]{Debarre2001Higher}.
    
    \item Let $X^c$ be a smooth compactification of $X$ (so that $\Br_{\nr}(X) = \Br(X^c)$) and apply \cite[Proposition 1.1]{CTS2000Fibration} (which is a consequence of Harari's formal lemma) to the inclusion $X \subseteq X^c$. This result says that $\ol{X(\Abb_k)^A} = X^c(k_\Omega)^{\Br_{\nr}(X)}$ in $X^c(k_\Omega)$, {\em a fortiori} $\ol{X(\Abb_k)^A} = X(k_\Omega)^{\Br_{\nr}(X)}$ in $X(k_\Omega)$. 
\end{enumerate}
\end{proof}

\begin{thm} [Theorem \ref{customthm:DescentConnected}] \label{thm:DescentConnected}
    Let $k$ be a number field, $X$ a smooth rationally connected $k$-variety, and $G$ a connected linear algebraic $k$-group. If  $f:Y \to X$ is a torsor under $G$, then
	\begin{equation*}
            X(k_\Omega)^{\Br_{\nr}(X)} = \ol{\bigcup_{[\sigma] \in \H^1(k,G)} f^{\sigma}(Y^{\sigma}(k_\Omega)^{\Br_{\nr}(Y^{\sigma})})}
	\end{equation*}
    in the topological space $X(k_\Omega)$. In particular, if the varieties $Y^{\sigma}$ satisfy $\BMHP$ (resp. $\BMWA$) for all $[\sigma] \in \H^1(k,G)$, then it is also the case for $X$.
\end{thm}
\begin{proof}
    The inclusion ``$\supseteq$'' is obvious by functoriality and continuity of the pairing \eqref{eq:BrauerManinPairingNr}, so let us show the inclusion ``$\subseteq$''. Let $A \subseteq \Br(X)$ denote the preimage of $\Br_{\nr}(Y)$ by the map $f^\ast: \Br(X) \to \Br(Y)$. By Theorem \ref{thm:DescentNonReductive}, we have
	\begin{equation*}
		X(\Abb_k)^A = \bigcup_{[\sigma] \in \H^1(k,G)} f^{\sigma}(Y^{\sigma}(\Abb_k)^{\Br_{\nr} (Y^{\sigma})}).
	\end{equation*}
    Thus, it remains to show that $X(k_\Omega)^{\Br_{\nr}(X)} = \ol{X(k_\Omega)^A}$. In view of Lemma \ref{lemm:Formal} \ref{lemm:Formal2}, we are reduced to show that the group $A/\Br(k)$ is finite. Since $X$ is rationally connected, it follows from the Graber--Harris--Starr theorem \cite{GHS2003Rational} that $Y$ is also rationally connected (see \cite[Corollaire 3.1]{Debarre2003Rational}). By Lemma \ref{lemm:Formal} \ref{lemm:Formal1}, the group $\Br_{\nr}(Y)/\Br(k)$ is finite. Thus, to establish the finiteness of $A/\Br(k)$, it is enough to show that $\Ker(\Br(X) \to \Br (Y))$ is finite. Using Sansuc's units-Picard-Brauer exact sequence \cite[(6.10.1)]{Sansuc1981Arithmetique}, we are reduced to prove the finiteness of $\Pic(G)$. The same sequence tells us that this property is stable under extensions of groups. Hence, we may deal separately with the following cases.

    \begin{itemize}
        \item If $G$ is unipotent, then $\Pic(G) = 0$ since the $k$-variety $G$ is isomorphic to an affine space.
    
        \item If $G$ is a torus, then $\Pic(G) = \H^1(k,\wh{G})$ by \cite[Lemme 6.9 (ii)]{Sansuc1981Arithmetique}. Let $K/k$ be a finite Galois extension splitting $G$, then $\H^1(K,\wh{G}) = 0$, so $\H^1(k,\wh{G}) = \H^1(\Gal(K/k),\wh{G})$ by the inflation-restriction sequence. This latter is torsion and finitely generated, hence finite. 
 
        \item If $G$ is semisimple, by \cite[Lemme 6.9 (iii)]{Sansuc1981Arithmetique}, $\Pic(G)$ injects into the character module of the algebraic fundamental group $\pi_1(\bar{G}) = \Ker(\bar{G}^{\sconnect} \to \bar{G})$, which is finite. 
    \end{itemize} 
\end{proof}
Having established Theorem \ref{thm:DescentConnected} for torsors, we can now generalize it to the situation where no torsors are given {\em a priori}. Let us bring the notion of nonabelian descent types from Subsection \ref{subsec:Type} into play.

\begin{thm} \label{thm:DescentConnectedType}
    Let $X$ be a smooth quasi-projective geometrically integral variety over a number field $k$. Let $L = (\bar{G},\kappa)$ be a $k$-kernel, where $\bar{G}$ is a connected linear algebraic $\bar{k}$-group, and $\lambda = [\bar{Y},E]$ a descent type on $X$ bound by $L$. If $X$ is rationally connected, then
	\begin{equation*}
		X(k_\Omega)^{\Br_{\nr}(X)} = \ol{\bigcup_{\substack{f: Y \to X \\ \dtype(Y) = \lambda}} f(Y(k_\Omega)^{\Br_{\nr}(Y)})}
	\end{equation*}
    in the topological space $X(k_\Omega)$. In particular, if the torsors $f: Y \to X$ of descent type $\lambda$ satisfy $\BMHP$ (resp. $\BMWA$), then it is also the case for $X$. By convention, this hypothesis holds vacuously if there are no such torsors at all.
\end{thm}
\begin{proof}
    By Lemma \ref{lemm:DescentType} \ref{lemm:DescentType3}, two $X$-torsors of descent type $\lambda$ differ by twisting by a Galois cycle. Hence, deducing the present theorem from Theorem \ref{thm:DescentConnected} amounts to showing that torsors of descent type $\lambda$ {\em exist} whenever $X(k_\Omega)^{\Br_{\nr}(X)} \neq \varnothing$. To this end, the condition $X(k_\Omega)^{\Be(X)} \neq \varnothing$ is in fact sufficient, where the subgroup $\Be(X) \subseteq \Br_{\nr}(X)$ was defined in \eqref{eq:FirstObstruction}. Indeed, this has been proved in \cite[Proposition 3.10]{Linh2025Type}.
\end{proof}

As an application of Theorem \ref{thm:DescentConnectedType}, we give a direct proof of Borovoi's theorem on the Brauer--Manin obstruction for homogeneous spaces with connected stabilizers \cite[Theorems 2.2 and 2.3]{Borovoi1996BrauerManin}, using Sansuc's result for connected linear algebraic groups.

\begin{coro} \label{coro:BorovoiThm} 
    Let $k$ be a number field. Then, the property $\BMHP$ and $\BMWA$ hold for homogeneous spaces of connected linear algebraic groups over $k$ with connected geometric stabilizers.
\end{coro}
\begin{proof}
    Let $G$ be a connected linear algebraic group over $k$ and let $X$ be a homogeneous space of $G$. Fix a $\bar{k}$-point of $X$ and let $\bar{H} \subseteq \bar{G}$ be its stabilizer, which we assume to be connected. As in \cite[Proposition 3.3]{HS2002Nonabelian}, the subgroup $E_X \subseteq \SAut(\bar{G}/X)$ consisting of semilinear $s$-automorphisms (for some $s \in \Gamma_k$) for which $\alpha(x \cdot g) = x \cdot \tensor[^s]{g}{}$ for all $x \in X(\bar{k})$ and $g \in G(\bar{k})$ fits in an exact sequence
        \begin{equation*}
            1 \to \bar{H}(\bar{k}) \to E_X \to \Gamma_k \to 1
        \end{equation*}
    of topological groups, which defines a $k$-kernel $L_X$ as well as a descent type $\lambda = [\bar{G},E_X]$ on $X$ bound by $L_X$. 
    The torsors $Y \to X$ of descent type $\lambda$ (if they exist) are left torsors under $G$ (they are $k$-forms of the variety $G$ twisted by cocycles with values in $\bar{H}(\bar{k})$). According to \cite[Corollaires 8.7 and 8.13]{Sansuc1981Arithmetique}, these torsors satisfy $\BMHP$ and $\BMWA$. In fact, Sansuc's original result concerning the property $\BMHP$ requires the absence of $E_8$ factors, since his proof used the result of Kneser--Harder (see Th\'eor\`eme 4.2 in {\it loc. cit.}) which excluded the $E_8$ case. This problem, however, was later resolved by Chernousov \cite{Chernousov1989E8} (see also \cite[Theorems 6.4 and 6.6]{PR1994Group}). It follows from Theorem \ref{thm:DescentConnectedType} that the properties $\BMHP$ and $\BMWA$ hold for $X$. 
\end{proof}

\subsection{Reduction to the finite case} \label{subsec:DescentConnectedReduction}

Next, let us show that Conjecture \ref{conj:Descent} can be reduced to the case of torsors under finite group schemes. Actually, instead of working directly with torsors, it is more convenient to replace them with descent types. 

In what follows, when $G$ is an algebraic group over $k$, we denote by $G^{\circ}$ its identity component, and $G^{\f}:=\pi_0(G) = G/G^{\circ}$. When $L = (\bar{G},\kappa)$ is a $k$-kernel, the outer Galois action $\kappa$ induces a canonical outer action on $\bar{G}^{\f}$ (since $\bar{G}^{\circ}$ is characteristic in $\bar{G}$). This gives rise to a $k$-kernel $L^{\f} = (\bar{G},\kappa^{\f})$ and a surjective morphism $\pi: L \to L^{\f}$ of $k$-kernels. If $X$ is a nonempty variety and $\lambda \in \DType(X,L)$, we shall write $\lambda^{\f}:=\pi_\ast \lambda \in \DType(X,L^{\f})$ ({\em cf.} Lemma \ref{lemm:DescentTypePushforward} \ref{lemm:DescentTypePushforward1}), which can then be thought as a finite descent type on $X$ in the sense of Harpaz--Wittenberg, {\em i.e.}, $\lambda^{\f}$ is completely determined by the $\bar{X}$-torsor $\bar{Z}:=\bar{Y}/\bar{G}^{\circ}$ under $\bar{G}^{\f}$ (see p. \pageref{FiniteDescentTYpe}).

\begin{thm} \label{thm:DescentFiniteType}
     Let $k$ be a number field, $L = (\bar{G},\kappa)$ a $k$-kernel with $\bar{G}$ linear (not necessarily connected), $X$ a smooth quasi-projective $k$-variety, and $\lambda = [\bar{Y},E] \in \DType(X,L)$ a descent type with $\bar{Y}$ rationally connected (in particular, $X$ is rationally connected). Then, one has
        \begin{equation} \label{eq:DescentFiniteType}
            \ol{\bigcup_{\substack{f: Y \to X \\ \dtype(Y) = \lambda}} f(Y(k_\Omega)^{\Br_{\nr}(Y)})} = \ol{\bigcup_{\substack{g: Z \to X \\ \dtype(Z) = \lambda^{\f}}} g(Z(k_\Omega)^{\Br_{\nr}(Z)})}
        \end{equation}
    in the topological space $X(k_\Omega)$. Here, we use the convention that both sides can be empty (for example, when there are no torsors of mentioned descent types at all). Furthermore, consider the following statements.
    \begin{enumerate}
        \item \label{thm:DescentFiniteType1} If the $X$-torsors of descent type $\lambda^{\f}$ satisfy $\BMHP$, then $X$ itself satisfies $\BMHP$.
        \item \label{thm:DescentFiniteType2} If the $X$-torsors of descent type $\lambda^{\f}$ satisfy $\BMWA$, then $X$ itself satisfies $\BMWA$.
        \item \label{thm:DescentFiniteType3} If the $X$-torsors of descent type $\lambda$ satisfy $\BMHP$, then $X$ itself satisfies $\BMHP$.
        \item \label{thm:DescentFiniteType4} If the $X$-torsors of descent type $\lambda$ satisfy $\BMWA$, then $X$ itself satisfies $\BMWA$.
    \end{enumerate} 
    Then, \ref{thm:DescentFiniteType1} implies \ref{thm:DescentFiniteType3} and \ref{thm:DescentFiniteType2} implies \ref{thm:DescentFiniteType4}.
\end{thm}
\begin{proof}
    The inclusion ``$\subseteq$'' of \eqref{eq:DescentFiniteType} follows from the following fact. If $Y \xrightarrow{f} X$ is a torsor under a $k$-form $G$ of $L$, of descent type $\lambda$, then $Z:=Y/G^{\circ} \xrightarrow{g} X$ is a torsor under the $k$-form $G^{\f}$ of $L^{\f}$, of descent type $\lambda^{\f}$ (see {\em e.g.} \cite[Construction 2.6]{Linh2025Type}). Conversely, let $x_\Omega \in X(k_\Omega)$ be a family lying in the right hand side of \eqref{eq:DescentFiniteType}. Then, for any neighborhood $\Uscr \subseteq X(k_\Omega)$ of $x_\Omega$, there exists a torsor $g: Z \to X$ of descent type $\lambda^{\f}$ such that $g(Z(k_\Omega)^{\Br_{\nr}(Z)}) \cap \Uscr \neq \varnothing$, or $Z(k_\Omega)^{\Br_{\nr}(Z)} \cap g^{-1}(\Uscr) \neq \varnothing$. By Lemma \ref{lemm:DescentTypePushforward} \ref{lemm:DescentTypePushforward2}, we obtain a $k$-kernel $L_Z = (\bar{G}^{\circ}, \kappa_Z)$ and a descent type $\lambda_Z = [\bar{Y},E_Z] \in \DType(Z,L_Z)$. By Theorem \ref{thm:DescentConnectedType}, there exists a torsor $h: Y \to Z$ of descent type $\lambda_Z$ such that $h(Y(k_\Omega)^{\Br_{\nr}(Y)}) \cap g^{-1}(\Uscr) \neq \varnothing$, or $g(h(Y(k_\Omega)^{\Br_{\nr}(Y)})) \cap \Uscr \neq \varnothing$. By Lemma \ref{lemm:DescentTypePushforward} \ref{lemm:DescentTypePushforward3}, the composite $g \circ h: Y \to X$ is a torsor of descent type $\lambda$. Hence, $x_\Omega$ lies in the left hand side of \eqref{eq:DescentFiniteType}. The identity \eqref{eq:DescentFiniteType} is thus proved.

    Now, we suppose that \ref{thm:DescentFiniteType1} holds true and proceed to show the contraposition of \ref{thm:DescentFiniteType3}. Assume that the property $\BMHP$ fails for $X$. Under the assumption that \ref{thm:DescentFiniteType1} holds true, there exists a torsor $Z \to X$ of descent type $\lambda^{\f}$ that fails the property $\BMHP$, that is, $Z(k_\Omega)^{\Br_{\nr}(Z)} \neq \varnothing$ but $Z(k) = \varnothing$. In particular, the right hand side of \eqref{eq:DescentFiniteType} is nonempty. Hence, the left hand side is also nonempty, that is, there exists a torsor $Y \to X$ of descent type $\lambda$ such that $Y(k_\Omega)^{\Br_{\nr}(Y)} \neq \varnothing$. But, since $X(k) = \varnothing$, we have necessarily that $Y(k) = \varnothing$, {\em i.e.}, $Y$ fails the property $\BMHP$. This means \ref{thm:DescentFiniteType3} holds true.
    
    Finally, we suppose that \ref{thm:DescentFiniteType2} holds true and proceed to show \ref{thm:DescentFiniteType4} directly. Assume that every torsor $Y \to X$ of descent type $\lambda$ satisfies $\BMWA$. Let $g: Z \to X$ be any torsor of descent type $\lambda^{\f}$. Assume that there is a family $z_\Omega \in Z(k_\Omega)^{\Br_{\nr}(Z)}$, and let $\Uscr \subseteq Z(k_\Omega)$ be any neighborhood of it. Again, by Lemma \ref{lemm:DescentTypePushforward} \ref{lemm:DescentTypePushforward2} and Theorem \ref{thm:DescentConnectedType}, we obtain a torsor $h: Y \to Z$ under a $k$-form of $\bar{G}^{\circ}$ such that $h(Y(k_\Omega)^{\Br_{\nr}(Y)}) \cap \Uscr \neq \varnothing$, and the composite $g \circ h: Y \to X$ is a torsor of descent type $\lambda$. Under our assumption, one has $Y(k) \cap h^{-1}(\Uscr) \neq \varnothing$, {\em a fortiori} $Z(k) \cap \Uscr \neq \varnothing$. This means $Z$ satisfies $\BMWA$. Under the assumption that \ref{thm:DescentFiniteType2} holds true, the variety $X$ also satisfies $\BMWA$, which proves \ref{thm:DescentFiniteType4}.
\end{proof}

In Theorem \ref{thm:DescentFiniteType}, the quasi-projective assumption is required since we were working with descent types instead of torsors. If a torsor $Y \to X$ is given, it is possible to drop this condition, allowing one to recover the following.

\begin{thm} [Theorem \ref{customthm:DescentFinite}] \label{thm:DescentFinite}
    Let $k$ be a number field, $G$ a linear (not necessarily connected) algebraic $k$-group, $X$ be a smooth $k$-variety, and $f: Y \to X$ a torsor under $G$, with $Y$ rationally connected (in particular, $X$ is rationally connected). Put $Z:=Y/G^{\circ}$, which gives a torsor $g: Z \to X$ under $G^{\f}$ and a torsor $h: Y \to Z$ under $G^{\circ}$. Then, one has
        \begin{equation} \label{eq:DescentFinite}
            \ol{\bigcup_{[\sigma] \in \H^1(k,G)} f^{\sigma}(Y^{\sigma}(k_\Omega)^{\Br_{\nr}(Y^{\sigma})})} = \ol{\bigcup_{[\tau] \in \H^1(k,G^{\f})}  g^{\tau}(Z^{\tau}(k_\Omega)^{\Br_{\nr}(Z^{\tau})})}
        \end{equation}
    in the topological space $X(k_\Omega)$. Furthermore, consider the following statements.
    \begin{enumerate}
        \item \label{thm:DescentFinite1} If $Z^\tau$ satisfies $\BMHP$ for every $\tau \in \Z^1(k,G^{\f})$, then $X$ itself satisfies $\BMHP$.
        \item \label{thm:DescentFinite2} If $Z^\tau$ satisfies $\BMWA$ for every $\tau \in \Z^1(k,G^{\f})$, then $X$ itself satisfies $\BMWA$.
        \item \label{thm:DescentFinite3} If $Y^\sigma$ satisfies $\BMHP$ for every $\sigma \in \Z^1(k,G)$, then $X$ itself satisfies $\BMHP$.
        \item \label{thm:DescentFinite4} If $Y^\sigma$ satisfies $\BMWA$ for every $\sigma \in \Z^1(k,G)$, then $X$ itself satisfies $\BMWA$.
    \end{enumerate} 
    Then, \ref{thm:DescentFinite1} implies \ref{thm:DescentFinite3} and \ref{thm:DescentFinite2} implies \ref{thm:DescentFinite4}.
\end{thm}
\begin{proof}
    If $X$ is quasi-projective, then the theorem follows from Theorem \ref{thm:DescentFiniteType} and Lemma \ref{lemm:DescentType} \ref{lemm:DescentType3}. In the general case, take any dense affine open subset $U \subseteq X$, and let $V:=f^{-1}(U)$ and $W:=g^{-1}(U)$. Then, the identity \eqref{eq:DescentFiniteType} from Theorem \ref{thm:DescentFiniteType} becomes
        \begin{equation*}
            \ol{\bigcup_{[\sigma] \in \H^1(k,G)} f(V^{\sigma}(k_\Omega)^{\Br_{\nr}(V^{\sigma})})} = \ol{\bigcup_{[\tau] \in \H^1(k,G^{\f})}  g(W^{\tau}(k_\Omega)^{\Br_{\nr}(W^{\tau})})},
        \end{equation*}
    which implies \eqref{eq:DescentFinite}, because $V^{\sigma}(k_\Omega)^{\Br_{\nr}(V^{\sigma})}$ (resp. $W^{\tau}(k_\Omega)^{\Br_{\nr}(W^{\tau})}$) is dense in $Y^{\sigma}(k_\Omega)^{\Br_{\nr}(Y^{\sigma})}$ (resp. $Z^{\tau}(k_\Omega)^{\Br_{\nr}(Z^{\tau})}$) for all $\sigma \in \Z^1(k,G)$ (resp. $\tau \in \Z^1(k,G^{\f})$), thanks to Lemma \ref{lemm:Formal}. The implications ``\ref{thm:DescentFinite1} $\Rightarrow$ \ref{thm:DescentFinite3}'' and ``\ref{thm:DescentFinite2} $\Rightarrow$ \ref{thm:DescentFinite4}'' follow from the respective implications in Theorem \ref{thm:DescentFiniteType} (applied to $U$) and the birational invariance of the properties $\BMHP$ and $\BMWA$.
\end{proof}

The ``finite case'' of Conjecture \ref{conj:Descent} is hence the most difficult one. The inductive argument used in the proofs of Theorems \ref{thm:DescentFiniteType} and \ref{thm:DescentFinite} was used in a similar way by Harpaz and Wittenberg in their work \cite{HW2024Supersolvable} on supersolvable descent. The base case ({\em i.e.}, ``cyclic descent'') relies on the Graber--Harris--Starr theorem \cite{GHS2003Rational}. If we were able to prove Conjecture \ref{conj:Descent} (for torsors) in the finite abelian case, then we would obtain the same result for descent types thanks to \cite[Lemma 3.7]{HW2024Supersolvable}. Using the same sort of d\'evissage argument, the descent conjecture for the descent types whose underlying $\bar{k}$-group are finite solvable (as an abstract group!) would then be proved. In particular, this would imply that homogeneous spaces of connected linear algebraic groups with finite solvable geometric stabilizers have the properties $\BMHP$ and $\BMWA$ (it suffices to repeat the proof of Corollary \ref{coro:BorovoiThm}). A direct consequence of this would be Shafarevich's theorem on solvable Galois groups. Unfortunately, at least to the author's knowledge, a proof of Conjecture \ref{conj:Descent} in the finite abelian case is currently out of reach. It should be noted that if $G$ is a connected linear algebraic group over a number field $k$ such that $G^{\ssimple} = G^{\sconnect}$ and $X$ is a homogeneous space of $G$ with finite abelian geometric stabilizers, then $X$ satisfies $\BMHP$ and $\BMWA$ (see \cite[Theorems 2.2 and 2.3]{Borovoi1996BrauerManin}).

Just like in the cyclic case, where Harpaz and Wittenberg applied the sophisticated theorem of Graber--Harris--Starr, a proof of the descent conjecture for $G$ finite abelian should use the rational connected assumption in an essential way. As a matter of fact, this assumption is {\em absolutely indispensable}, even when $G = \Zbb/2$. Indeed, in \cite[Theorem 2.10, Proposition 2.11]{HS2005Enriques}, Harari and Skorobogatov constructed a double cover $f: Y \to X$, where $X$ is an Enriques surfaces $X$ over $k = \Qbb$, and $Y$ is a Kummer surface, such that $X(k_\Omega)^{\Br(X)} \not \subseteq \bigcup_{[\sigma] \in \H^1(k,\Zbb/2)} f^{\sigma}(Y^{\sigma}(k_\Omega))$ (by properness, this latter is a closed subset of $X(k_\Omega)$). Worse, the work of Balestrieri--Berg--Manes--Park--Viray \cite[Theorem 1.2]{BBMPV2016Enriques} gave an Enriques surface $X$ over $k=\Qbb$, and a double cover $f: Y \to X$, where $Y$ is a $K3$ surface, such that $X(k) = Y^{\sigma}(k_\Omega)^{\Br(Y^{\sigma})} = \varnothing$ for all $\sigma \in \H^1(k,\Zbb/2)$, but $X(k_\Omega)^{\Br(X)} \neq \varnothing$ (that is, the descent varieties $Y^{\sigma}$ vacuously satisfy $\BMHP$ and $\BMWA$, but $X$ does not).

\section{The descent conjecture for zero-cycles} \label{sec:DescentCycle}

In this section, we prove Theorem \ref{customthm:DescentCycle}. For convenience, all of the subsequent results shall be formulated in the language of descent types. They can be easily converted into ones for torsors (dropping the quasi-projective assumption) using the same trick deployed in the proof of Theorem \ref{thm:DescentFinite}.

\subsection{Brauer--Manin obstruction for zero-cycles} \label{subsec:DescentCycleBM}

We recall some facts on the Chow groups of zero-cycles.

{\bf Zero-cycles.}  The standard reference for the Chow groups is \cite{Fulton1984Intersection}. When $X$ is a {\em smooth} variety over a field $k$, denote by $\Z_0(X)$ the group of zero-cycles on $X$, which is functorial contravariantly with respect to flat morphisms of relative dimension $0$, and covariantly with respect to any morphism. Its quotient $\CH_0(X)$ by rational equivalence is the {\em zero-dimensional Chow group} of $X$, which is functorial covariantly with respect to proper morphisms, and contravariantly with respect to generically finite morphisms (between smooth varieties). For any field extension (not necessarily finite) $K/k$, we have a restriction homomorphism $\res_{K/k}: \Z_0(X) \to \Z_0(X_K)$, which preserves degree. If $[K:k]$ is finite, we also have a corestriction homomorphism $\N_{K/k}: \Z_0(X_K) \to \Z_0(X)$ (also called the {\em $(K/k)$-norm}) which multiplies the degree by $[K:k]$. Both of these maps induce homomorphisms on Chow groups, and the composite $\N_{K/k} \circ \res_{K/k}: \Z_0(X) \to \Z_0(X)$ is the multiplication by $[K:k]$. 

{\bf Local (co-)restrictions over number fields.} If $k$ is a number field, put $\Z_{0,\Omega}(X):=\prod_{v \in \Omega_k}\Z_0(X_v)$. For $v \in \Omega_k$, denote by $\loc_v: \Z_0(X) \to \Z_0(X_v)$ the restriction map $\res_{k_v/k}$, and let $\loc_\Omega = ((\loc_v)_{v \in \Omega_k}): \Z_0(X) \to \Z_{0,\Omega}(X)$. The {\em local norm}
    \begin{equation} \label{eq:DescentCycleLocalNorm}
        \N_{K/k}^{\Omega}: \Z_{0,\Omega}(X_{K}) \to \Z_{0,\Omega}(X)
    \end{equation}
is defined as follows. For each family $x_\Omega = (x_w)_{w \in \Omega_K} \in \Z_{0,\Omega}(X_K)$, we define $\N_{K/k} (x_\Omega) \in \Z_{0,\Omega}(X)$ to be the family whose $v$-component, for each $v \in \Omega_k$, is $(\N_{K/k}^{\Omega}(x_\Omega))_v:=\sum_{w \in \Omega_K,\,w|v}\N_{K_w/k_v}(x_w)$. One checks with ease that the construction \eqref{eq:DescentCycleLocalNorm} is functorial and compatible with the global norm $\N_{K/k}$. Similarly, we define the {\em local restriction}
        \begin{equation} \label{eq:DescentCycleLocalRestriction}
            \res_{K/k}^{\Omega}: \Z_{0,\Omega}(X) \to \Z_{0,\Omega}(X_K)
	\end{equation}
as follows. For each family $x_\Omega = (x_v)_{v \in \Omega_k} \in Z_{0,\Omega}(X)$ and each place $w \in \Omega_K$ lying over a place $v \in \Omega_k$, let $(\res_{K/k}^{\Omega}(x_\Omega))_w:=\res_{K_w/k_v}(x_v)$. Then, \eqref{eq:DescentCycleLocalRestriction} is functorial and compatible with the global restriction $\res_{K/k}$. Furthermore, $\N_{K/k}^{\Omega} \circ \res_{K/k}^{\Omega}$ is the multiplication by $[K:k]$ on $\Z_{0,\Omega}(X)$.

{\bf Brauer--Manin obstruction.} A complete exposition on the Brauer--Manin obstruction in the context of zero-cycles can be found in \cite[\S{15.1}]{CTS2021Brauer}. Let $X$ be a smooth geometrically integral variety over a number field $k$. We have the local pairings
    \begin{equation*} \label{eq:ZeroCyclePairingLocal}
		\pair{-,-}_v: \Z_0(X_v) \times \Br(X) \to \Qbb/\Zbb,
    \end{equation*}
defined for all $v \in \Omega_k$, which factors through a pairing $\CH_0(X_v) \times \Br(X) \to \Qbb/\Zbb$ whenever $X$ is proper. If $v \in \Omega_k^{\infty}$, this pairing vanishes on the image of the norm $\N_{\bar{k}_v/k_v}: \CH_0(\bar{X}_v) \to \CH_0(X_v)$. By restricting to the unramified elements of $\Br(X)$, we obtain from \eqref{eq:ZeroCyclePairingLocal} a pairing
    \begin{equation} \label{eq:BrauerManinPairingCycle} \tag{BM${}_{\text{cyc}}$}
        \pair{-,-}_{\BM}: \Z_{0,\Omega}(X) \times \Br_{\nr}(X) \to \Qbb/\Zbb, \quad ((x_v)_{v \in \Omega_k}, \alpha) \mapsto \sum_{v \in \Omega_k}\pair{x_v,\alpha}_v.
    \end{equation}
The pairing \eqref{eq:BrauerManinPairingCycle} vanishes on the image of $\loc_\Omega$ and on $\Img(\Br(k) \to \Br(X))$, thanks to the global reciprocity law. When $X$ is proper, we let
    \begin{equation*}
        \CH_0'(X_v):= \begin{cases}
            \CH_0(X_v) & \text{if } v \in \Omega_k^{\f}, \\
            \Coker(\N_{\bar{k}_v/k_v}) & \text{if } v \in \Omega_k^{\infty}.
        \end{cases}
    \end{equation*}
Put $\CH_{0,\Omega}:=\prod_{v \in \Omega_k} \CH_0'(X_v)$. Then, the pairing \eqref{eq:BrauerManinPairingCycle} induces a pairing
	\begin{equation*}
		\CH_{0,\Omega}(X) \times \Br(X) \to \Qbb/\Zbb,
	\end{equation*}
which vanishes on the image of $\loc_\Omega$, hence a map $\CH_{0,\Omega}(X) \to \Br(X)^D$ fitting in a complex
	\begin{equation} \label{eq:ConjectureEPrime} \tag{E$'$}
		\CH_0(X) \xrightarrow{\loc_\Omega} \CH_{0,\Omega}(X) \to \Br(X)^D.
	\end{equation}
Applying the functor $A \mapsto \varprojlim_n(A/n)$ to \eqref{eq:ConjectureEPrime} (noting that the group $\Br(X)$ is torsion) yields the complex \eqref{eq:ConjectureE}. Recall that the conjecture (E) ({\em i.e.}, Conjecture \ref{conj:E}) predicts that \eqref{eq:ConjectureE} is exact for any smooth proper geometrically integral variety $X$. Note that the exactness of \eqref{eq:ConjectureE} is a stable birational invariant of smooth proper geometrically integral varieties, because the groups $\CH_0(X)$ and $\Br(X)$ are so (see \cite[Remarques 1.1 (vi)]{Wittenberg2012Fibration}, \cite[Example 16.1.11]{Fulton1984Intersection}, and \cite[Proposition 6.3]{CTC1979Rational}). 

{\bf Weak topology on $\Z_{0,\Omega}(X)$.} \label{TopologyZeroCycles} Following Liang \cite[\S 1.2]{Liang2013Arithmetic}, we consider the problem of weak approximation for zero-cycles. Let $X$ be a (not necessarily proper) smooth geometrically integral variety defined over a number field $k$. Let $X^c$ be a smooth compactification of $X$. For each $v \in \Omega_k$, we equip the group $\Z_0(X_{v})$ with the initial topology defined by the composites
    \begin{equation*}
        \Z_0(X_{v}) \subseteq \Z_0(X_{v}^c) \to \CH_0'(X_{v}^c) \to \CH_0'(X_{v}^c)/n,
    \end{equation*} 
for $n \ge 1$, where $\CH_0'(X_{v}^c)/n$ is equipped with the discrete topology. By birational invariance of the group $\CH_0$, this topology is independent of the choice of the smooth compactification $X^c$. If $K/k_v$ is any finite extension and $P, P' \in X(K)$ are two points which are close to each other, then $[P]$ and $[P']$ are close in $\Z_0(X_v)$ \cite[Lemme 1.8]{Wittenberg2012Fibration}. Finally, we equip the group $\Z_{0,\Omega}(X) = \prod_{v \in \Omega_k} \Z_0(X_{v})$ with the product topology. Hence, two families $x_{\Omega} = (x_v)_{v \in \Omega_k}$ and $x'_{\Omega} = (x'_v)_{v \in \Omega_k}$ are close to each other if for a large (but finite) subset $S \subseteq \Omega_k$ and a large integer $n \ge 1$, the local zero-cycles $x_v$ and $x'_v$ have the same image in $\CH'_0(X^c_{k_v})/n$ for every $v \in S$.

{\bf Descent.} In \cite{BB2024Descent}, Balestrieri and Berg developed a formalism allowing one to study descent theory in the context of zero-cycles. In short, given a torsor $f: Y \to X$ under an algebraic group $G$ over $k$, we consider all finite extension $K/k$ and all the twisted torsors $Y^{\sigma}_K \to K$, where $[\sigma] \in \H^1(K,G)$. Since the notion of descent type (introduced in Subsection \ref{subsec:Type}) classifies torsors up to twisting by a Galois cocycle, it is more convenient to translate the formalism of Balestrieri--Berg into one in this new language. Let $\lambda$ be a descent type on $X$ bound by a $k$-kernel $L$. We fix an algebraic closure $\bar{k}$ of $k$, and let $\Tscr_\lambda$ be the set of pairs $(K,[Y])$, where $K \subseteq \bar{k}$ is a finite extension of $k$, and $[Y]$ is the class of a torsor $Y \to X_{K}$ of descent type $\dtype(Y) = \lambda_K:=\res_{K/k}(\lambda) \in \DType(X_K,L_K)$. Following \cite[Definition 2.2]{BB2024Descent}, we define the ``recombining'' map associated with $\lambda$ as
\begin{equation} \label{eq:DescentCycleRecombinationGlobal}
	\rec_\lambda = \sum_{(K,[f: Y \to X_{K}]) \in \Tscr_\lambda} \N_{K/k} \circ f_\ast: \bigoplus_{(K,[Y]) \in \Tscr_\lambda} \Z_0(Y) \to \Z_0(X).
\end{equation}
On the local level, one defines
\begin{equation} \label{eq:DescentCycleRecombinationLocal}
	\rec_{\lambda,\Omega} = \sum_{(K,[f: Y \to X_{K}]) \in \Tscr_\lambda} \N_{K/k}^{\Omega} \circ f_\ast: \bigoplus_{(K,[Y]) \in \Tscr_\lambda} \Z_{0,\Omega}(Y) \to \Z_{0,\Omega}(X),
\end{equation}
where $\N_{K/k}^{\Omega}: \Z_{0,\Omega}(X_K) \to \Z_{0,\Omega}(X)$ is the local norm map \eqref{eq:DescentCycleLocalNorm}. Denote by $\Z_{0,\Omega}(X)^{\Br_{\nr}(X)}$ the subgroup of $\Z_{0,\Omega}(X)$ consisting of the families orthogonal to $\Br_{\nr}(X)$ relative to the pairing \eqref{eq:BrauerManinPairingCycle}. By functoriality of this pairing, \eqref{eq:DescentCycleRecombinationLocal} restricts to a homomorphism
\begin{equation} \label{eq:DescentCycleRecombinationBrauer}
	\rec_{\lambda,\Omega}^{\Br}: \bigoplus_{(K,[Y]) \in \Tscr_\lambda} \Z_{0,\Omega}(Y)^{\Br_{\nr}(Y)} \to \Z_{0,\Omega}(X)^{\Br_{\nr}(X)}.
\end{equation}

Here is the main result of this section.

\begin{thm} [Theorem \ref{customthm:DescentCycle}]  \label{thm:DescentCycle}
    Let $X$ be a smooth quasi-projective variety over a number field $k$. Let $\lambda = [\bar{Y},E]$ be a descent type on $X$ bound by a $k$-kernel $L = (\bar{G},\kappa)$, with $\bar{Y}$ rationally connected (the existence of such a descent type imply that $X$ is rationally connected) and $\bar{G}$ linear. For any smooth compactification $X^c$ of $X$, one has
		\begin{equation*}
			\Z_{0,\Omega}(X)^{\Br_{\nr}(X)} = \ol{\loc_\Omega(\Z_0(X^c)) + \Img(\rec_{\lambda,\Omega}^{\Br})}
		\end{equation*}
    in the set $\Z_{0,\Omega}(X)$ equipped with the topology defined in p. \pageref{TopologyZeroCycles}. In particular, if for every finite extension $K/k$ and every torsor $Y \to X_{K}$ of descent type $\dtype(Y) = \lambda_K$, Conjecture \ref{conj:E} holds for the smooth compactifications of $Y$, then it is also the case for the smooth compactifications of $X$.
\end{thm}

The strategy for the proof of Theorem \ref{thm:DescentCycle} is as follows. In Subsection \ref{subsec:DescentCycleConnected}, we combine Liang's method (passing from rational points to zero-cycles) in \cite[\S 3.2]{Liang2013Arithmetic} with Theorem \ref{thm:DescentConnectedType} to deal with the case where $\bar{G}$ is connected. In particular, this can be applied to the case where $\bar{G} = \bar{T}$ is a torus. Next, using the fibration argument of Harpaz--Wittenberg in \cite[\S 3.1 and 3.2]{HW2024Supersolvable} for rational points to deal with the case where $\bar{G}$ is a group of multiplicative type in Subsection \ref{subsec:DescentCycleMultiplicative} (this requires the previous case of torsors under a torus). In particular, this implies the result in the case where $\bar{G}$ is finite abelian. An inductive argument then yields the result in the case where $\bar{G}$ is finite solvable. The case where $\bar{G}$ is finite is treated in Subsection \ref{subsec:DescentCycleGeneral}, by passing to its Sylow subgroups (this uses an argument similar to one in \cite[\S 7.2]{HW2020Galois}). Finally, we deploy (again) a d\'evissage argument to combine the results in the finite and connected cases.

\subsection{The connected case} \label{subsec:DescentCycleConnected}

In this subsection, we prove the following result.

\begin{prop} \label{prop:DescentCycleConnected}
    Let $X$ be a smooth quasi-projective rationally connected variety over a number field $k$. Let $\lambda$ be a descent type on $X$ bound by a $k$-kernel $L = (\bar{G},\kappa)$, with $\bar{G}$ linear and connected. Then, one has
	\begin{equation*}
		\Z_{0,\Omega}(X)^{\Br_{\nr}(X)} = \ol{\bigcup_{(K, [f: Y \to X_K]) \in \Tscr_\lambda} \N_{K/k}^{\Omega}(f_\ast (Y(K_\Omega)^{\Br_{\nr}(Y_K)}))},
	\end{equation*}
	where $\Z_{0,\Omega}(X)$ is equipped with the topology defined in p. \pageref{TopologyZeroCycles}.
\end{prop}

First, we show the following.

\begin{lemm} \label{lemm:DescentCycleConnected}
    Let $X$ be a smooth {\em proper} geometrically integral variety over a number field $k$ such that the N\'eron--Severi group of $\bar{X}$ is torsion-free and $\H^1(\bar{X},\Ocal_{\bar{X}}) = \H^2(\bar{X},\Ocal_{\bar{X}}) = 0$ (this is the case if $X$ is rationally connected). For any $x_\Omega = (x_v)_{v \in \Omega_k} \in \Z_{0,\Omega}(X)^{\Br(X)}$, any finite set $S \subseteq \Omega_k$ of places, and any integer $n \ge 1$, there exist a finite extension $K/k$ and a family $P_\Omega = (P_w)_{w \in \Omega_K} \in X_K(K_\Omega)^{\Br(X_K)}$ such that $(\N_{K/k}^{\Omega}([P_\Omega]))_v$ and $x_v$ have the same image in $\CH_0'(X)/n$ for all  $v \in S$.
\end{lemm}
\begin{proof}
    Note that, for (smooth proper) rationally connected varieties, the complex \eqref{eq:ConjectureEPrime} is exact if and only if \eqref{eq:ConjectureE} is so (see \cite[Rappel 1.3]{HW2020Galois}). For $X = \Spec(k)$, the complex \eqref{eq:ConjectureE} is exact because it is the dual of the Albert--Brauer--Noether--Hasse sequence from global class field theory:
        \begin{equation*}
            0 \to \Br(k) \to \bigoplus_{v \in \Omega_k} \Br(k_v) \xrightarrow{\sum} \Qbb/\Zbb \to 0.
        \end{equation*}
    It follows that \eqref{eq:ConjectureEPrime} is also exact for $X = \Spec(k)$. Since $\CH_0'(\Spec(k_v)) = \Zbb$ for $v \in \Omega_k^{\f}$ and $\CH_0'(\Spec(\Rbb)) = \Zbb/2$, any family $(x_v)_{v \in \Omega_k} \in \Z_{0,\Omega}(X)^{\Br(X)}$ has ``the same'' local degree in the following sense. There exists an integer $d \in \Zbb$ such that $\deg_{k_v}(x_v) = d$ for all $v \in \Omega_k^{\f}$ and $\deg_{k_v}(x_v) \equiv d \pmod 2$ for all real places $v \in \Omega_k^{\infty}$. For $v \in \Omega_k^{\infty}$, we take any point $\bar{P}_v \in X(\bar{k}_v)$ and add a suitable multiple of the $0$-cycle $\N_{\bar{k}_v/k_v}(\bar{P}_v)$ to $x_v$ (which does not affect the class of $x_v$ in $\CH_0'(X_{k_v})$). Hence, we may assume that $\deg_{k_v}(x_v) = d$ for all $v \in \Omega_k$.
    
    The rest of the proof is essentially done in \cite[Proof of Theorem 3.2.1]{Liang2013Arithmetic} if $d = 1$, and in {\em loc. cit.}, Proposition 3.4.1 if $d \neq 1$. At the end of these proofs, we obtain a finite extension $K/k$ and a family $P_\Omega = (P_w)_{w \in \Omega_K} \in X_K(K_\Omega)^{\Br(X_K)}$ satisfying the stated condition.
\end{proof}

\begin{proof} [Proof of Proposition \ref{prop:DescentCycleConnected}]
    Let $X^c$ be a smooth compactification of $X$ and let $x_\Omega \in \Z_{0,\Omega}(X)^{\Br(X^c)}$. By Lemma \ref{lemm:DescentCycleConnected}, we find a finite extension $K/k$ and a family $P_\Omega \in X_K^c(K_\Omega)^{\Br(X^c)}$ such that $\N_{K/k}^{\Omega}([P_\Omega]) \in \Z_{0,\Omega}(X^c)^{\Br(X^c)}$ is sufficiently close to $x_\Omega$. By Lemma \ref{lemm:Formal}, we may assume that $P_\Omega \in X_K(K_\Omega)^{\Br(X^c)}$. By Theorem \ref{thm:DescentConnectedType}, there exists a torsor $f: Y \to X_K$ of descent type $\lambda_K$ and a family $Q_\Omega \in Y(K_\Omega)^{\Br_{\nr}(Y)}$ such that $f(Q_\Omega)$ is sufficiently close to $P_\Omega$, hence $\N_{K/k}^{\Omega}(f_\ast[Q_\Omega])$ is sufficiently close to $x_\Omega$.
\end{proof}

\subsection{The multiplicative type case} \label{subsec:DescentCycleMultiplicative}

In this subsection, we adapt the argument from \cite[\S 3.2]{HW2024Supersolvable} to show the following.
 
\begin{prop} \label{prop:DescentCycleMultiplicative}
    Let $X$ be a smooth quasi-projective variety over a number field $k$. Let $M$ be a $k$-group of multiplicative type and let $\lambda = [\bar{Y},E] \in \DType(X,M)$, with $\bar{Y}$ rationally connected (in particular, $X$ is rationally connected). Then, for any smooth compactification $X^c$ of $X$, one has
	\begin{equation*}
		\Z_{0,\Omega}(X)^{\Br_{\nr}(X)} = \ol{\loc_\Omega(\Z_0(X^c)) + \bigcup_{(K, [f: Y \to X_K]) \in \Tscr_\lambda} \N_{K/k}^{\Omega}(f_\ast (Y(K_\Omega)^{\Br_{\nr}(Y_K)}))}
	\end{equation*}
    in the set $\Z_{0,\Omega}(X)$ equipped with the topology defined in p. \pageref{TopologyZeroCycles}.
\end{prop} 

We require the following analogue of \cite[Theorem 3.6]{HW2024Supersolvable} in the context of zero-cycles.

\begin{lemm} \label{lemm:DescentCycleMultiplicative}
    Let $Z$ be a smooth geometrically integral variety equipped with a dominant morphism $\pi: Z \to \Abb^n_k$, for some $n \ge 1$. Let $U \subseteq \Abb^n_k$ be a dense open subset over which $\pi$ is smooth, and suppose that $Z_q = \pi^{-1}(q)$ is rationally connected for all $q \in U$. Then, for any smooth compactification $Z^c$ of $Z$ and any family $z_{\Omega} \in \Z_{0,\Omega}(Z)^{\Br(Z^c)}$, there exist a global zero-cycle $z \in \Z_0(Z^c)$, a closed point $q \in U$, and a family $z'_\Omega \in \Z_{0,\Omega}(Z_q)^{\Br_{\nr}(Z_q)}$ such that $\loc_\Omega(z) + \N_{k(q)/k}(z'_\Omega)$ is sufficiently close to $z_\Omega$ in $\Z_{0,\Omega}(Z^c)$ (equipped with the topology defined in p. \pageref{TopologyZeroCycles}).
\end{lemm}
\begin{proof}
    Since the group $\CH_0$ is a birational invariant of smooth proper geometrically integral variety, it is enough to prove the lemma for an {\em arbitrarily chosen} smooth compactification $Z^c$. Hence, by Nagata's theorem, we may assume that $\pi$ extends into a proper morphism $\pi^c: Z^c \to \Pbb^n_k$, which is smooth over $U$ (after shrinking $U$ if necessary). 
	
    We proceed by induction on $n$. For $n = 1$, this follows from \cite[Proposition 8.7]{HW2016Fibration} (note that  $\Z_{0,\Omega}(Z_q)^{\Br(Z_q^c)}$ is dense in $\Z_{0,\Omega}(Z_q^c)^{\Br(Z_q^c)}$ by continuity of \eqref{eq:BrauerManinPairingCycle} and Lemma \ref{lemm:Formal} \ref{lemm:Formal1}). 
	
    Let $n \ge 2$ and let $p: \Abb^n_k \to \Abb^1_k$ be the first projection. Let $V \subseteq \Abb^1_k$ be a dense open subset over which $p \circ \pi$ is smooth, such that $U \cap p^{-1}(r) \neq \varnothing$ for every $r \in V$. For such $r$, the fibre $Z_{r} = (p \circ \pi)^{-1}(r)$ is equipped with a dominant morphism to $p^{-1}(r) \cong \Abb^{n-1}_k$ with rationally connected generic fibre, hence is rationally connected by the theorem of Graber--Harris--Starr (see  \cite[Corollaire 3.1]{Debarre2003Rational}). Given a family $z_{\Omega} \in \Z_{0,\Omega}(Z)^{\Br(Z^c)}$, by the case $n = 1$ applied to $p \circ \pi$, we find a global zero-cycle $w \in \Z_0(Z^c)$, a closed point $r \in V$, and a family $w'_\Omega \in \Z_{0,\Omega}(Z_{r})^{\Br(Z_r^c)}$ such that $\loc_\Omega(w) + \N_{k(r)/k} (w'_\Omega)$ is sufficiently close to $z_\Omega$. By the induction hypothesis applied to the restriction $\pi|_{Z_r}: Z_r \to p^{-1}(r)$, we obtain a global zero-cycle $w^1 \in \Z_0(Z^c_r)$, a closed point $q \in U \cap \pi^{-1}(r)$, and a family $z'_\Omega \in \Z_{0,\Omega}(Z_q)^{\Br(Z_q^c)}$ such that $\loc_\Omega(w^1) + \N_{k(q)/k(r)} (z'_\Omega)$ is sufficiently close to $w'_\Omega$. Let $z := w + \N_{k(r)/k}(w^1)$, then $\loc_\Omega(z) + \N_{k(q)/k}(z'_\Omega)$ is sufficiently close to $z_\Omega$.
\end{proof}

\begin{proof} [Proof of Proposition \ref{prop:DescentCycleMultiplicative}]
    Let $x_\Omega = (x_v)_{v \in \Omega_k} \in \Z_{0,\Omega}(X)^{\Br(X^c)}$. There is an exact sequence
	\begin{equation}
 \label{eq:DescentCycleMultiplicative1}
		1 \to M \xrightarrow{\iota} T \to Q \to 1,
	\end{equation}
    where $T$ is a $k$-torus and $Q$ is a quasi-trivial $k$-torus (see {\em e.g.} \cite[Lemma 3.8]{HW2024Supersolvable}). Let us apply Proposition \ref{prop:DescentCycleConnected} to the pushforward descent type $\iota_\ast \lambda \in \DType(X,T)$ ({\em cf.} Lemma \ref{lemm:DescentTypeAbelian}). We then obtain a finite extension $k'/k$, a torsor $g: Z \to X_{k'}$ under $T_{k'}$ with $\dtype(Z) = \iota_\ast\lambda_{k'}$, and a family $P'_\Omega \in Z(k'_\Omega)^{\Br_{\nr}(Z)}$ such that $\N^{\Omega}_{k'/k}(g_\ast [P'_\Omega])$ is sufficiently close to $x_\Omega$. Furthermore, by the exact sequence \eqref{eq:DescentTypeAbelian} (also from Lemma \ref{lemm:DescentTypeAbelian})), we have a commutative diagram
	\begin{equation} \label{eq:DescentCycleMultiplicative2}
		\xymatrix{
			\H^1(k',M) \ar[r]^{p^\ast} \ar@{->>}[d]^{\iota_\ast} & \H^1(X_{k'},M) \ar[rr]^-{\dtype} \ar[d]^{\iota_\ast} && \DType(X_{k'},M_{k'}) \ar[r]^-{\partial} \ar[d]^{\iota_\ast} & \H^2(k',M) \ar@{^{(}->}[d]^{\iota_\ast} \\
			\H^1(k',T) \ar[r]^{p^\ast} & \H^1(X_{k'},T) \ar[rr]^-{\dtype} && \DType(X_{k'},T_{k'}) \ar[r]^-{\partial} & \H^2(k',T)
		}
	\end{equation}
    with exact rows, where $p: X \to \Spec(k)$ denotes the structure morphism. Since $\H^1(k',Q) = 0$, the left arrow from \eqref{eq:DescentCycleMultiplicative2} is surjective and the right arrow there is injective, by virtue of the long exact sequence associated with \eqref{eq:DescentCycleMultiplicative1}. Since $\iota_\ast \lambda_{k'} = \dtype(Z)$, one has $\iota_\ast \partial(\lambda_{k'}) = \partial(\iota_\ast \lambda_{k'}) = 0$. Thus, $\partial(\lambda_{k'}) = 0$ by injectivity of $\iota_\ast: \H^2(k',M) \to \H^2(k',T)$. It follows that there exists a torsors $f: Y \to X_{k'}$ with $\dtype(Y) = \lambda_{k'}$. By exactness of the bottom row, one has $\iota_\ast [Y] = [Z] - p^\ast [\tau]$ for some $[\tau] \in \H^1(k',T)$. By surjectivity of $\iota_\ast: \H^1(k',M) \to \H^1(k',T)$, we may write $[\tau] = \iota_\ast [\mu]$ for some $[\mu] \in \H^1(k',M)$. Twisting $Y$ by a cocycle $\mu$ representing $[\mu]$ (which yields a torsor $Y^{\mu} \to X_{k'}$ of descent type $\dtype(Y^{\mu}) = \lambda_{k'}$), we may assume that $\iota_\ast[Y] = [Z]$. In other words, we have an isomorphism $Y \times_{k'}^{M_{k'}} T_{k'} \cong Z$ of $X_{k'}$-torsors under $T_{k'}$.  
	 
    By Nagata's theorem \cite{Nagata1963Resolution}, we may extend $g$ to a proper morphism $g^c: Z^c \to X^c_{k'}$, where $Z^c$ is a smooth compactification of $Z$. Let us now follow the argument in \cite[p. 797]{HW2024Supersolvable}. The second projection $Y \times_{k'} T_{k'} \to T_{k'}$ induces a morphism 
        \begin{equation*}
            \pi: Z \cong Y \times_{k'}^{M_{k'}} T_{k'} \to T_{k'}/M_{k'} = Q_{k'},
        \end{equation*}
    which has the following property. For every overfield $K/k'$ and every $K$-point $q \in Q(K)$, if $\delta: Q(K) \to \H^1(K,M)$ is the connecting homomorphism associated with \eqref{eq:DescentCycleMultiplicative1} and $\mu$ is a cocycle representing $[\mu] = \delta(q)$, then the fibre $Z_q = \pi^{-1}(q)$ is isomorphic to $Y^\mu_K$ as $X_K$-torsors under $M_K$ (in particular, $Z_q$ is rationally connected). Indeed, if $\varpi: T \to Q$ denotes the projection, then $[\mu] = [\varpi^{-1}(q)] \in \H^1(K,M)$ and hence $[Z_q] = [\pi^{-1}(q)] = [Y_K \times_{K}^{M_K} \varpi^{-1}(q)] = [Y_K^{\mu}] \in \H^1(X_K,M)$ by the very definition of twisting.
     
    Since $Q$ is a quasi-trivial torus, it is a dense open subset of an affine space. Hence, we may apply Lemma \ref{lemm:DescentCycleMultiplicative} to obtain a global zero-cycle $z \in \Z_0(Z^c)$, a closed point $q \in Q_{k'}$, and a family $z'_\Omega \in \Z_{0,\Omega}(Z_q)^{\Br_{\nr}(Z_q)}$ such that $\loc_\Omega(z) + \N_{k'(q)/k'}(z'_\Omega)$ is sufficiently close to $[P'_\Omega]$. Twisting $Y_{k'(q)}$ by a cocycle $\mu$ representing $\delta(q) \in \H^1(k'(q),M)$, we may assume that $\delta(q) = 0$. As mentioned above, there is an isomorphism $Z_q \cong Y_{k'(q)}$ of $X_{k'(q)}$-torsors under $M_{k'(q)}$. Thus, we obtain a family $y_\Omega \in \Z_{0,\Omega}(Y_{k'(q)})^{\Br_{\nr}(Y_{k'(q)})}$ such that $f_\ast y_\Omega = g_\ast z'_\Omega$ in $\Z_{0,\Omega}(X_{k'(q)})$. Let $x := \N_{k'/k}(g^c_\ast z) \in \Z_0(X^c)$. Then $\loc_\Omega(x) + \N_{k'(q)/k}^{\Omega}(f_\ast y_\Omega)$ is sufficiently close to $\N_{k'/k}^{\Omega}(g_\ast [P'_\Omega])$, hence also to $x_\Omega$.
\end{proof}

\subsection{The general case} \label{subsec:DescentCycleGeneral}

We prove Theorem \ref{thm:DescentCycle} in this subsection. First, we start with the finite solvable case. Recall that a finite descent type on a smooth geometrically integral variety $X$ is an isomorphism class of a (connected) finite \'etale cover $\bar{Y} \to \bar{X}$ such that $\bar{Y}$ is Galois over $X$ ({\em cf.} p. \pageref{FiniteDescentTYpe}).

\begin{lemm} \label{lemm:DescentCycleSolvable}
    Let $X$ be a smooth quasi-projective variety over a number field $k$. Let $\lambda = [\bar{Y}]$ be a rationally connected finite descent type on $X$ (in particular, $X$ is rationally connected) such that the finite group $\bar{G} = \Aut(\bar{Y}/\bar{X})$ is solvable. For any smooth compactification $X^c$ of $X$, one has
	\begin{equation*}
		\Z_{0,\Omega}(X)^{\Br_{\nr}(X)} = \ol{\loc_\Omega(\Z_0(X^c)) + \bigcup_{(K, [f: Y \to X_K]) \in \Tscr_\lambda} \N_{K/k}^{\Omega}(f_\ast (Y(K_\Omega)^{\Br_{\nr}(Y_K)}))}
	\end{equation*}
	in the set $\Z_{0,\Omega}(X)$ equipped with the topology defined in p. \pageref{TopologyZeroCycles}.
\end{lemm}
\begin{proof}
    We proceed by induction on the derived length $n$ of $\bar{G}$. If $n = 0$, there is nothing to prove. 
	
    For $n \ge 1$, let $x_\Omega \in \Z_{0,\Omega}(X)^{\Br(X^c)}$. The derived subgroup $[\bar{G},\bar{G}]$ of $\bar{G}$ is solvable of derived length $n-1$. Let $\bar{Z} := \bar{Y}/[\bar{G},\bar{G}]$, which defines a finite descent type $\lambda^{\ab}$ on $X$, and $\Aut(\bar{Z}/\bar{X}) = \bar{G}^{\ab}$ is finite abelian (hence of multiplicative type). By Proposition \ref{prop:DescentCycleMultiplicative}, we find a global zero-cycle $x \in \Z_0(X^c)$, a finite extension $k'/k$, a torsor $g: Z \to X_{k'}$ of descent type $\dtype(Z) = \lambda_{k'}^{\ab}$, and a family $P'_\Omega \in Z(k'_\Omega)^{\Br_{\nr}(Z)}$ such that $\loc_\Omega(x) + \N^{\Omega}_{k'/k}(g_\ast [P'_\Omega])$ is sufficiently close to $x_\Omega$. 

    By Nagata's theorem \cite{Nagata1963Resolution}, $g$ extends into a proper morphism $g^c: Z^c \to X^c_{k'}$, where $Z^c$ is a smooth compactification of $Z$. Now, $\bar{Y}$ is Galois over $Z$ since it is Galois over $X$, hence $\lambda_Z = [\bar{Y}]$ is a finite descent type on $Z$, with $\Aut(\bar{Y}/\bar{Z}) = [\bar{G},\bar{G}]$. By the induction hypothesis, there is a global zero-cycle $z \in \Z_0(Z^c)$, a finite extension $K/k'$, a torsor $h: Y \to Z_K$ of descent type $\res_{K/k'}(\lambda_{Z})$, and a family $P_\Omega \in \Z_{0,\Omega}(Y)^{\Br_{\nr}(Y)}$ such that $\loc_\Omega(z) + \N^{\Omega}_{K/k'}(h_\ast [P_\Omega])$ is sufficiently close to $[P'_\Omega]$. The composite $f = g_K \circ h: Y \to X_K$ is then a torsor of descent type $\dtype(Y) = \lambda_K$. We have $x + \N_{k'/k}(g_\ast^c z) \in \Z_{0,\Omega}(X^c)$ and $\loc_\Omega(x + \N_{k'/k}(g_\ast^c z)) + \N^{\Omega}_{K/k}(f_\ast [P_\Omega])$ is sufficiently close to $x_\Omega$.
\end{proof}

We can now deal with the finite case using the strategy in \cite[\S 5]{HW2020Galois}.

\begin{prop} \label{prop:DescentCycleFinite}
	Let $X$ be a smooth quasi-projective variety over a number field $k$. Let $\lambda = [\bar{Y}]$ be a rationally connected finite descent type on $X$ (in particular, $X$ is rationally connected). Then
	\begin{equation*}
		\Z_{0,\Omega}(X)^{\Br_{\nr}(X)} = \ol{\loc_\Omega(\Z_0(X^c)) + \Img(\rec_{\lambda,\Omega}^{\Br})}
	\end{equation*}
	in the set $\Z_{0,\Omega}(X)$ equipped with the topology defined in p. \pageref{TopologyZeroCycles}, for any smooth compactification $X^c$ of $X$. Here, the ``recombining'' map $\rec_{\lambda,\Omega}^{\Br}$ was defined in \eqref{eq:DescentCycleRecombinationBrauer}. More precisely, for any $x_\Omega \in \Z_{0,\Omega}(X)^{\Br(X^c)}$ and any prime number $p$, there is an integer $n_p$ prime to $p$ such that
	\begin{equation*}
		n_p x_\Omega \in \ol{\loc_\Omega(\Z_0(X^c)) + \bigcup_{(K, [f: Y \to X_K]) \in \Tscr_{\lambda}} \N_{K/k}^{\Omega}(f_\ast (Y(K_\Omega)^{\Br_{\nr}(Y)}))}.
	\end{equation*}
\end{prop}
\begin{proof}
	Let $\bar{G}:=\Aut(\bar{Y}/\bar{X})$. Let $x_\Omega \in \Z_{0,\Omega}(X)^{\Br(X^c)}$ and let $p$ be any prime number. By \cite[Lemme 5.6]{HW2020Galois}, we have a commutative diagram with exact rows
		\begin{equation*}
			\xymatrix{
				1 \ar[r] & \bar{H} \ar[r] \ar@{^{(}->}[d] & E \ar[r] \ar@{^{(}->}[d] & \Gamma_{k'} \ar[r] \ar@{=}[d] & 1 \\
				1 \ar[r] & \bar{G} \ar[r] \ar@{=}[d] & \Aut(\bar{Y}/X_{k'}) \ar[r] \ar@{^{(}->}[d] & \Gamma_{k'} \ar[r] \ar@{^{(}->}[d] & 1\\
				1 \ar[r] & \bar{G} \ar[r] & \Aut(\bar{Y}/X) \ar[r] & \Gamma_{k} \ar[r] & 1,
			}
		\end{equation*}
	where $\bar{H}$ (resp. $E$, resp. $\Gamma_{k'}$) is a $p$-Sylow subgroup of $\bar{G}$ (resp. $\Aut(\bar{Y}/X)$, resp. $\Gamma_k$) in the profinite sense (see \cite[Chapitre I, \S 1.4]{Serre1994Galois}). The subgroup $E$ is of index $n_{p}' := [\bar{G}:\bar{H}]$ in $\Aut(\bar{Y}/X_{k'})$, hence corresponds to an \'etale cover $Z' \to X_{k'}$ of degree $n_{p}'$, over which $\bar{Y}$ is Galois. The cover $Z' \to X_{k'}$ is the base change to $k'$ of an \'etale cover $Z \to X_{\ell}$, for some finite extension $\ell$ of $k$ contained in $k'$ (hence of degree $n_p'':=[\ell:k]$ prime to $p$). Recall from \eqref{eq:DescentCycleLocalRestriction} that we have a local restriction map $\res_{\ell/k}^{\Omega}: \Z_{0,\Omega}(X) \to \Z_{0,\Omega}(X_\ell)$. By functoriality of the pairing \eqref{eq:BrauerManinPairingCycle}, one has $\res_{\ell/k}^{\Omega}(x_\Omega) \in \Z_{0,\Omega}(X_\ell)^{\Br(X^c_\ell)}$. 
	
    By Nagata's theorem \cite{Nagata1963Resolution}, we may extend the $Z \to X_\ell$ to a morphism $\pi: Z^c \to X^c_\ell$, which is generically finite of degree $n_p'$, where $Z^c$ is a smooth compactification of $Z$. Applying the refined Gysin pullback $\pi^!$ \cite[Definition 8.1.2]{Fulton1984Intersection} to the rational equivalence class of $\res_{\ell/k}^{\Omega}(x_\Omega)$ yields a family $z_\Omega = (z_w)_{w \in \Omega_\ell} \in \Z_{0,\Omega}(Z^c)$ such that $\pi_\ast z_w$ is rationally equivalent to $n_p' (\res_{\ell/k}^{\Omega}(x_\Omega))_w$ for all $w \in \Omega_\ell$. Since $\res_{\ell/k}^{\Omega}(x_\Omega) \in \Z_{0,\Omega}(X_\ell)^{\Br(X^c_\ell)}$, we have $z_\Omega \in  \Z_{0,\Omega}(Z^c)^{\Br(Z^c)}$ by (covariant) functoriality of the pairing \eqref{eq:BrauerManinPairingCycle}. We may assume that $z_\Omega \in \Z_{0,\Omega}(Z)^{\Br(Z^c)}$ by continuity of \eqref{eq:BrauerManinPairingCycle} and Lemma \ref{lemm:Formal} \ref{lemm:Formal1}. Now, since $\bar{Y}$ is Galois over $Z$, it defines a finite descent type $\lambda_Z$ on $Z$, with $\Aut(\bar{Y}/\bar{Z}) = \bar{H}$. Furthermore, any $Z$-torsor of descent type $\lambda_Z$ is an $X_\ell$-torsor of descent type $\lambda_\ell$ and an $X$-torsor of descent type $\lambda$. Since $\bar{H}$ is a $p$-group ({\em a fortiori} solvable), we have
		\begin{equation*}
			z_\Omega \in \ol{\loc_\Omega(\Z_0(Z^c)) + \bigcup_{(K, [g: Y \to Z_K]) \in \Tscr_{\lambda_Z}} \N_{K/\ell}^{\Omega}(g_\ast (Y(K_\Omega)^{\Br_{\nr}(Y)}))}
		\end{equation*}
	by Lemma \ref{lemm:DescentCycleSolvable}. Applying $\pi_\ast$ to both side, we obtain 
	\begin{equation*}
		n_p' \res_{\ell/k}^{\Omega}(x_\Omega) \in \ol{\loc_\Omega(\Z_0(X^c_\ell)) + \bigcup_{(K, [f: Y \to X_K]) \in \Tscr_{\lambda_\ell}} \N_{K/\ell}^{\Omega}(f_\ast (Y(K_\Omega)^{\Br_{\nr}(Y)}))}.
	\end{equation*}
	Finally, applying the local norm $\N_{\ell/k}^{\Omega}$ to both side, one gets
	\begin{equation*}
		n_p x_\Omega \in \ol{\loc_\Omega(\Z_0(X^c)) + \bigcup_{(K, [f: Y \to X_K]) \in \Tscr_{\lambda}} \N_{K/k}^{\Omega}(f_\ast (Y(K_\Omega)^{\Br_{\nr}(Y)}))},
	\end{equation*}
	where $n_p:=n_p' n_p''$ is an integer prime to $p$.
	
	Let us now show that $x_\Omega$ lies in the closure of $\loc_\Omega(\Z_0(X^c)) + \Img(\rec_{\lambda,\Omega}^{\Br})$. Fix an integer $n \ge 1$ and a finite set $S \subseteq \Omega_k$ of places. For each prime $p | n$, we have seen that there are an integer $n_p$ prime to $p$, a global zero-cycle $x^p \in \Z_0(X^c)$, and a family $z_\Omega^p = (z_v^p)_{v \in \Omega_k} \in \Img(\rec_{\lambda,\Omega}^{\Br})$ such that $n_px_v$ and $\loc_v(x^p) + z^p_v$ have the same image in $\CH_0'(X^c_{v})/n$ for all $v \in S$. The integers $n_p$, together with $n$, are mutually coprime, hence there are integers $d,d_p$ such that $1 = dn + \sum_{p|n} d_pn_p$.	Let $x := \sum_{p|n} d_p x^p \in \Z_0(X^c)$ and $z_\Omega := (z_v)_{v \in \Omega_k} =\sum_{p|n} d_p z^p_\Omega \in \Img(\rec_{\lambda,\Omega}^{\Br})$. Then 
	\begin{equation*}
		x_v = dnx_v + \sum_{p|n} d_p n_p x_v \quad \text{and} \quad \loc_v(x) + z_v
	\end{equation*} have the same image in $\CH_0'(X^c_{v})/n$ for all $v \in S$. This proves the proposition.
\end{proof}

We can now combine Propositions \ref{prop:DescentCycleConnected} and \ref{prop:DescentCycleFinite} to obtain Theorem \ref{thm:DescentCycle} ({\em i.e.} Theorem \ref{customthm:DescentCycle}).

\begin{proof} [Proof of Theorem \ref{thm:DescentCycle}]
    Let $L^{\f} := (\bar{G}^{\f},\kappa^{\f})$ be the natural $k$-kernel (recall that $\bar{G}^{\f} = \pi_0(\bar{G}) = \bar{G}/\bar{G}^\circ$), and let $\lambda^{\f} := [\bar{Z}] \in \DType(X,L^{\f})$ be the pushforward of $\lambda$ along the  morphism $L \to L^{\f}$ induced by the projection $\bar{G} \to \bar{G}^{\f}$ (hence $\bar{Z} = \bar{Y}/\bar{G}^\circ$, {\em cf.} Lemma \ref{lemm:DescentTypePushforward} \ref{lemm:DescentTypePushforward1}). For each $x_\Omega \in \Z_0(X)^{\Br_{\nr}(X)}$, applying Proposition \ref{prop:DescentCycleFinite} to $\lambda^{\f}$ yields a global zero-cycle $x \in \Z_0(X^c)$, a finite set $\{(k_i, [g_i: Z_i \to X_{k_i}])\}_{i \in \Ical} \subseteq \Tscr_{\lambda^{\f}}$, and a family $z_\Omega^i \in \Z_0(Z_i)^{\Br_{\nr}(Z_i)}$ for each $i \in \Ical$, such that $\loc_\Omega(x) + \sum_{i \in \Ical} \N_{k_i/k}^{\Omega}((g_i)_{\ast} z_\Omega^i)$ is sufficiently close to $x_\Omega$. 
	
    For $i \in \Ical$, by Lemma \ref{lemm:DescentTypePushforward} \ref{lemm:DescentTypePushforward2}, the torsor $Z_i \to X_{k_i}$ gives rise to a $k_i$-kernel $L_i = (\bar{G}^\circ,\kappa_i)$ and a descent type $\lambda_i = [\bar{Y},E_i] \in \DType(Z_i, L_i)$. By Proposition \ref{prop:DescentCycleConnected}, there exist a finite extension $K_i/k_i$, a torsor $h_i: Y_i \to (Z_i)_{K_i}$ of descent type $\res_{K_i/k_i}(\lambda_i)$, and a family $P^i_\Omega \in Y_i((K_i)_\Omega)^{\Br_{\nr}(Y_i)}$ such that $\N_{K_i/k_i}^{\Omega} ((h_i)_\ast[P^i_\Omega])$ is sufficiently close to $z_\Omega^i$. The composite $f_i = (g_i)_{K_i} \circ h_i: Y_i \to X_{K_i}$ is a torsor of descent type $\dtype(Y_i) = \res_{K_i/k}(\lambda)$ by Lemma \ref{lemm:DescentTypePushforward} \ref{lemm:DescentTypePushforward3}, and $\loc_\Omega(x) + \sum_{i \in \Ical} \N_{K_i/k}^{\Omega}((f_i)_\ast [P^i_\Omega])$ is sufficiently close to $x_\Omega$. The theorem is proved.
\end{proof}

Theorem \ref{thm:DescentCycle} is a generalization of the following result, proved in \cite[Th\'eor\`eme A]{HW2020Galois}.

\begin{coro}
    Conjecture \ref{conj:E} holds for the smooth compactifications of homogeneous spaces of connected linear algebraic groups over number fields.
\end{coro} 
\begin{proof}
    Let $X$ a homogeneous space of a connected linear algebraic group $G$ over a number field $k$. As in the proof of Corollary \ref{coro:BorovoiThm}, we have a descent type $\lambda = [\bar{G}, E_X]$ on $X$. For any finite extension $K/k$, a torsor $Y \to X_{K}$ of type $\lambda_K$ is a left $K$-torsor under $G_K$, hence satisfies $\BMHP$ and $\BMWA$ by \cite[Corollaires 8.7 and 8.13]{Sansuc1981Arithmetique} and \cite{Chernousov1989E8}. It follows from Liang's theorems \cite[Theorems A and B]{Liang2013Arithmetic} that Conjecture \ref{conj:E} holds for their smooth compactifications. By virtue of Theorem \ref{thm:DescentCycle}, Conjecture \ref{conj:E} also holds for the smooth compactifications of $X$.
\end{proof}

{\bf Acknowledgements.} This work is a part of the author's Ph. D. thesis. I would like to express my sincerest thanks to my Ph. D. advisor, David Harari, for bringing up this problem and for his guidance. I am specially grateful to Olivier Wittenberg for his valuable comments and suggestions on the manuscript. I appreciate Cristian David Gonz\'alez-Avil\'es for the interesting discussion on the abelianization map. I also thank the Laboratoire de Math\'ematiques d’Orsay, Universit\'e Paris-Saclay and the Institut de Math\'ematiques de Jussieu-Paris Rive Gauche, CNRS for the excellent working condition. This research funded by a ``Contrat doctorant sp\'ecifique normalien'' from \'Ecole normale sup\'erieure de Paris and as well as by the postdoctoral program of Fondation Sciences Math\'ematiques de Paris.

\bibliographystyle{alpha}
\bibliography{ref}

\newcommand{\etalchar}[1]{$^{#1}$}
\begin{thebibliography}{CTSSD87b}

\bibitem[BB24]{BB2024Descent}
Francesca Balestrieri and Jennifer Berg.
\newblock {Descent and \'etale-Brauer obstructions for $0$-cycles}.
\newblock {\em {International Mathematics Research Notices}}, 2024(16):11712--11748, 2024.

\bibitem[BBM{\etalchar{+}}16]{BBMPV2016Enriques}
Francesca Balestrieri, Jennifer Berg, Michelle Manes, Jennifer Park, and Bianca Viray.
\newblock {Insufficiency of the Brauer--Manin obstruction for rational points on Enriques surfaces}.
\newblock In {\em Directions in Number Theory}, volume~3 of {\em Association for Women in Mathematics Series}, pages 1--31. Springer, Cham, 2016.

\bibitem[BGA14]{BGA2014Fundamental}
Mikhail~V. Borovoi and Cristian~D. Gonz\'alez-Avil\'es.
\newblock {The algebraic fundamental group of a reductive group scheme over an arbitrary base scheme}.
\newblock {\em Central European Journal of Mathematics}, 12(4):545--558, 2014.

\bibitem[BK00]{BK2000Brauer}
Mikhail~V. Borovoi and Boris~E. Kunyavskii.
\newblock {Formulas for the unramified Brauer group of a principal homogeneous space of a linear algebraic group}.
\newblock {\em Journal of Algebra}, 225(2):804--821, 2000.

\bibitem[Bor93]{Borovoi1993H2}
Mikhail~V. Borovoi.
\newblock {Abelianization of the second nonabelian Galois cohomology}.
\newblock {\em Duke Mathematical Journal}, 72(1):217--239, 1993.

\bibitem[Bor96]{Borovoi1996BrauerManin}
Mikhail~V. Borovoi.
\newblock {The Brauer--Manin obstructions for homogeneous spaces with connected or abelian stabilizer}.
\newblock {\em Journal f{\"u}r die reine und angewandte Mathematik}, 1996(473):181--194, 1996.

\bibitem[Bor98]{Borovoi1998Reductive}
Mikhail~V. Borovoi.
\newblock {\em {Abelian Galois Cohomology of Reductive Groups}}, volume 132 of {\em Memoirs of the American Mathematical Society}.
\newblock American Mathematical Society, 1998.

\bibitem[BvH09]{BvH2009Picard}
Mikhail~V. Borovoi and Joost van Hamel.
\newblock {Extended Picard complexes and linear algebraic groups}.
\newblock {\em Journal f\"ur die reine und angewandte Mathematik}, 2009(627):53--82, 2009.

\bibitem[CDX19]{CDX2019Compare}
Yang Cao, Cyril Demarche, and Fei Xu.
\newblock {Comparing descent obstruction and Brauer--Manin obstruction for open varieties}.
\newblock {\em Transactions of the American Mathematical Society}, 371(12):8625--8650, 2019.

\bibitem[Che54]{Chevalley1954Group}
Claude Chevalley.
\newblock {On algebraic group varieties}.
\newblock {\em Journal of the Mathematical Society of Japan}, 6(3-4):303--324, 1954.

\bibitem[Che89]{Chernousov1989E8}
Vladimir~I. Chernousov.
\newblock {On the Hasse principle for groups of type $E_8$}.
\newblock {\em Doklady Akademii Nauk SSSR}, 306(25):1059--1063, 1989.

\bibitem[Con11]{Conrad2011Reductive}
Brian Conrad.
\newblock {Reductive group schemes (SGA3 Summer school)}, 2011.
\newblock preprint, 390 pages, \href{https://math.stanford.edu/~conrad/papers/luminysga3smf.pdf}{https://math.stanford.edu/{$\sim$}conrad/papers/luminysga3smf.pdf}.

\bibitem[Con12]{Conrad2012Adelic}
Brian Conrad.
\newblock {Weil and Grothendieck approaches to adelic points}.
\newblock {\em L'Enseignement Math\'ematique}, 58(1-2):61--97, 2012.

\bibitem[CT95]{CT1995Chow}
Jean-Louis Colliot-Th\'el\`ene.
\newblock {L'arithm\'etique du groupe de Chow des z\'ero-cycles}.
\newblock {\em Journal de th\'eorie des nombres de Bordeaux}, 7(1):51--73, 1995.

\bibitem[CTC79]{CTC1979Rational}
Jean-Louis Colliot-Th\'el\`ene and Daniel Coray.
\newblock {L'\'equivalence rationnelle sur les points ferm\'es des surfaces rationnelles fibr\'ees en coniques}.
\newblock {\em Compositio Mathematica}, 39(3):301--332, 1979.

\bibitem[CTS81]{CTS1981Chow}
Jean-Louis Colliot-Th\'el\`ene and Jean-Jacques Sansuc.
\newblock {On the Chow groups of certain rational surfaces: a sequel to a paper of S. Bloch}.
\newblock {\em Duke Mathematical Journal}, 48(2):421--447, 1981.

\bibitem[CTS87]{CTS1987Descente}
Jean-Louis Colliot-Th\'el\`ene and Jean-Jacques Sansuc.
\newblock {La descente sur les vari\'et\'es rationnelles, II}.
\newblock {\em Duke Mathematical Journal}, 54(1):375–492, 1987.

\bibitem[CTS00]{CTS2000Fibration}
Jean-Louis Colliot-Th\'el\`ene and Alexei~N. Skorobogatov.
\newblock Descent on fibrations over $\mathbf{P}^1_k$ revisited.
\newblock {\em Mathematical Proceedings of the Cambridge Philosophical Society}, 128(3):383–393, 2000.

\bibitem[CTS13]{CTS2013Reduction}
Jean-Louis Colliot-Th\'el\`ene and Alexei~N. Skorobogatov.
\newblock {Good reduction of the Brauer--Manin obstruction}.
\newblock {\em Transactions of the American Mathematical Society}, 365(2):579--590, 2013.

\bibitem[CTS21]{CTS2021Brauer}
Jean-Louis Colliot-Th\'el\`ene and Alexei~N. Skorobogatov.
\newblock {\em The Brauer--Grothendieck Group}, volume~71 of {\em Ergebnisse der Mathematik und ihrer Grenzgebiete. 3. Folge / A Series of Modern Surveys in Mathematics}.
\newblock Springer International Publishing, 2021.

\bibitem[CTSSD87a]{CTSSD1987I}
Jean-Louis Colliot-Th\'el\`ene, Jean-Jacques Sansuc, and Sir~Peter Swinnerton-Dyer.
\newblock {Intersections of two quadrics and Ch\^atelet surfaces, I}.
\newblock {\em Journal f{\"u}r die reine und angewandte Mathematik}, 1987(373):37--107, 1987.

\bibitem[CTSSD87b]{CTSSD1987II}
Jean-Louis Colliot-Th\'el\`ene, Jean-Jacques Sansuc, and Sir~Peter Swinnerton-Dyer.
\newblock {Intersections of two quadrics and Ch\^atelet surfaces, II}.
\newblock {\em Journal f{\"u}r die reine und angewandte Mathematik}, 1987(374):72--168, 1987.

\bibitem[Deb01]{Debarre2001Higher}
Olivier Debarre.
\newblock {\em {Higher-Dimensional Algebraic Geometry}}.
\newblock Universitext. Springer New York, 2001.

\bibitem[Deb03]{Debarre2003Rational}
Olivier Debarre.
\newblock Vari\'et\'es rationnellement connexes.
\newblock In {\em S\'eminaire Bourbaki : volume 2001-2002, expos\'es 894-908}, number 290 in Ast\'erisque, pages 243--266. Soci\'et\'e math\'ematique de France, 2003.
\newblock Talk no. 905.

\bibitem[Dem70a]{Demazure1970Reductive1}
Michel Demazure.
\newblock {Groupes r\'eductifs - g\'en\'eralit\'es}.
\newblock In {\em Schemas en Groupes. S\'eminaire de Geometrie Algebrique du Bois-Marie 1962-1964 (SGA 3) - Tome 3 (Structure des sch\'emas en groupes r\'eductifs)}, volume 153 of {\em Lecture Notes in Mathematics}, pages 1--26. Springer Berlin Heidelberg, 1970.

\bibitem[Dem70b]{Demazure1970Reductive2}
Michel Demazure.
\newblock {Groupes r\'eductifs : d\'eploiements, sous-groupes, groupes quotients}.
\newblock In {\em Schemas en Groupes. S\'eminaire de Geometrie Algebrique du Bois-Marie 1962-1964 (SGA 3) - Tome 3 (Structure des sch\'emas en groupes r\'eductifs)}, volume 153 of {\em Lecture Notes in Mathematics}, pages 109--175. Springer Berlin Heidelberg, 1970.

\bibitem[Dem09]{Demarche2009Etale}
Cyril Demarche.
\newblock {Obstruction de descente et obstruction de Brauer--Manin \'etale}.
\newblock {\em Algebra and Number Theory}, 3(2):237--254, 2009.

\bibitem[Dem11a]{Demarche2011Approximation}
Cyril Demarche.
\newblock {Le d{\'e}faut d'approximation forte dans les groupes lin{\'e}aires connexes}.
\newblock {\em {Proceedings of the London Mathematical Society}}, 102(3):563--597, 2011.

\bibitem[Dem11b]{Demarche2011PT}
Cyril Demarche.
\newblock {Suites de Poitou--Tate pour les complexes de tores \`a deux termes}.
\newblock {\em {International Mathematics Research Notices}}, 2011(1):135--174, 2011.

\bibitem[DLA19]{DLA2019Reduction}
Cyril Demarche and Giancarlo Lucchini~Arteche.
\newblock {Le principe de Hasse pour les espaces homog\`enes : r\'eduction au cas des stabilisateurs finis}.
\newblock {\em Compositio Mathematica}, 158(8):1568--1593, 2019.

\bibitem[FSS98]{FSS1998H2}
Yuval Flicker, Claus Scheiderer, and Ramdorai Sujatha.
\newblock {Grothendieck's theorem on non-abelian $H^2$ and local-global principles}.
\newblock {\em Journal of the American Mathematical Society}, 11(3):731–750, 1998.

\bibitem[Ful84]{Fulton1984Intersection}
William Fulton.
\newblock {\em {Intersection Theory}}, volume~2 of {\em {Ergebnisse der Mathematik und ihrer Grenzgebiete: a series of modern surveys in mathematics. Folge 3}}.
\newblock Springer New York, NY, 2 edition, 1984.

\bibitem[GA12]{GA2012Abelianization}
Cristian~D. Gonz\'alez-Avil\'es.
\newblock Quasi-abelian crossed modules and nonabelian cohomology.
\newblock {\em Journal of Algebra}, 369(1):235--255, 2012.

\bibitem[GA18]{GA2018UPic}
Cristian~D. Gonz\'alez-Avil\'es.
\newblock {The units-Picard complex and the Brauer group of a product}.
\newblock {\em Journal of Pure and Applied Algebra}, 222(9):2746--2772, 2018.

\bibitem[GA19]{GA2019Reductive}
Cristian~D. Gonz\'alez-Avil\'es.
\newblock {The units-Picard complex of a reductive group scheme}.
\newblock {\em Journal of the Ramanujan Mathematical Society}, 34(2):191--230, 2019.

\bibitem[GHS03]{GHS2003Rational}
Tom Graber, Joe Harris, and Jason Starr.
\newblock Families of rationally connected varieties.
\newblock {\em Journal of the American Mathematical Society}, 16(1):57--67, 2003.

\bibitem[Gir71]{Giraud1971Cohomology}
Jean Giraud.
\newblock {\em {Cohomologie non ab\'elienne}}, volume 179 of {\em Grundlehren der mathematischen Wissenschaften}.
\newblock Springer Berlin Heidelberg, 1971.

\bibitem[Gro67]{EGAIV4}
Alexander Grothendieck.
\newblock {\'El\'ements de g\'eom\'etrie alg\'ebrique : {IV.} {\'Etude} locale des sch\'emas et des morphismes de sch\'emas, {Quatri\`eme} partie}.
\newblock {\em Publications Math\'ematiques de l'IH\'ES}, 32:5--361, 1967.

\bibitem[Gro70a]{Grothendieck1970Multiplicatif1}
Alexandre Grothendieck.
\newblock {Caract\'erisation et classification des groupes de type multiplicatif}.
\newblock In {\em Schemas en Groupes. S\'eminaire de Geometrie Algebrique du Bois-Marie 1962-1964 (SGA 3) - Tome 2 (Groupes de type nultiplicatif, et structure des schemas en groupes generaux)}, volume 152 of {\em Lecture Notes in Mathematics}, pages 63--106. Springer Berlin Heidelberg, 1970.

\bibitem[Gro70b]{Grothendieck1970Multiplicatif2}
Alexandre Grothendieck.
\newblock {Groupes de type multiplicatif : homomorphismes dans un sch\'ema en groupes}.
\newblock In {\em Schemas en Groupes. S\'eminaire de Geometrie Algebrique du Bois-Marie 1962-1964 (SGA 3) - Tome 2 (Groupes de type nultiplicatif, et structure des schemas en groupes generaux)}, volume 152 of {\em Lecture Notes in Mathematics}, pages 31--62. Springer Berlin Heidelberg, 1970.

\bibitem[Gro70c]{Grothendieck1970Multiplicatif3}
Alexandre Grothendieck.
\newblock {Groupes diagonalisables}.
\newblock In {\em Schemas en Groupes. S\'eminaire de Geometrie Algebrique du Bois-Marie 1962-1964 (SGA 3) - Tome 2 (Groupes de type nultiplicatif, et structure des schemas en groupes generaux)}, volume 152 of {\em Lecture Notes in Mathematics}, pages 1--30. Springer Berlin Heidelberg, 1970.

\bibitem[Gro73]{Grothendieck1973Foncteurs}
Alexandre Grothendieck.
\newblock {Foncteurs fibres, supports, \'etude cohomologique des morphismes finis}.
\newblock In {\em Th\'eorie des Topos et Cohomologie Etale des Sch\'emas. S\'eminaire de G\'eom\'etrie Algébrique du Bois-Marie 1963-1964 (SGA 4) - Tome 3}, volume 305 of {\em Lecture Notes in Mathematics}, pages 366--412. Springer Berlin Heidelberg, 1973.

\bibitem[Har94]{Harari1994Fibration}
David Harari.
\newblock {M\'ethode des fibrations et obstruction de Manin}.
\newblock {\em Duke Mathematical Journal}, 75(1):221--260, 1994.

\bibitem[Hir64]{Hironaka1964Resolution}
Heisuke Hironaka.
\newblock {Resolution of singularities of an algebraic variety over a field of characteristic zero: I}.
\newblock {\em Annals of Mathematics}, 79(1):109–203, 1964.

\bibitem[HS02]{HS2002Nonabelian}
David Harari and Alexei~N. Skorobogatov.
\newblock {Non-abelian cohomology and rational points}.
\newblock {\em Compositio Mathematica}, 130(3):241--273, 2002.

\bibitem[HS05]{HS2005Enriques}
David Harari and Alexei~N. Skorobogatov.
\newblock {Non-abelian descent and the arithmetic of Enriques surfaces}.
\newblock {\em International Mathematics Research Notices}, 2005(52):3203--3228, 2005.

\bibitem[HS13]{HS2013Descent}
David Harari and Alexei~N. Skorobogatov.
\newblock {Descent theory for open varieties}.
\newblock In {\em Torsors, \'Etale Homotopy and Applications to Rational Points}, volume 405 of {\em London Mathematical Society Lecture Note Series}, pages 250--279. Cambridge University Press, 2013.

\bibitem[HW16]{HW2016Fibration}
Yonatan Harpaz and Olivier Wittenberg.
\newblock {On the fibration method for zero-cycles and rational points}.
\newblock {\em {Annals of Mathematics}}, 183(1):229--295, 2016.

\bibitem[HW20]{HW2020Galois}
Yonatan Harpaz and Olivier Wittenberg.
\newblock {Z{\'e}ro-cycles sur les espaces homog{\`e}nes et probl{\`e}me de Galois inverse}.
\newblock {\em {Journal of the American Mathematical Society}}, 33(3):775--805, 2020.

\bibitem[HW24]{HW2024Supersolvable}
Yonatan Harpaz and Olivier Wittenberg.
\newblock {Supersolvable descent for rational points}.
\newblock {\em Algebra \& Number Theory}, 18(4):787--814, 2024.

\bibitem[Kot86]{Kottwitz1986Trace}
Robert~E. Kottwitz.
\newblock {Stable trace formula: Elliptic singular terms}.
\newblock {\em Mathematische Annalen}, 275:365--400, 1986.

\bibitem[KS86]{KS1986Arithmetic}
Kazuya Kato and Shuji Saito.
\newblock {Global class field theory of arithmetic schemes}.
\newblock In {\em Applications of Algebraic K-theory to Algebraic Geometry and Number Theory, Part I, II}, volume~55 of {\em Contemporary Mathematics}, pages 255--331. American Mathematical Society, 1986.

\bibitem[Lia13]{Liang2013Arithmetic}
Yongqi Liang.
\newblock {Arithmetic of $0$-cycles on varieties defined over number fields}.
\newblock {\em Annales scientifiques de l'\'Ecole Normale Sup\'erieure}, 4e s{\'e}rie, 46(1):35--56, 2013.

\bibitem[Lin25]{Linh2025Type}
Nguyen~M. Linh.
\newblock {Nonabelian descent types}.
\newblock {\em Journal of Algebra}, 661:380--429, 2025.

\bibitem[Man71]{Manin1971Brauer}
Yuri~I. Manin.
\newblock {Le groupe de Brauer--Grothendieck en g\'eom\'etrie diophantienne}.
\newblock In {\em Actes du Congr{\`e}s international des Math{\'e}maticiens (Nice, 1970)}, pages 401--411. Gauthier-Villars, Paris, 1971.

\bibitem[Mil80]{Milne1980Etale}
James~S. Milne.
\newblock {\em Etale Cohomology}, volume~33 of {\em Princeton Mathematical Series}.
\newblock Princeton University Press, 1980.

\bibitem[Nag63]{Nagata1963Resolution}
Masayoshi Nagata.
\newblock {A generalization of the imbedding problem of an abstract variety in a complete variety}.
\newblock {\em Journal of Mathematics of Kyoto University}, 3(1):89--102, 1963.

\bibitem[Nis55]{Nishimura1955Rational}
Hajime Nishimura.
\newblock {Some remark on rational points}.
\newblock {\em Memoirs of the College of Science, University Kyoto, Series A. Mathematics}, 29(2):189--192, 1955.

\bibitem[NSW08]{NSW2008Cohomology}
J\"urgen Neukirch, Alexander Schmidt, and Kay Wingberg.
\newblock {\em Cohomology of Number Fields}, volume 323 of {\em Grundlehren der mathematischen Wissenschaften}.
\newblock Springer Berlin Heidelberg, 2008.

\bibitem[Poo10]{Poonen2010Points}
Bjorn Poonen.
\newblock {Insufficiency of the Brauer--Manin obstruction applied to \'etale covers}.
\newblock {\em Annals of Mathematics}, 171(3):2157--2169, 2010.

\bibitem[PR94]{PR1994Group}
Vladimir~P. Platonov and Andrei~S. Rapinchuk.
\newblock {\em Algebraic Groups and Number Theory}, volume 139 of {\em Pure and Applied Mathematics}.
\newblock AcademicPress, Inc., Boston, MA, 1994.

\bibitem[San81]{Sansuc1981Arithmetique}
Jean-Jacques Sansuc.
\newblock {Groupe de Brauer et arithm\'etique des groupes alg\'ebriques lin\'eaires sur un corps de nombres}.
\newblock {\em Journal f\"ur die reine und angewandte Mathematik}, 1981(327):12--80, 1981.

\bibitem[Ser94]{Serre1994Galois}
Jean-Pierre Serre.
\newblock {\em Cohomologie Galoisienne: Cinqui\`eme \'edition, r\'evis\'ee et compl\'et\'ee}, volume~5 of {\em Lecture Notes in Mathematics}.
\newblock Springer-Verlag Berlin Heidelberg, 1994.

\bibitem[Sko01]{Skorobogatov2001Torsors}
Alexei~N. Skorobogatov.
\newblock {\em Torsors and Rational Points}, volume 144 of {\em Cambridge Tracts in Mathematics}.
\newblock Cambridge University Press, 2001.

\bibitem[Sko09]{Skorobogatov2009Descent}
Alexei~N. Skorobogatov.
\newblock {Descent obstruction is equivalent to {\'e}tale Brauer–-Manin obstruction}.
\newblock {\em {Mathematische Annalen}}, 344:501--510, 2009.

\bibitem[Sto06]{Stoll2006Descent}
Michael Stoll.
\newblock {Finite descent obstructions and rational points on curves}.
\newblock {\em Algebra and Number Theory}, 1(4):349--391, 2006.

\bibitem[vH03]{vanHamel2003Cycle}
Joost van Hamel.
\newblock {The Brauer--Manin obstruction for zero-cycles on Severi–-Brauer fibrations over curves}.
\newblock {\em Journal of the London Mathematical Society}, 68(2):317--337, 2003.

\bibitem[Wei94]{Weibel1994Homological}
Charles~A. Weibel.
\newblock {\em {An Introduction to Homological Algebra}}, volume~38 of {\em Cambridge Studies in Advanced Mathematics}.
\newblock Cambridge University Press, 1994.

\bibitem[Wei16]{Wei2016Open}
Dasheng Wei.
\newblock {Open descent and strong approximation}, 2016.
\newblock preprint, 13 pages, \url{https://arxiv.org/abs/1604.00610}.

\bibitem[Wit12]{Wittenberg2012Fibration}
Olivier Wittenberg.
\newblock {Z\'ero-cycles sur les fibrations au-dessus d'une courbe de genre quelconque}.
\newblock {\em Duke Mathematical Journal}, 161(11):2113--2166, 2012.

\bibitem[Wit24]{Wittenberg2024ParkCity}
Olivier Wittenberg.
\newblock {Park City lecture notes: around the inverse Galois problem}, 2024.
\newblock preprint, 32 pages, \url{https://arxiv.org/abs/2302.13719}, to appear in IAS/Park City Mathematics Series, AMS.

\end{thebibliography}
\end{document}